\documentstyle[leqno,11pt,epsf]{report}
\font\tenmsbm=msbm10  \font\sevenmsbm=msbm7
\font\fivemsbm=msbm5
\newfam\msbmfam
\textfont\msbmfam=\tenmsbm
\scriptfont\msbmfam=\sevenmsbm
\scriptscriptfont\msbmfam=\fivemsbm
\def\db#1{{\fam\msbmfam\relax#1}}
\def\db#1{{\fam\msbmfam\relax#1}}

\newtheorem{theorem}{Theorem}
\newtheorem{proposition}[theorem]{Proposition}
\newtheorem{lemma}[theorem]{Lemma}
\newtheorem{definition}[theorem]{Definition}
\newtheorem{corollary}[theorem]{Corollary}
\newtheorem{remark}[theorem]{Remark}
\newtheorem{ansatz}[theorem]{Ansatz}

\newcounter{lfdnum}[chapter]
\newcommand{\lfd}{\stepcounter{lfdnum}}

\begin{document}

\newcommand{\D}{\displaystyle}
\newcommand{\SC}{\scriptstyle}
\newcommand{\T}{\textstyle}
\newcommand{\be}{\lfd \begin{equation}}
\newcommand{\ee}{\end{equation}}
\newcommand{\bea}{\begin{eqnarray*}}
\newcommand{\eea}{\end{eqnarray*}}
\newcommand{\beqa}{\begin{eqnarray}}
\newcommand{\eeqa}{\end{eqnarray}}
\newcommand{\mref}[1]{(\ref{#1})}
\newcommand{\bd}{\begin{displaymath}}
\newcommand{\ed}{\end{displaymath}}
\newcommand{\bef}{\lfd \begin{figure}[htb]}
\newcommand{\eef}{\end{figure}}

\newcommand{\summa}{\sum_{k=1}^{\infty}}
\newcommand{\summb}{\sum_{k=2}^{\infty}}
\newcommand{\rhof}{\rho_{F,c}}

\newcommand{\raua}{{\cal L} _{0,w}}
\newcommand{\raub}{{\cal L} _{\sigma ,w}}
\newcommand{\raubo}{{\cal L} _{\sigma _0,w}}
\newcommand{\rauc}{\ell _{1,w}}
\newcommand{\nora}{\| _{0,w}}
\newcommand{\norb}{\| _{\sigma ,w}}
\newcommand{\norbo}{\| _{\sigma _0,w}}
\newcommand{\norc}{\| _{\ell _{1,w}}}

\newcommand{\tue}{T _{u,\epsilon}}
\newcommand{\bue}{B _{u,\epsilon}}
\newcommand{\vue}{\tilde{v}}

\newcommand{\ganz}{{\db Z}}
\newcommand{\ganzz}{{\scriptsize {\db Z}}}
\newcommand{\natz}{{\db N}}
\newcommand{\relz}{{\db R}}
\newcommand{\comz}{{\db C}}

\newcounter{pictyp}

\newenvironment{proof}{\medskip {\bf Proof : }}{\begin{flushright}
                       $\Diamond$ \end{flushright} \medskip}

\title{Forced Lattice Vibrations -- a videotext.}

\author{Percy Deift, Courant Institute\\
  \and 
Thomas Kriecherbauer, University of Augsburg\\
  \and 
Stephanos Venakides, Duke University}

\maketitle

\chapter{Introduction}
\label{e}

We begin with a description of recent numerical and analytical
results that are closely related to the results of this paper.

In 1978 Holian and Straub \cite{HS} conducted an extensive series of
numerical experiments on a driven, semi-infinite lattice
\be
\ddot{x}_n=F(x_{n-1}-x_n)-F(x_n-x_{n+1}),\qquad n=1,2,\ldots,
\label{e.5}
\ee
with initial conditions
\be
x_n(0)=nd,\quad\dot{x}_n(0)=0,\quad
n=1,2,\ldots,\,\,\,d\,\,\,{\rm constant},
\label{e.10}
\ee
for a variety of force laws $F$, and in the case that the velocity
of the driving particle $x_0$ is fixed,
\be
x_0(t)=2at,\qquad t\geq0,\quad a>0.
\label{e.15}
\ee
They discovered, in particular, a striking new phenomenon -- the
existence of a critical ``shock'' strength $a_{{\rm crit}}$.  If
$a<a_{{\rm crit}}=a_{{\rm crit}}(F)$, then in the frame moving with the
particle $x_0$, they observed behavior similar to that shown in 
Figure \ref{fe8}.
\lfd
\begin{figure}[hpbt]
\leavevmode \epsfysize=6.0cm
\epsfbox{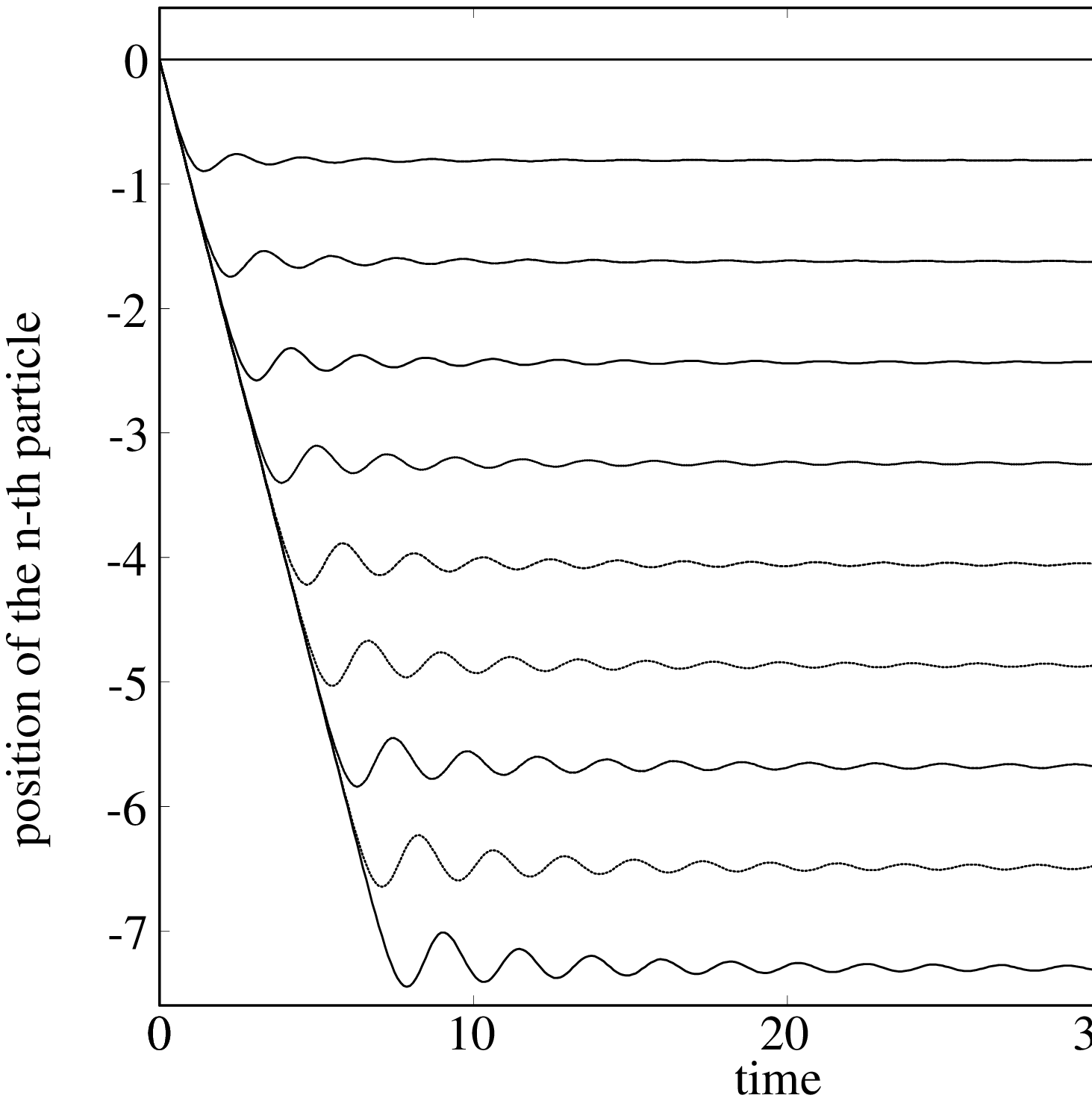}
\unitlength=1cm
\caption{{\em Motion of the first ten particles of a lattice described  
by the above system (1.1) -- (1.3) with $F(x) = e^x, d = 0, a=.5$, 
in the frame of $x_0$ (case $a < a_{{\rm crit}}$).}}
\label{fe8}
\eef
Thus the particles come to rest in a regular lattice behind the
driver.  However, if $a>a_{{\rm crit}}=a_{{\rm crit}}(F)$, then, again in the
frame of the driver, they observed behavior as in Figure \ref{fe9}.
\lfd
\begin{figure}[hpbt]
\leavevmode \epsfysize=6.0cm
\epsfbox{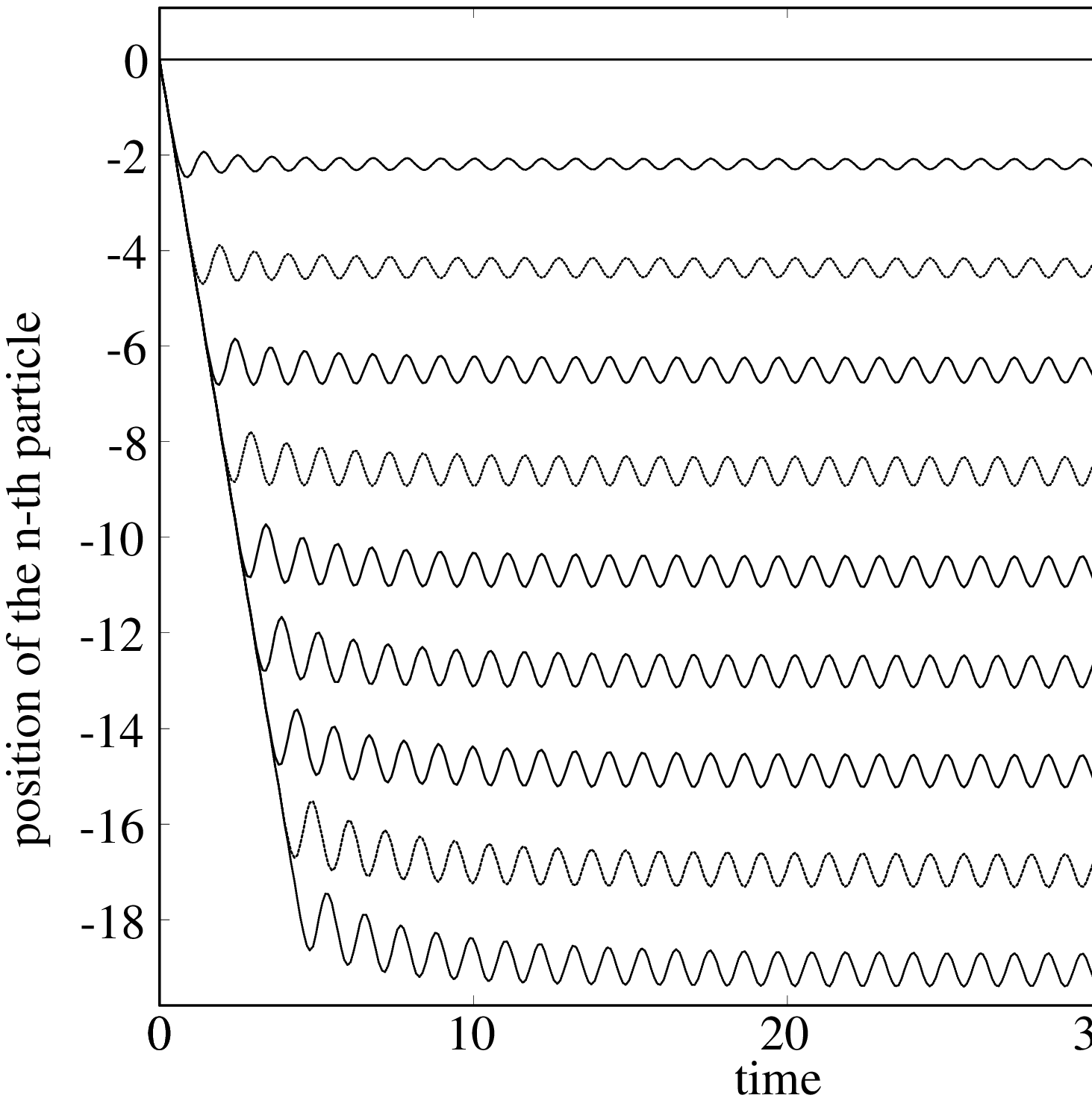}
\unitlength=1cm
\caption{{\em Motion of the first ten particles of a lattice described  
by the above system (1.1) -- (1.3) with $F(x) = e^x, d = 0, a=2$, 
in the frame of $x_0$ (case $a > a_{{\rm crit}}$).}}
\label{fe9}
\eef 
Now the particles do {\bf not} come to rest behind the driver, but
execute an on-going binary oscillation (i.e.
$x_n(t+T)=x_n(t),\,x_n(t)=x_{n+2}(t)+{\rm const.}$).
This is a marvelous, fundamentally nonlinear phenomenon; if
$F$ is linear, the effect is absent. 
\par
This phenomenon has now been observed for many different singular
and nonsingular, nonlinear force laws $F$, but an explanation of
the phenomenon from first principles in the general case has not
yet been given.  We believe that the phenomenon is present for 
a very wide class of
genuinely nonlinear $F$ (in particular, if $F'' > 0$) with $F'>0$. 
Observe that if $F'>0$, then the
force on the $n^{\rm th}$ particle $\ddot x_n=F(x_{n-1}-x_n)-F(x_n-
x_{n+1})$ is negative for $x_n>\frac{1}{2}\,(x_{n-1}+x_{n+1})$ and positive
for $x_n<\,\frac{1}{2}\,(x_{n-1}+x_{n+1})$.  Thus the only equilibrium
configuration is the regular lattice $x_n=c n, c \in \relz$, and moreover, in
this case, all the forces are restoring.
\par
In 1981, Holian, Flaschka and McLaughlin [HFM] considered the shock
problem in the special case in which $F$ is an exponential
$F(x)=e^x$, the so-called Toda shock problem.  They considered
this case because the Toda equation
\be
\ddot{x}_n=e^{x_{n-1}-x_n}-e^{x_n-x_{n+1}},
\label{e.30}
\ee
with appropriate boundary conditions, is well-known to be
completely integrable (a fact discovered by Flaschka \cite{Fla} and Manakov
\cite{Man};
see also \cite{H}) and hence there was the possibility of solving
(\ref{e.5}) -- (\ref{e.15}) explicitly and so explaining the phenomena observed
by Holian and Straub in the special case where $F(x) = e^x$. 
However, the driven
system (\ref{e.5}) -- (\ref{e.15}) is non-autonomous and it was not clear a
priori that the (formal) integrability of the Toda equation could
be used to analyze the system.  For example, the one-dimensional
oscillator $\ddot{x}+ax+bx^3=0$ is certainly integrable; however,
the driven oscillator $\ddot{x}+ax+bx^3=f(t)$, the so-called Duffing
system, is far from ``integrable'' and requires highly
sophisticated techniques for its analysis.  Nevertheless, 
Holian, Flaschka
and McLaughlin (\cite{HFM}) realized that if they went into the frame of
the driver, so that (\ref{e.10}), (\ref{e.15}) become
\be
x_n(0)=nd,\quad\dot{x}_n(0)=-2a,\quad n\geq 1,
\label{e.35}
\ee
\be
x_0(t)\equiv0,
\label{e.40}
\ee
and doubled up the system
\be
x_n(t)\equiv-x_{-n}(t),\quad n<0,
\label{e.45}
\ee
then the full system $\{x_n\}^\infty_{n=-\infty}$ solves the
{\bf autonomous} Toda equations
\be
\ddot{x}_n=e^{x_{n-1}-x_n}-e^{x_n-x_{n+1}},\quad-\infty<n<\infty,
\label{e.50}
\ee
with initial conditions
\be
x_n(0)=dn,\quad\dot{x}_n(0)=-2a(\mbox{ sgn } n),\quad-\infty<n<\infty.
\label{e.55}
\ee
But the solutions of these equations lie in a class to which the
method of inverse scattering applies.  To see what is involved we
use Flaschka's variables, 
\be
a_n=-\dot{x}_n/2,\quad
b_n=\frac{1}{2}\,\,e^{\frac{1}{2}(x_n-x_{n+1})},\,\,\,-\infty<n<\infty,
\label{e.60}
\ee
and arrange these variables into a doubly-infinite tridiagonal
matrix
\be
\widetilde{L}=
\left(
\begin{array}{lllllll}
\ddots&\ddots&\ddots\\
&b_{-2}&a_{-1}&b_{-1}&&\bigcirc\\
&&b_{-1}&a_0&b_0\\
&&&b_0&a_1&b_1\\
&\bigcirc&&&b_1&\ddots&\ddots\\
&&&&&\ddots
\end{array}
\right),
\label{e.65}
\ee
which represents the state of the system at any given time,
with companion matrix
\be
\widetilde{B}=
\left(
\begin{array}{lllllll}
\,\,\quad\ddots&\,\,\quad\ddots&\,\,\quad\ddots\\
&-b_{-2}&\quad0&\quad b_{-1}&&\bigcirc\\
&&-b_{-1}&\quad0&\quad b_0\\
&&&-b_0&\quad0&\quad b_1\\
&\bigcirc&&&-b_1&\quad0&\quad\ddots\\
&&&&&\ddots&\ddots
\end{array}
\right).
\label{e.70}
\ee
\medskip
Then, remarkably, (\ref{e.50}) is equivalent to the so-called Lax-pair
system
\be
\frac{d\tilde L}{dt}=[\widetilde{B},\widetilde{L}]=\widetilde{B}\widetilde
{L}-\widetilde{L}\widetilde{B}.
\label{e.75}
\ee
Thus the Toda equations are equivalent to an iso-spectral
deformation of the matrix operator $\widetilde L$. Inverse
scattering theory tells us that one can solve (\ref{e.75}), and hence
(\ref{e.5}) -- (\ref{e.15}), through the scattering map for $\widetilde{L}$.
Rescaling time, one sees that it is sufficient to consider the case
where the initial spacing $d=0$.  Then at $t=0$,
\be
a_n=a \mbox{ sgn } n,\quad b_n=\frac{1}{2},
\label{e.80}
\ee
and one sees that the essential spectrum of $\widetilde{L}$ is given
by two bands (cf Figure \ref{fe1}).
\bef
\leavevmode \epsfysize=0.6cm
\epsfbox{fe1.ps}
\caption{{\em The spectrum of $\tilde{L}(0)$}}
\label{fe1}
\eef
The bands overlap if and only if $a<1$.  Holian, Flaschka
and McLaughlin observed that supercritical behavior occurred for the
Toda shock problem only if the gap was open.  Hence they
identified $a_{{\rm crit}}(F=e^x)=1$.  Using the inverse method they
were able to calculate a number of other features of the Toda
shock problem, such as the speed and the form of the shock front,
and also the form of the binary oscillations.  The problem of how to 
extract detailed information about the long-time behavior of
the Toda shock problem from knowledge of the initial data using
the rather formidable formulae of inverse scattering theory,
however, remained open.
\par
In the early 80's, a very important development took place in the
analysis of infinite-dimensional integrable systems in the form of
the calculation by Lax and Levermore (\cite{LL}) of the leading order
asymptotics for the zero-dispersion limit of the Kortweg de Vries
equation, in which the weak limit of the solution as the
dispersion coefficient tends to zero is derived and the small
scale oscillations that arise are averaged out.   
This was followed in the late 80's by the calculation
of Venakides \cite{V} for the higher order terms in the Lax-Levermore
theory which produces the detailed structure of the small scale
oscillations.   These developments raised the possibility of being able
to analyze the inverse scattering formulae for the solution of the
Toda shock problem effectively as $t\rightarrow\infty$, and in
\cite{VDO}, Venakides, Deift and Oba proved the following result in the
supercritial case $a>1$:
\par
In addition to the shock speed $v_s$ calculated by Holian, Flaschka and
McLaughlin, there is a second speed $0<v_0<v_s$.
\par
In the frame moving with the driver, as $t\rightarrow\infty$,\hfill
\break
\label{pg1}
\begin{itemize}
\item for $0<n/t<v_0$, the lattice converges
to a binary oscillation $x_n(t+T)=x_n(t),\,\,\,x_{n+2}(t)=x_n(t)+$
constant, (cf Figure \ref{fe9}).  The band structure corresponding to
the asymptotic solution is $[-a-1,\,-a+1]\,\cup\,[a-1,\,a+1]$.
The binary oscillation is connected to the driver $x_0(t)\equiv0$,
through a boundary layer, in which the local disturbance due to the
driver decays exponentially in $n$.  
\medskip
\item for $v_0<n/t<v_s$,\quad the asymptotic motion
is a modulated, single-phase, quasi-periodic Toda wave with band
structure $[-a-1,\,\gamma(n/t)]\,\cup\,[a-1,\,a+1]$, where $\gamma
(n/t)$ varies monotonically from $-a+1$ to $-a-1$ as $n/t$
increases from $v_0$ to $v_s$.
\medskip
\item for $n/t>v_s$, the deviation of the particles from
their initial motion $-2at$ is exponentially small.  The influence
from the shock has not yet been felt.  As noted in \cite{HFM}, for
$n/t\sim v_s$, the motion of the lattice is described by a Toda
solution with associated spectrum $\{-a-1\}\,\cup\,[a-1,\,a+1]$.
\end{itemize}
\par
In 1991, again using the techniques in \cite{LL} and \cite{V}, Kamvissis 
(\cite{Kam})
showed that in the subcritical case $a<1$, in the frame moving
with the driver, as $t\rightarrow\infty$, the oscillatory motion
behind the shock front dies down to a quiescent regular lattice
with spacing $x_{n+1}-x_n\rightarrow-2\log(1+a)$, (cf Figure \ref{fe8}).
\vspace{.2in}
\newline
{\bf A ``Thermodynamic'' Remark.}

It is easy to see that the average spacing of the binary
oscillation of the asymptotic state in the case $a > 1$ is 
given by $- \ln 4a$. Thus the average spacing of the asymptotic 
states is given by
\begin{eqnarray*}
-2 \ln (1+a), &\mbox{ for }& a < 1, \\
- \ln 4a, &\mbox{ for }& a > 1.
\end{eqnarray*}
Observe that these expressions and their first derivatives
agree at $a=1$. Thus
we may say that the density of the asymptotic state has a  
\underline{second order} phase transition at $a = 1$.
\par
As in \cite{LL} and \cite{V}, the above results are not fully rigorous and
rely on certain (reasonable) asymptions that have not yet been
justified from first principles.  In particular, as in \cite{LL},  
the contribution
of the reflection coefficient associated with the band
$[a-1,\,a+1]$ is ignored.  Also, as in \cite{V}, an ansatz is needed to
control the long-time behavior of certain integrals.  
Recently in \cite{DMV}, the authors, again
using the approach of \cite{LL} and \cite{V}, circumvented the first
difficulty by considering finite dimensional approximations to the
lattice of length $\ell(t)\gg t$, but they still need the above
mentioned ansatz in order to re-derive the results
in \cite{VDO}.  In \cite{BK}, Bloch and Kodama consider the Toda shock
problem, both in the subcritical and the supercritical cases, from
the point of view of Whitham modulation theory in which the
validity of a modulated wave form for the solution is assumed a
priori, and the parameters of the modulated wave form are
calculated explicitly.
More recently in \cite{GN}, Greenberg and Nachman have considered the shock
problem for a general force law in the weak shock limit; they are able to 
describe many aspects of the solution, including the modulated wave 
region where they use a KdV-type continuum limit.
\par
In a slightly different direction, motivated by the so-called 
von-Neumann problem arising in the computation of shock fronts using
discrete approximations, Goodman and Lax \cite{GL} and Hou and Lax 
\cite{HL} observed and analyzed features strikingly similar to those in \cite{HS}.
Finally, 
in an interesting series of papers spanning the 80's,
Kaup and his collaborators, \cite{Kp}, \cite{KN}, \cite{WK}, use various
integrable features of the non-autonomous system 
(\ref{e.5}) -- (\ref{e.15}) to
gain valuable insight into the Toda shock problem. We will return
to these papers below.
\par
In this paper we consider the driven lattice (\ref{e.5}), (\ref{e.10}) in the
case where the uniform motion of the driving particle $x_0$ is
periodically perturbed\footnote{The more general initial value
$x_n(0)=dn,\,\dot{x}_n(0)=0$, can clearly be converted to (\ref{e.90})
by shifting the argument of $F,\,F(\cdot)\rightarrow F(\cdot-d):$
in the case of Toda, as noted above, this shift converts into a
rescaling of the time.}  
\be
\left\{\begin{array}{ll}
\ddot{x}_n=F(x_{n-1}-x_n)-F(x_n-x_{n+1}),&n\geq1,\\
x_n(0)=\dot{x}_n(0)=0,&n\geq1,\\
x_0(t)=2at+h(\gamma t).
\end{array}
\right.
\label{e.90}
\ee
Here $h(\cdot)$ is periodic with period $2\pi$ and the frequency
$\gamma>0$ is constant.  We restrict our attention to the case
where the average value of the velocity of the driver
$\overline{\dot{x}_0}=2a$ is subcritical, i.e. $a<a_{{\rm crit}}$.  (For some
discussion of the supercritical case $a>1$, see Problem 3 at the end of
the Introduction below).  
Again we consider a variety of force laws $F$, but
henceforth we restrict our attention to forces which are real
analytic and monotone increasing in the region of interest.
\par
Typically we observed the following phenomena:
\par
In the frame moving with the average velocity $2a$ of the driver,
as $t\rightarrow\infty$, the asymptotic motion of the particles
behind the shock front, is $\frac{2\pi}{\gamma}$-periodic in time,
\be
x_n(t+\frac{2\pi}{\gamma})=x_n(t),\quad0<n\ll t.
\label{e.95}
\ee
Moreover, there is a sequence of thresholds,\newline
\parbox{.45in}{\be \label{e.100} \ee} 
\parbox{5.05in}{
\bea
\begin{array}{ll}
\gamma_1=\gamma_1(a,h,F)>\gamma_2=\gamma_2(a,h,F)>\cdots>&\gamma_k=\gamma
_k(a,h,F)>\cdots>0,\\
&\gamma_k\rightarrow0\quad{\rm as}\quad k\rightarrow\infty. 
\end{array}
\eea
}
\begin{itemize}
\item If $\gamma>\gamma_1$, there exist constants $c,d$
such that $x_n-cn-d$ converges exponentially to zero as
$n\rightarrow\infty$, (cf Figure \ref{fc1}).  In other words, the effect
of the oscillatory component of the driver does not propagate into
the lattice and away from the boundary $n=0$.  The lattice behaves
in a similar way to the subcritical case of constant driving
considered by Holian, Flaschka and McLaughlin.
\medskip
\item If $\gamma_1>\gamma>\gamma_2$, then the asymptotic
motion is described by a travelling wave
\be
x_n(t)=c_1n+X_1(\beta_1n+\gamma t),\qquad1\ll n\ll t,
\label{e.105}
\ee
transporting energy away from the driver $x_0$.  (See Figure \ref{fc2}).
Here $c_1=c_1(a,h,F,\gamma)$ and
$X_1(\cdot)=X_1(\cdot\,;a,h,F,\gamma)$ is a $2\pi$-periodic
function.
\medskip
\item More generally, if $\gamma_k>\gamma>\gamma_{k+1}$,
a multi-phase wave emerges which is well-described by the
wave form
\newline
\parbox{.45in}{\be \label{e.110} \ee}
\parbox{4.9in}{
\bd
x_n=c_kn+X_k(\beta_1n+\gamma t,\,\beta_2n+2\gamma
t,\ldots,\beta_kn+k\gamma t),\quad1\ll n\ll t,
\ed
}
again transporting energy away from the driver (see Figure \ref{fc3} for the
case $k=2$).
Here $c_k=c_k(a,h,F,\gamma)$ and
$X_k(\cdot,\ldots,\cdot ;a,h,F,\gamma)$ is $2\pi$-periodic in
each of its $k$ variables.  
\par
Thus, at the phenomenological level,
we see that the periodically driven lattice behaves like a long,
heavy rope which one shakes up and down at one end.
\end{itemize}
\vspace{.2in}
{\bf Remark:}
\newline
We have restricted our experiments to the asymptotic region $1 \ll n
\ll t$. However, we expect that the solution also exhibits many
interesting phenomena when studied as a function of $n / t$. For
example (see discussion on page \pageref{pg1}), we expect that 
for $\gamma _k > \gamma > \gamma _{k+1}$, there will be a sequence
of $2k$ speeds $s_1 > s_2 > \ldots > s_{2k}$, with the property
that for t large,
\begin{itemize}
\item 
for $s_2 < n/t < s_1$, the solution is a modulated one-phase wave,
\item 
for $s_3 < n/t < s_2$, the solution is a pure one-phase wave,
\item 
for $s_4 < n/t < s_3$, the solution is a modulated two-phase wave,
\end{itemize}
and so on, until
\begin{itemize}
\item 
for $s_{2k} < n/t < s_{2k - 1}$, the solution is a modulated $k$-phase wave,
\end{itemize}
and
\begin{itemize}
\item 
for $1/t \ll n/t < s_{2k}$, the region studied in this paper, 
we have a pure $k$-phase wave.
\end{itemize}
Note from the Figures \ref{fc4} -- \ref{fc6} that the 
above phenomena are present for
both small and large values of the amplitude of the periodic
component $h$ of the driver.  In the linear case, $F(x)=\alpha
x,\,\,\alpha>0$, the origin of the thresholds
$\gamma_1>\gamma_2>\ldots$ is simple to understand.  
\medskip
The solution of the lattice equations
\newline
\parbox{.45in}{\be \label{e.115} \ee}
\parbox{5.05in}{
\bea
\ddot{x}_n&=&F(x_{n-1}-x_n)-F(x_n-x_{n+1})=\alpha(x_{n+1}+x_{n-1}-2x_n),
\quad n\geq1,\\
x_n(0)&=&\dot{x}_n(0)=0,\quad n\geq1,\\
x_0(t)&=&2at+\sum\limits_{m\in\db Z}b_m\,\,e^{i\gamma
mt},\qquad b_{-m}=\bar{b}_m,
\eea
}
is easy to evaluate using Fourier methods and one sees that as
$t\rightarrow\infty$,
\be
x_n(t)\longrightarrow
2a(t-n)+\sum\limits_m\,\,b_m\,z^n_m\,e^{i\gamma mt},\quad0<n\ll t,
\label{e.120}
\ee
where $z_m,\,|z_m|\leq1$, is the root of
\be
z^2_m+\left(\frac{(\gamma m)^2}{\alpha}-2\right)z_m+1=0,\label{e.125}
\ee
chosen, in the case $|z_m|=1$, such that the energy is transported
{\bf away} from the driver.  Observe that if $m_0>0$ is the
largest integer $m$ for which $(\gamma m)^2/\alpha-2\leq2$, then
$|z_m|=1$ for $0\leq|m|\leq m_0$, and $|z_m|<1$ for $|m|>m_0$.
Inserting this information into (\ref{e.120}), we find that, away from the
driver, an $m_0$-phase wave propagates through the lattice in the
region $0<n\ll t$.  Thus the threshold values of $\gamma_k$ are
given, in this case, by
\be
\gamma_k=\frac{2\sqrt\alpha}{k},\quad k=1,2,\ldots\,\,\,.\label{e.130}
\ee
As we will see below, the above calculations are useful in
understanding the asymptotic state of the solution of (\ref{e.90}) as
$t\rightarrow\infty$ in the case that $h$ is small.
\par
In the case of the Toda lattice, when the driving is constant the
doubling-up trick converts the shock problem into an iso-spectral
deformation (\ref{e.75}) for the operator $\widetilde{L}$.  When $h$ is
non-zero, it is no longer clear how to convert the shock problem
(1.18) into an integrable form (although recent results of Fokas
and Its \cite{FI} suggest that this may still be possible to do).  As a
tool for analyzing (\ref{e.90}) in the Toda case we consider, rather,
the Lax pair of operators
\be
L=
\left(
\begin{array}{llll}
a_1&b_1&&0\\
b_1&a_2&b_2\\
&b_2&\ddots&\ddots\\
0&&\ddots\\
\end{array}
\right),\qquad
B=
\left(
\begin{array}{llll}
\,\,\,0&b_1&&0\\
-b_1&\,\,\,0&b_2\\
&-b_2&\ddots&\ddots\\
0&&\ddots\\
\end{array}
\right)\label{e.135}
\ee
for the semi-infinite lattice
$a_n=-\dot{x}_n/2,\,\,b_n=\frac{1}{2}\,e^{\frac{1}{2}(x_n-x_{n+1})},\,\,n\geq1$.
\par
A straightforward calculation shows that $L$ solves the equation
\be
\dot{L}=[B,L]-2\,b^2_0(t)P,\quad
b_0=\frac{1}{2}\,e^{\frac{1}{2}(x_0-x_1)},\label{e.140}
\ee
which we think of as a forced Lax system.  Here
$P=(P_{ij})_{i,j\geq1}$, is a matrix operator with
$P_{ij}=0$ unless $i=j=1$, and $P_{11}=1$.  The equation
describes a motion that is almost, but not quite, an iso-spectral
deformation of $L$.  As $t$ evolves, the essential spectrum of
$L(t)$ remains fixed,
$\sigma_{ess}(L(t))=\sigma_{ess}(L(0))=[a-1,a+1]$, but eigenvalues
may ``leak out'' from the continuum.  This is true, in particular,
in the case of constant driving $x_0=2at$, as was first observed
by Kaup and Neuberger \cite{KN}.
\par
In the case of constant driving with $a<1$, what happens to
$\sigma(L(t))$?   
\lfd
\begin{figure}[hpbt]
\leavevmode \epsfysize=6.0cm
\epsfbox{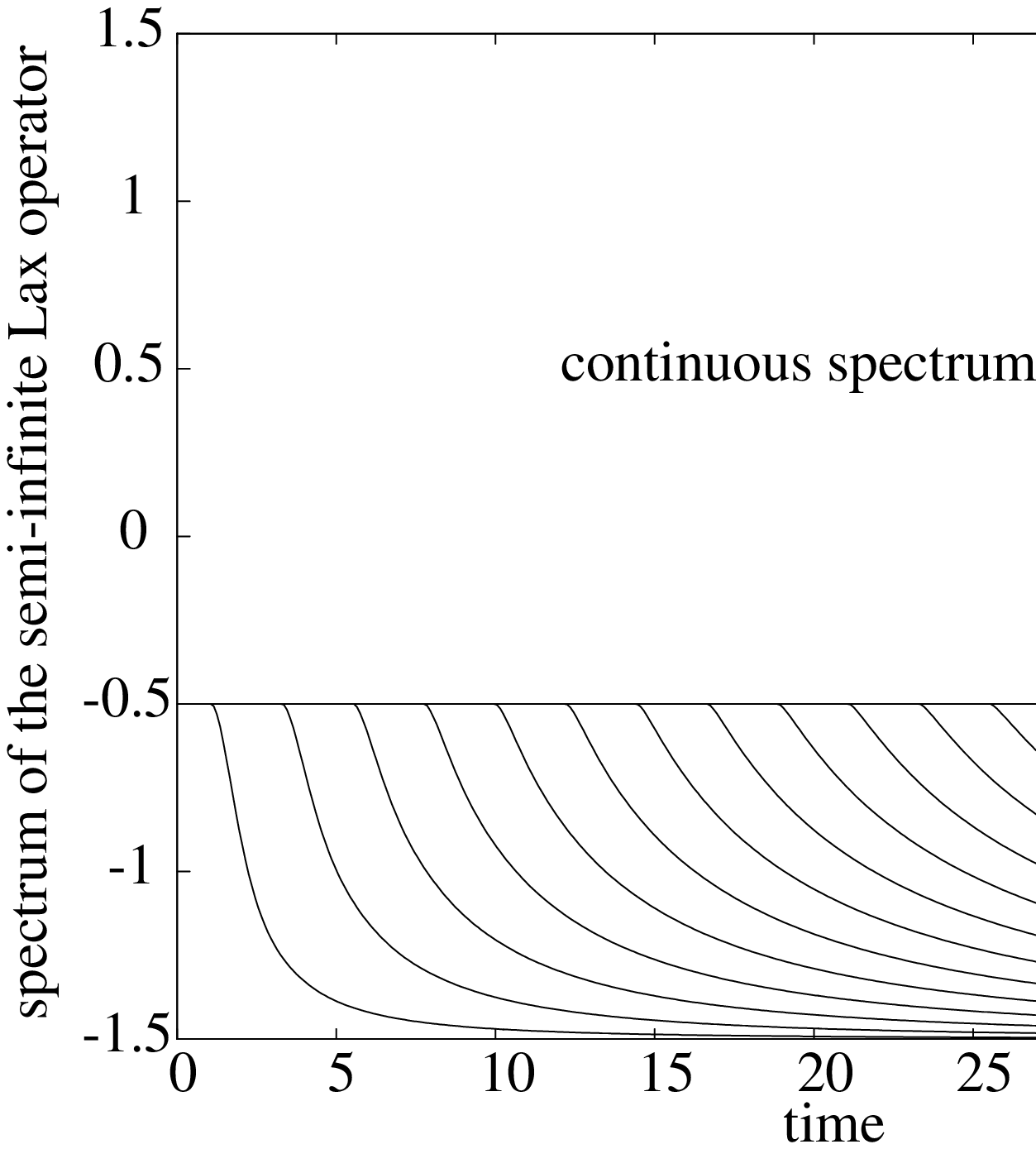}
\unitlength=1cm
\caption{{\em Evolution of $\sigma (L(t))$; driver: $x_0(t) = t$}}
\label{fe3}
\eef
We see in Figure \ref{fe3} that the eigenvalues emerge from 
the band $[a-1,a+1]$ and
eventually fill the larger band
$[-a-1,a+1]=[-a-1,-a+1]\cup[a-1,a+1]$.  (In the case $a>1$, the
bands $[-a-1,-a+1],\,[a-1,a+1]$ are disjoint and the spectrum of
$L(t)$ fills these two bands separated by a gap).
Thus this ``ghost'' band, which appeared as an
artifact of the solution procedure through the introduction of the
doubled-up operator $\widetilde{L}$, now emerges in real form,
populated by eigenvalues emerging from the original band
$[a-1,a+1]$. 
We learn from Figure \ref{fe3} that for $a<1$ 
there is no gap in the spectrum at
$t=\infty$, and (hence) there are no oscillations.
\par
In the periodically driven case, $x_0=2at+h(\gamma t)$, where
$\gamma>\gamma_1$ (and $a<1$), we find a similar picture to Figure
\ref{fe3} for the evolution of $\sigma(L(t))$ which is displayed in
Figure \ref{fe4} at some later time, so that more eigenvalues are present
than in Figure \ref{fe3}.
\lfd
\begin{figure}[hpbt]
\leavevmode \epsfysize=6.0cm
\epsfbox{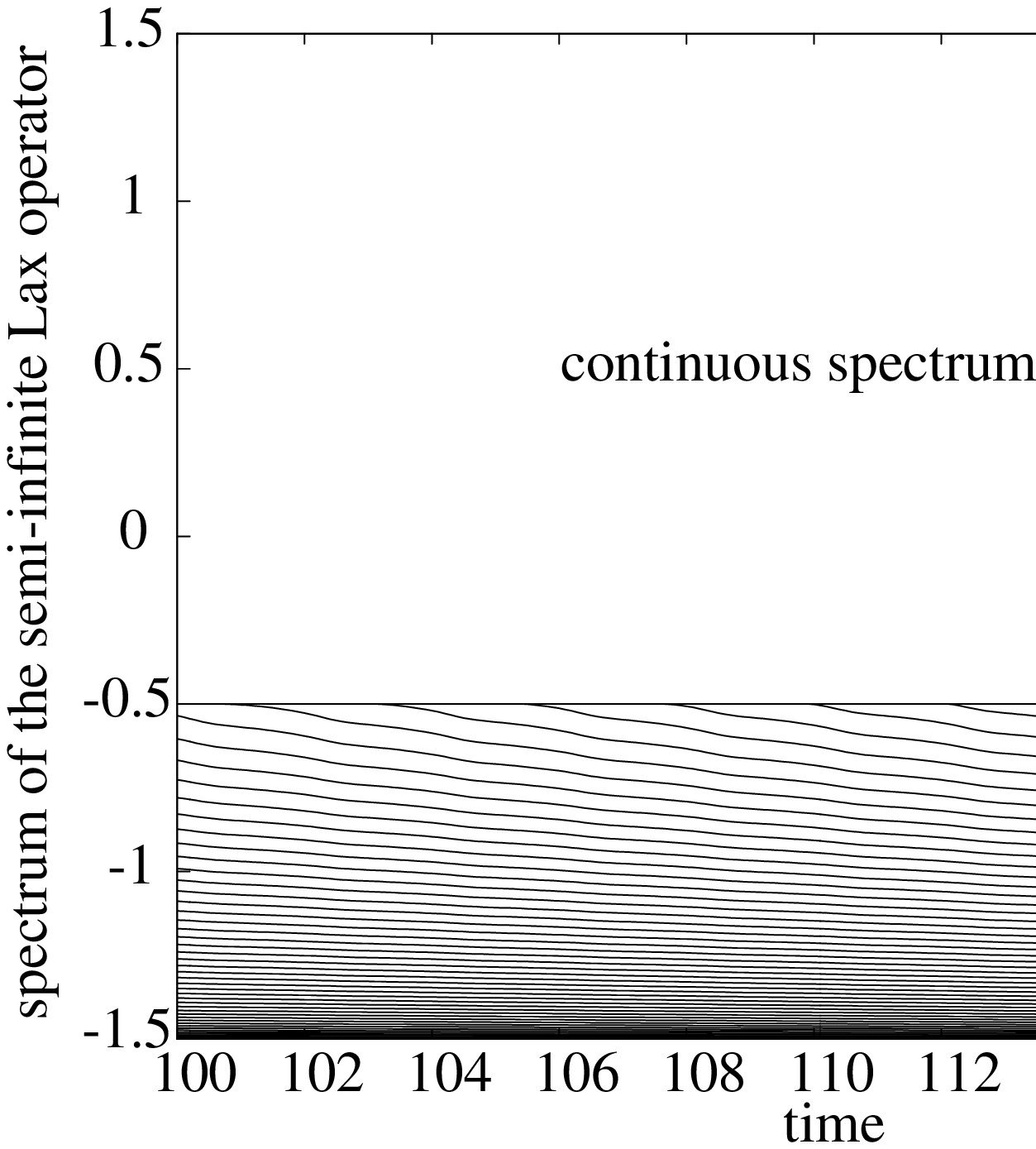}
\unitlength=1cm
\caption{{\em Evolution of $\sigma (L(t))$; driver: $x_0(t) = t + 
0.2(\sin \gamma t + 0.5 \cos 2 \gamma t), 
\gamma = 3.1 > \gamma _1$ }}
\label{fe4}
\eef
As $t\rightarrow\infty,\,\sigma(L(t))$ again converges to a single
band and no travelling wave emerges.  However, if
$\gamma_1>\gamma>\gamma_2$, we find different behavior
(Figure \ref{fe5}). 
\lfd
\begin{figure}[hpbt]
\leavevmode \epsfysize=6.0cm
\epsfbox{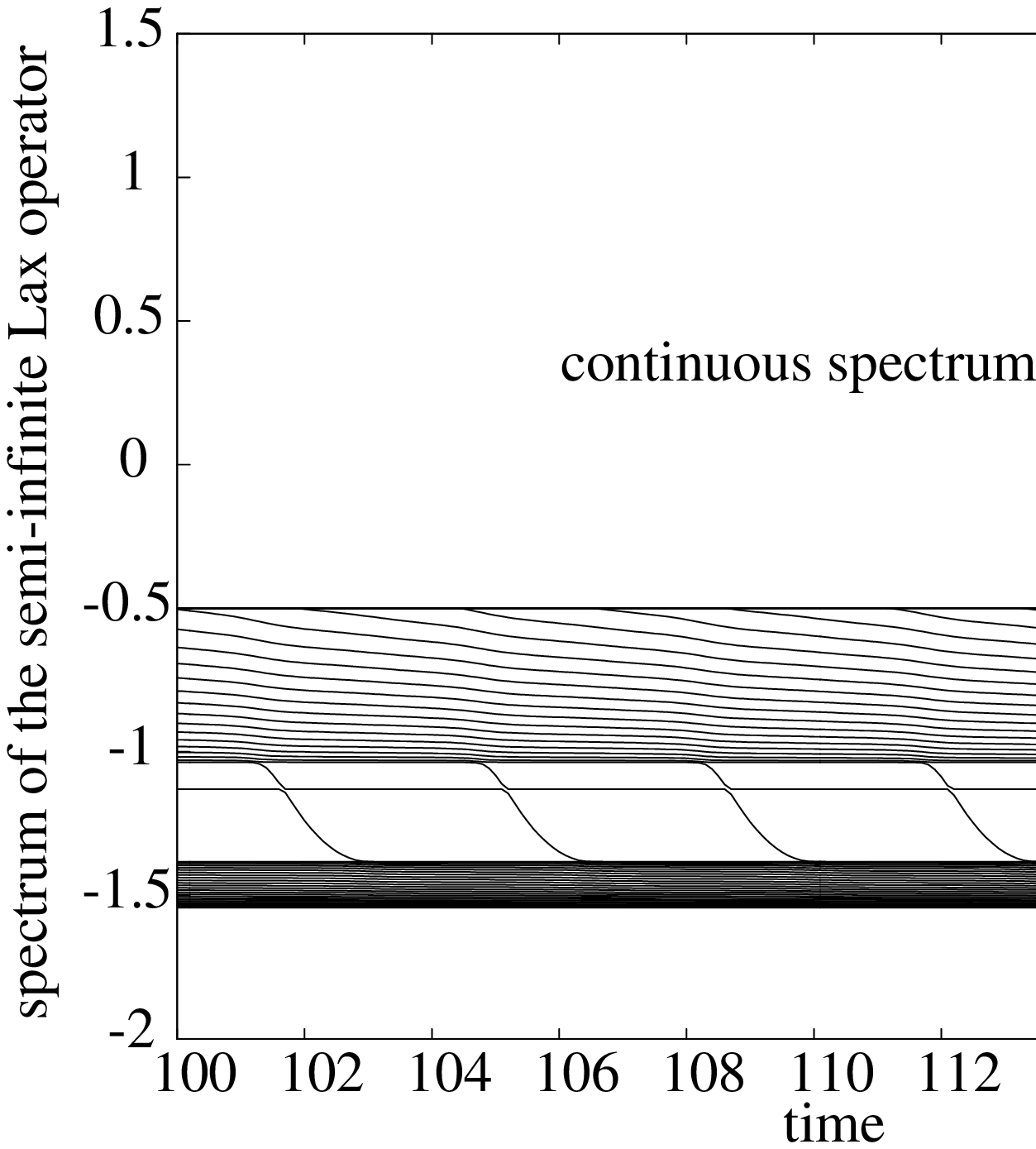}
\unitlength=1cm
\caption{{\em Evolution of $\sigma (L(t))$; driver: $x_0(t) = t + 
0.2(\sin \gamma t + 0.5 \cos 2 \gamma t), 
\gamma = 1.8, \gamma _1 > \gamma > \gamma _2$}}
\label{fe5}
\eef
We see that $\sigma(L(t))$ converges to two bands separated by a gap, and
a single phase wave emerges.  For $\gamma_2>\gamma>\gamma_3$, we
see from Figure \ref{fe6}
\lfd
\begin{figure}[hpbt]
\leavevmode \epsfysize=6.0cm
\epsfbox{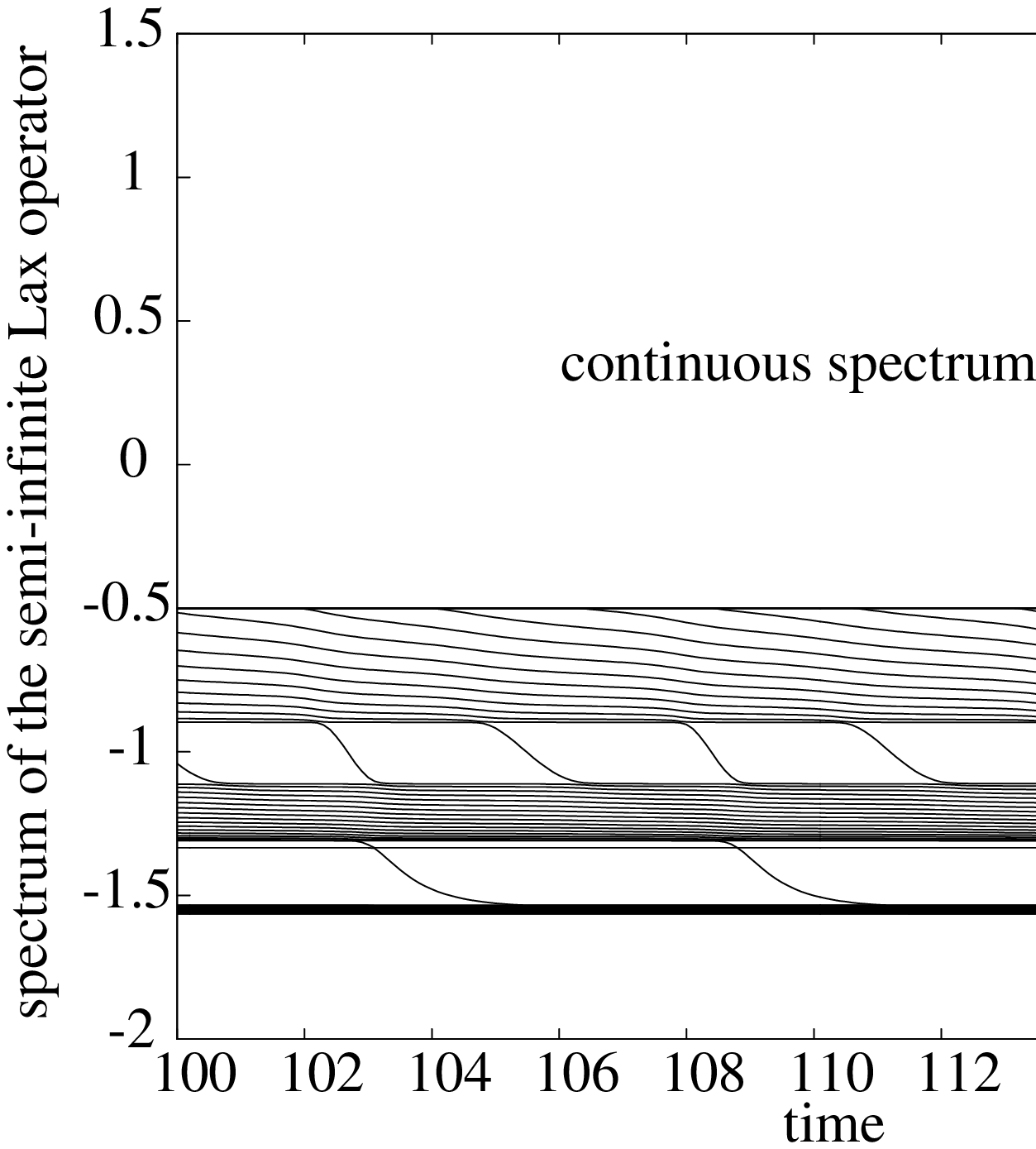}
\unitlength=1cm
\caption{{\em Evolution of $\sigma (L(t))$; driver: $x_0(t) = t + 
0.2(\sin \gamma t + 0.5 \cos 2 \gamma t), 
\gamma = 1.1, \gamma _2 > \gamma > \gamma _3$}}
\label{fe6}
\eef
that $\sigma(L(t))$ converges to three bands separated by two gaps,
and a two phase wave emerges, etc.
\vspace{.2in}
\newline
{\bf Remarks on the eigenvalues in the gap.}
\begin{itemize}
\item[(1)]
The eigenvalues which can be observed in the gaps of the 
spectrum of the semi-infinite Lax operator $L(t)$ 
(compare with Figures \ref{fe5} and \ref{fe6}) 
are of two different
types. They are either constant in time or they move down from the lower 
edge of one band to the upper edge of the next band. Eigenvalues of the 
second kind can be understood from the corresponding $g$-gap solution 
They are connected to the zeros of a theta function,
which is used in the construction of the $g$-gap solution 
(see Chapter \ref{g}).
Eigenvalues which are constant in time can be interpreted as follows.
Numerical computations show that they correspond to eigenvectors which are
moving out as $t \rightarrow \infty$. Hence the eigenvalues do not
survive in the spectrum of the limiting operator $L$, which
corresponds to a lattice where all particles perform time periodic
motion. In other words, from the spectral theoretic point of view,
this is an example of the general phenomenon 
that under strong convergence of operators
the spectrum is not necessarily conserved.
\item[(2)]
Figures \ref{fe5} and \ref{fe6} give the impression 
that the eigenvalue branches 
which come down cross the eigenvalues, that remain constant in time.
This, of course, is not possible as the symmetric, tridiagonal
operator $L$ cannot have double eigenvalues. Instead, a close look
demonstrates, that a billard ball collision is taking place as 
shown in Figure \ref{fe2}. It is possible to analyze the interaction
of $\lambda _k$ and $\lambda _{k+1}$ in detail by using equation
(\ref{s1.15}) of Chapter \ref{s} for $j=k$ and $j=k+1$, together
with the asymptotic assumption that $|\lambda _k - \lambda _{k+1}|$
is much smaller than the distance between any two other eigenvalues,
but we do not present the details. 
\end{itemize}
\bef
\leavevmode \epsfysize=4.0cm
\epsfbox{fe2.ps}
\caption{{\em Close look at a ``collision'' of two eigenvalues}}
\label{fe2}
\eef 
\par
For $\lambda<\inf\sigma_{ess}(L(0))$, 
an interesting quantity to compute is 
\be
J(\lambda)=\mathop{\lim}\limits_{t\rightarrow\infty}\,\,\,\frac{\sharp
\,\,\{{\rm eigenvalues}\,\,{\rm of}\,\,L(t)\,\,{\rm that}\,\,{\rm
are}
\,\,<\lambda\}}{t}.\label{e.165}
\ee
Clearly $J(\lambda)$ represents the asymptotic flux of  
eigenvalues of $L(t)$ across the value $\lambda$. In Chapter \ref{s} we will
extend the definition of $J$ to all values of $\lambda$.

It is observed numerically that $J(\lambda)$ indeed exists and
for $\gamma_2>\gamma>\gamma_3$, say, we find that $J(\lambda)$
looks as displayed in Figure \ref{fe7}
\lfd
\begin{figure}[hbt]
\leavevmode \epsfysize=6.0cm
\epsfbox{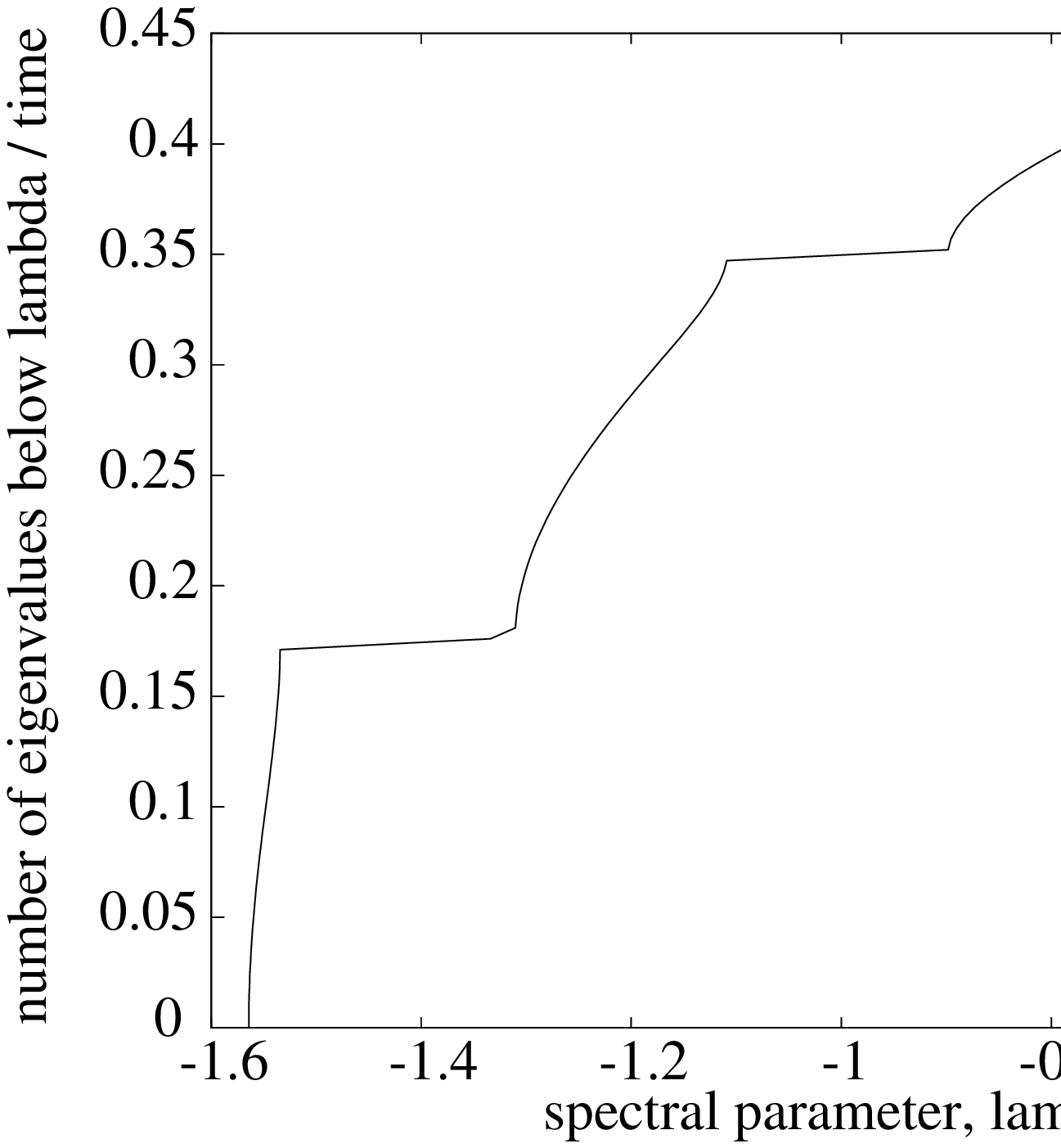}
\unitlength=1cm
\caption{{\em Numerically observed $J(\lambda)$ in the case
$\gamma _2 > \gamma > \gamma _3, \gamma = 1.1$}}
\label{fe7}
\eef
Thus $J(\lambda)$ is constant in the gaps and indeed we observe
more generally that  
\be
J(\lambda)=\frac{j\gamma}{2\pi}\quad
{\rm for}\quad
\lambda\,\,{\rm in}\,\,{\rm the}\,\, j^{\rm th}\,\,{\rm gap}.
\label{e.170}
\ee
This is a very intriguing fact, reminiscent of the Johnson-Moser
gap labelling theorem \cite{JM} in the spectral theory of
one-dimensional Schr\"odinger operators with almost periodic
potentials (see also the analogous gap labelling theorem for
Jacobi matrices \cite{B}, \cite{S}).
\par
Finally we are at the stage where we can describe our analytical
results, whose goal is to explain the above numerical experiments.
In the Toda case with constant driving, it was possible, using
the exact formulae of inverse scattering, to show that the
solution of the initial value problem converges as
$t\rightarrow\infty$ to the binary motion if $a>1$, and to a
quiescent lattice if $a<1$.  In the present case, where we no
longer have these formulae, our goals are more modest and we
restrict our attention to a description of the observed attractor.
Our results are the following:
\par
\noindent
{\bf I.  Strongly nonlinear case.}
\par
Here we consider (\ref{e.90}) in the case of the Toda lattice without
any smallness restriction on the size of the oscillatory component
$h$ of the driver $x_0$.  The main result is Theorem \ref{ts3.5} below in
which we show how to compute the normalized density of state
$J(\lambda)$ through the solution of a linear integral equation,
once the number and endpoints of the bands in $\sigma(L(\infty))$
are known.  This linear equation, in turn, can be solved
explicitly via an associated Riemann-Hilbert problem.
\par
At this stage it is not clear how to relate the number and
endpoints of the bands to the parameters of the problem
$a,\gamma,h$.  To test Theorem \ref{ts3.5} in any given situation, one
reads off the discrete information given by the number and endpoints
of the bands from the numerical experiment, and then compares the
solution of the integral equation with the normalized density of
states $J(\lambda)$ obtained directly from definition (1.33) using
the numerically computed eigenvalues of $L(t)$ at large times.
The numerical and analytical solutions for $J(\lambda)$ agree to
very high order:  see Appendix \ref{c} (Figures \ref{fc7}, \ref{fc8}) 
for further details.
\par
The proof of this result proceeds by deriving an equation of motion
(see (\ref{s1.15})) for the eigenvalues of a truncated version of $L(t)$ of
size $N\gg t$ as $t\rightarrow\infty$.  The continuum limit of the
time average of these equations, leads to the linear integral
equation (\ref{s2.30}) for $J(\lambda)$.
\par
\noindent
{\bf II.  Weakly nonlinear case.}
\par
Here we consider general $F$, but the periodic component $h$ is
now required to be suitably small.  From the numerical experiments
we see that if $h=0(\varepsilon)$, then as
$t\rightarrow\infty,\,\,x_n(t)$ converges to an asymptotic state
which is a $\frac{2\pi}{\gamma}$-time periodic solution
$x_{{\rm asymp},\,n}(t)$ with $x_{{\rm asymp},\,n}(t)=cn+0(\varepsilon)$ for
some lattice spacing $c$.  The goal here is prove that such time
periodic asymptotic states $x_{{\rm asymp},\,n}(t)$ indeed exist for
$\varepsilon$ small.  We proceed by linearizing around the
particular solution $x_{{\rm asymp},\,n}(t)=cn,\,\,n\geq0$, of the
equations $\ddot{x}_n=F(x_{n-1}-x_n)-F(x_n-x_{n+1}),\,\,n\geq1$,
and use various tools from implicit function theory.
\par
Our first result (Theorem \ref{tb3.1}) is a nonlinear version of the
classical linear method of separation of variables.  For example,
in solving the heat equation $u_t=u_{xx}$ on a half-line $x\geq0$
with boundary conditions at $x=0$, one proceeds by expressing the
solution as a combination
\bd
\int\,\,a_+(z)\,\,e^{-(tz^2+ixz)}dz+\int\,\,a_-(z)\,\,e^{-(tz^2-ixz)}
\ed
of elementary solutions $e^{-(tz^2\pm ixz)}$ of the heat equation
on the full line, and then choosing the parameters $a_+,a_-$ to
satisfy the boundary condition at $x=0$.  In the nonlinear case
(Theorem \ref{tb3.1}) we show that provided a sufficiently large parameter
family of travelling wave solutions of the doubly infinite lattice
\be
\ddot{x}_n=F(x_{n-1}-x_n)-F(x_n-x_{n+1}),\,\,\,-\infty<n<\infty,
\label{e.175}
\ee
exist, then (modulo technicalities, see Chapter \ref{m}) the parameters can
always be chosen to produce the desired asymptotic states
$x_{{\rm asymp},\,n}(t)$ of the driven semi-infinite problem.
\par
Thus the problem of the existence of the observed asymptotic states,
reduces to the problem of constructing sufficiently large
parameter families of travelling waves of the full lattice
equation (\ref{e.175}).  As we will see in Chapter \ref{m}, for
$\gamma_k>\gamma>\gamma_{k+1}\,\,k\geq1$, we will need
$2k$-parameter families of $k$-phase travelling waves of type
(\ref{e.110}) on the full lattice in order to construct the solution of
the driven lattice observed as $t\rightarrow\infty$ in the
numerical experiments.  If $\gamma>\gamma_1$ (see Section \ref{m4}) the
requirement of travelling wave solutions of (\ref{e.175}) trivializes,
and Theorem \ref{tb3.1} guarantees the existence of the desired asymptotic
states $x_{{\rm asymp},\,n}(t)$ of the driven lattice for sufficiently
small $h$ and general real analytic $F$ which are monotone in the
region of interest, and this explains Figure \ref{fc1}.  
\par
The next result (Theorem \ref{tl3.1}) concerns general $F$ in the case that
$\gamma_1>\gamma>\gamma_2$.  Here we show that a 2-parameter
family of one-phase travelling wave solutions of (\ref{e.175}) always
exist for general $F$.  Together with Theorem \ref{tb3.1}, this implies that
for $\gamma_1>\gamma>\gamma_2$ the desired states
$x_{{\rm asymp},\,n}(t)$ of the driven lattice exist, and this explains 
Figure \ref{fc2}.
This 2-parameter family is constructed by deriving an equation
for the Fourier coefficients of the travelling wave solution,
which can be solved by a Lyapunov-Schmidt decomposition.  The
infinite dimensional part does not pose any problems (see Lemma
\ref{ll2.1}) and the degenerate finite dimensional equations can be
solved by using certain symmetries of the equation (see Lemma
\ref{ll2.2}).
\par
If we try a similar construction for $m_0$-phase waves, $m_0>1$,
then we encounter in the infinite dimensional part of the 
Lyapunov-Schmidt decomposition, a small divisor problem related to the
small divisor problem occurring in \cite{CW}, where periodic 
solutions of the
nonlinear wave equation are constructed, and which we hope to
solve in the near future.  In the special case of Toda, however,
the family of travelling waves can be constructed explicitly.
Indeed in our third, and final, result (Theorem \ref{tg4.1}) we use the
integrability of the doubly infinite Toda lattice and show how the
well-known class of $g$-gap solutions contains a sufficiently 
large family
of travelling waves to apply to Theorem \ref{tb3.1} and so construct the
desired asymptotics states $x_{{\rm asymp},\,n}(t)$ of the driven
lattice for any $\gamma\in\db
R_+\backslash\{\gamma_1,\gamma_2,\cdots\}$.  
\\
\\
\\
Finally we want to pose four open problems, which are connected to our 
investigations, some of which were mentioned above.
\begin{itemize}
\item[(1)]
The ``critical shock'' phenomena.

As discussed in the very beginning of the introduction, Holian and 
Straub have numerically discovered a critical shock strength
$a_{{\rm crit}} (F)$ in the case of constant driving velocity
$x_0(t) = 2 a t$. For $a < a_{{\rm crit}} (F)$ the lattice comes to rest 
behind the shock front as $t \rightarrow \infty$, whereas 
for $a > a_{{\rm crit}} (F)$ the particles of the lattice will execute
binary oscillations as $t \rightarrow \infty$.
So far this result has been analytically explained in the case
of the Toda lattice ($F(x) = e^x $) (cf \cite{HFM}, \cite{VDO})
and can be seen to be absent for linear force functions by
explicit calculation.

The question is to find general conditions on the force $ F $
for which one can prove the existence of a critical shock strength.
\item[(2)]
Existence of multi-phase travelling waves.

Let $ F, c $ satisfy the general assumptions (cf Section \ref{m22}) and
let $ \gamma \in \relz $ satisfy
$\frac{2}{k} \sqrt{F'(-c)} > \gamma > \frac{2}{k+1} \sqrt{F'(-c)} $
for some $k \in \natz $.

Does there exist a smooth $2k$-real parameter family of solutions
\be
x_n(q)(t) = 
c n + \chi _k (q) (n \beta _1 (q) + \gamma t, \ldots,
n \beta _k (q) + k \gamma t), 
\label{e.500}
\end{equation}
for $q \in \relz^{2k}, q$ small, 
of the equation
\be
\ddot{x} _n (t) = F(x_{n-1} - x_n) - F(x_n - x_{n + 1}), n \in \ganz,
\label{e.505}
\end{equation}
where $\chi _k (q)$ is for each $q$ a function periodic in its $k$
arguments, $\chi _k (0) = 0$ and $D_q \chi _k (0)$ has maximal rank $2k$?
Note that these solutions exist in the case of the Toda lattice and 
are given by (\ref{g3.110}).
For general force functions the work of Craig and Wayne (\cite{CW})
indicates that it might be necessary to aim for a slightly weaker 
result, namely, that the smooth family of functions $x_n(q)(t)$ of
the form (\ref{e.500}) are solutions of (\ref{e.505}) only for
a Cantor set in the parameter space $q \in \relz ^{2k}$, which 
has almost full measure.
\item[(3)]
The case $a > a_{{\rm crit}}$.

In the paper we always assume that the driver is of the form
$x_0(t) = 2at +h(\gamma t)$ with $a < a_{{\rm crit}}$. 
In the case of the 
Toda lattice we have also conducted some experiments for 
$a > a_{{\rm crit}} (Toda) = 1$. We have made the following observation:
for $h$ small, the limiting operator $L(t), t\rightarrow \infty$ seems
to have infinitly many gaps and again we obtain a version of 
gap labelling. In fact, for all the gaps we have observed that one can 
write the numerically determined integrated spectral density
$J(\lambda)$ (cf \ref{e.165}) in the form
\be
J(\lambda) = \frac{j \gamma + k \omega}{2 \pi},\quad  j,k \in \ganz, \quad
\lambda \mbox{ lies in a gap, }
\label{e.510}
\ee
where $\gamma$ denotes the frequency of the driver and $\omega$
is given by the frequency of the time-asymptotic oscillations,
which are observed in the case that the driver has constant speed
$x_0(t) = 2at, a > 1$. 

Corresponding to our results in Chapter \ref{g},
we ask whether it is possible to construct solutions of the driven
semi-infinite Toda lattice with driver $x_0(t) =  2 a t + h(\gamma t)$,
such that the spectrum of the corresponding Lax-operator
has infinitely many gaps and bands.
\item[(4)]
Connection to the initial value problem.

All of our results were motivated by the numerically observed 
long-time behavior of a certain initial value problem (with 
shock initial data). However, so far we are not able to prove 
from first principles, that the solution of the initial value 
problem actually converges as $t \rightarrow \infty$ to one of the 
asymptotic states described in 
Chapters \ref{s}-\ref{g}. This basic problem remains open and,
alas, seems far from a resolution.
\end{itemize}
\noindent
{\bf Acknowledgments.}
\par
The work of the first author was supported in part by NSF Grant
No. DMS--9203771.  The work of the second author was supported in
part by an Alfred P. Sloan Dissertation Fellowship.  The work of
the third author was supported by ARO Grant No. DAAH04-93-G-0011
and by NSF Grant No. DMS-9103386-002. The authors would
also like to acknowledge the support of MSRI and the Courant
Institute in preparing this videotext.  Finally, the authors are
happy to acknowledge useful conversations with many of their
colleagues, and in particular, with Gene Wayne and Walter Craig.
Also we would like to thank Fritz Gesztesy for making available
his very useful notes (\cite{Ges}, \cite{Ges2}) on the application
of the theory of Riemann surfaces to integrable systems.

\chapter{An asymptotic calculation in the strongly nonlinear case}
\label{s}

\section{The evolution equations}
\label{s1}

We recall from the Introduction the Flaschka variables (see (\ref{e.60}))
\be 
a_n=-\frac{\dot x_n}{2},\quad b_n=\frac 12 e^{(x_n-x_{n+1})/2}, 
\quad  n=0,1,2,\dots,
\label{s1.5}
\ee
\begin{flushleft}
\hspace{2in}$a_n\rightarrow a, \quad b_n\rightarrow \frac{1}{2} \quad \mbox{ as } ~n\rightarrow +\infty,$
\end{flushleft}
and we note that the function $a_0(t)$ is the given time-periodic
forcing function.  For the semi-infinite Toda chain with $F(x)=e^{x}$, 
equation (\ref{e.5}) reduces to the perturbed Lax pair equation
\be
\frac{dL}{dt} + LB - BL = -\rho (t)P,
\label{s1.10}
\ee
where $L$ is the tridiagonal operator (cf (\ref{e.135}))
$$
L = \left( \begin{array}{ccccc}
a_1 & b_1 &&&0\\
b_1 & a_2 & b_2 & \\
& b_2 & a_3 & b_3\\
0&&\ddots & \ddots & \ddots
\end{array} \right) .
$$

\noindent $B$ is the antisymmetric tridiagonal operator given by
$$
B=\left( \begin{array}{ccccc}
0 & b_1&&&0\\
-b_1 & 0 & b_2\\
 & -b_2 & 0 & b_3\\
0&&\ddots & \ddots & \ddots
\end{array}
\right).
$$
P is the rank-one matrix given by:

$$
P_{ij} = \left\{ \begin{array}{lll}
1, & \mbox{ if } \qquad i=j=1\\
0, & \mbox{ otherwise }
\end{array}\right.$$
and $\rho(t)$ is the function:
\be
\rho(t) = 2b_0^2(t) = 2b_1^2(t) - \dot a_1(t),\quad \cdot = \frac{d}{dt} \;.
\label{s1.12}
\ee

The matrices $L$ and $B$ are semi-infinite.  We truncate the chain at some particle of 
very large index $N$, and work with the truncated finite matrices $L_N$ and $B_N$. 
The disturbance in the chain caused by the truncation, travels essentially with finite 
velocity. Only exponentially small effects display infinite speed. The bulk of the chain 
essentially 
does not feel the truncation until a time $t=O(N)$.  Thus, we expect that the finite 
system is a good approximation to the full semi-infinite system in the space-time 
region $1<<t<<N$ and $n<< N$.   In what follows, when we take the limit as $t\rightarrow\infty$, 
we always understand $t\to\infty, ~ N\to\infty, ~\frac{t}{N}\to 0$.  

We remark that on $\sigma_{ess}(L(t))=\sigma_{ess}(L(0))$, by standard spectral methods, 
the matrix $L_{N}(t)$ has a set of discrete eigenvalues tightly packed at densities of 
order $O(N)$.  On the other hand, we expect that the discrete spectrum of $L(t)$, 
which emerges from $\sigma_{ess}(L(0))$ as described in the introduction, is well 
approximated by the eigenvalues of ~$L_{N}(t)$ which lie outside 
$\sigma_{ess}(L(0))$.  This is because the associated eigenvectors are typically 
exponentially decreasing in $n$ and hence do not feel the truncation at $N$.

A word of explanation:  Typically an eigenvalue of $L_{N}(t)$ starts off as an 
eigenvalue of $L_{N}(0)$ lying in $\sigma_{ess}(L(0))$.  As $t$ increases, 
the eigenvalue moves with velocity $O(\frac{1}{N})$ until it emerges from 
$\sigma_{ess}(L(0))$.  It is \underline{only} after this point that the motion 
of the eigenvalue becomes relevant to the evolution of the discrete spectrum 
of $L(t)$.

Our strategy is to derive evolution equations for
\begin{enumerate}
\begin{description}
\item{(a)} the eigenvalues $\lambda_{j}^{N}$ of the truncated matrix $L_{N}$,
\item{(b)} the first entry $f_j$ of the $j^{th}$ eigenvector of $L_{N}~(j=1,\dots,N)$ 
when it is normalized to have Euclidean length equal to one.
\end{description}
\end{enumerate}

It is well known that the set $\{\lambda_{j},f_{j}\}^{N}_{j=1}$ determines the tridiagonal 
matrix $L_{N}$.
\lfd
\begin{theorem}
\label{ts1.5}
The evolution of the $\lambda_j$'s and $f_j$'s is given by:
\be
\left\{ \begin{array}{ll}
\frac 12\frac{d}{dt}\ln(-\dot\lambda_j)=\lambda_j-a_0(t)+
\sum^N_{\stackrel{i=1}{\scriptscriptstyle{i\ne j}}}
\frac{\dot \lambda_i}{\lambda_j-\lambda_i},\quad ~j=1,\dots,N\;,\\
f^2_j = \frac{-\dot\lambda_j}{\rho}\;,
\end{array}\right.
\label{s1.15}
\ee
where $\rho=2b_0^2(t) = -\sum^N_{i=1} \dot\lambda_i$. 
The initial values $\lambda_i(0)$ are the eigenvalues 
of $L_N$ at $t=0$ while the initial values $\dot\lambda_i(0)$ are given by
\be
\dot\lambda_i(0)=-2b_0^2(0)f_i^2(0).
\label{s1.17}
\ee
\end{theorem}
\begin{proof} 
Let $\Lambda$ be the diagonal matrix of the eigenvalues $\lambda_j$ of $L_N$ and let 
$\Psi$ be the orthogonal matrix whose $j^{th}$ column is the normalized eigenvector 
of $L_N$ corresponding to the eigenvalue $\lambda_j$. We have:
\be
L_{N}\Psi = \Psi\Lambda ,
\label{s1.20}
\ee
and we define the matrix $\Phi$ by:
\be
\Phi =\dot\Psi - B_{N}\Psi .
\label{s1.25}
\ee
Utilizing equations (\ref{s1.20}) and (\ref{s1.25}) and (\ref{s1.10}), 
we easily calculate:
\be
L_{N}\Phi - \Phi\Lambda=\Psi\dot\Lambda + \rho P\Psi .
\label{s1.30}
\ee
We now define the matrix $A= (a_{ij})$ by 
\be
A=\Psi^{-1} \Phi = \Psi^{T}\Phi.
\label{s1.35}
\ee
We calculate 
\bea
A+A^{T}&=&\Psi^T\Phi+\Phi^T\Psi \\
&=&\Psi^T(\dot\Psi-B_N\Psi)+(\dot\Psi^T-\Psi^TB_{N}^T)\Psi \\
&=&
\Psi^T\dot\Psi+\dot\Psi^T\Psi=\frac{d}{dt}(\Psi^T\Psi)=0.
\eea
Thus
\be
A+A^{T}=0;
\label{s1.40}
\ee
i.e. $A$ is antisymmetric. Using (\ref{s1.30}) we obtain
\be
L_{N}\Phi-\Phi\Lambda=L_{N}\Psi A-\Psi A\Lambda =\Psi(\Lambda A-A\Lambda).
\label{s1.45}
\ee
Comparing (\ref{s1.30}) with (\ref{s1.45}) we obtain easily:
\be
\dot\Lambda =[\Lambda,A]-\rho\Psi^T P\Psi .
\label{s1.50}
\ee
Let $f^T=(f_1,f_2,\dots,f_N)$ be the first row of $\Psi$. 
Then $\Psi^T P\Psi= ff^T$.  We insert this relation in (\ref{s1.50}),
\be
\dot\Lambda=[\Lambda,A]-\rho ff^T.
\label{s1.55}
\ee
Equating the diagonal elements on both sides we obtain:
\be
\dot\lambda_j =-\rho f^2_j,\qquad \sum^N_{j=1} \dot\lambda_j =-\rho.
\label{s1.60}
\ee
This proves the second relation in Theorem \ref{ts1.5}. Furthermore we note
that the first components $f_j$ of the eigenvectors of the tridiagonal matrix
$L_N$ do not vanish and we conclude by (\ref{s1.60}) that $- \dot{\lambda} _j
> 0$. Hence the first relation in Theorem \ref{ts1.5} is well defined.

Off the diagonal in (\ref{s1.55}) we have 
$\lambda_ia_{ij}-a_{ij}\lambda_j = \rho f_if_j$. 
On the other hand $a_{ii}=0$ by (\ref{s1.40}).  Thus
\be
\begin{array}{l}\left\{
\begin{array}{lll}
a_{ij} &=& \frac{\rho f_if_j}{\lambda_i-\lambda_j}, 
\quad \mbox{ when }\quad i\ne j,\\
a_{ii}&=&0.
\end{array} \right.
\end{array}
\label{s1.65}
\ee
We now calculate the evolution of $f_j$.  By (\ref{s1.25}):
\bd
\dot\Psi =\Phi+B_{N}\Psi=\Psi A+B_{N}\Psi.
\ed
Specializing to the first row we obtain
$\dot f^T=f^TA + B_{R_1}\Psi$, where $B_{R_1}$ 
is the first row of $B_N$.  This implies 
\bea
\dot f^T&=&f^TA+(L_{N}-a_1I)_{R_1}\Psi=f^TA+(L\Psi)_{R_1}-a_1f^T \\
&=& f^TA+(\Psi\Lambda)_{R_1}-a_1f^T = f^TA+f^T\Lambda-f^Ta_1,
\eea

or, taking transposes
\be
\dot f=(\Lambda-a_1I-A)f.
\label{s1.70}
\ee
But, by (\ref{s1.65}),
\be
A=\rho FDF,
\label{s1.75}
\ee
where $F$ is the invertible diagonal matrix  with entries $f_1,\dots,f_N$ 
and $D$ is the matrix $\left(\frac{1}{\lambda_i-\lambda_j}\right)$ with
zero entries on the diagonal. 
Thus $\dot f=(\Lambda-a_1I-\rho FDF)f$, and hence
\be
F^{-1}\dot f=F^{-1}\Lambda f-a_1F^{-1}f-\rho DFf.
\label{s1.80}
\ee
Now note that 
\begin{enumerate}
\item[(i)] $\rho Ff=\rho\left(\begin{array}{cc}f_1^2\\f_2^2\\ 
\vdots\end{array}\right)=-\left( \begin{array}{cc}\dot\lambda_1\\ 
\dot\lambda_2\\ \vdots\end{array}\right)$ by (\ref{s1.60}), 
\item[(ii)] $F^{-1}f=\left(\begin{array}{cc} 1\\1\\ 
\vdots\end{array}\right)$,~
$F^{-1}\Lambda=\Lambda F^{-1}$.
\end{enumerate}

Substituting in (\ref{s1.80}) we obtain
\be
\frac{\dot f_j}{f_j}=\lambda_j-a_1+
\sum\limits^N_{\stackrel{i=1}{\scriptscriptstyle{i\ne j}}}
\frac{\dot\lambda_i}{\lambda_j-\lambda_i}.
\label{s1.85}
\ee
The evolution of the $\lambda_i$'s in (\ref{s1.15}) is finally obtained 
by eliminating
$f_j$ between (\ref{s1.85}) and (\ref{s1.60}) and 
using the expression for $\rho$ given in (\ref{s1.12}).
\end{proof}

\section{The continuum limit of eigenvalue dynamics}
\label{s2}

The results of the numerical experiments described in the Introduction, 
(cf Figures \ref{fe3}-(\ref{fe6} lead us to consider the flux 
of eigenvalues of the matrix $L_{N}$ across a value $\lambda$.  
Noting that the eigenvalues of $L_{N}$ move toward lower values 
$(\dot{\lambda}_{i}=-\rho f^{2}_{i}<0)$, we define the eigenvalue 
flux at $\lambda$ averaged over a time interval $(t, t+T)$ by
\be
J_{t,T}(\lambda)=\frac{1}{T}~\mbox{ card }~\{{i}:
\lambda_{i}(t+T) < \lambda<\lambda_{i}(t), ~~i=1, 2\ldots N\}.
\label{s2.5}
\ee
We pose the following ansatz.
\lfd
\begin{ansatz}
\label{as2.5}
There exists a continuous, almost everywhere continuously differentiable 
function $J(\lambda)$ such that
\be
J_{t,T}(\lambda)\to J(\lambda),\quad\mbox{ when }
\quad T,t,N\to\infty \quad\mbox{ subject to }\quad T<t \ll N.
\label{s2.10}
\ee
\end{ansatz}
When $\lambda < inf~\sigma_{ess}(L)$, it is clearly true that
\be
J(\lambda)=\lim\limits_{t\to \infty} \D\frac{\mbox{ \# eigenvalues of }~ L(t)~\mbox{ that are smaller than }~\lambda}{t}~, 
\label{s2.15}
\ee
as defined in the Introduction.  The net gain in eigenvalues of $L_N$ of 
an interval $({\lambda},\hat{\lambda})$ over a long time $T$ is given 
asymptotically by $[J(\hat{\lambda})-J(\lambda)]T.$  Dividing by 
$(\hat{\lambda}-\lambda)T$ and letting $\hat{\lambda}\to\lambda$ 
we obtain that the asymptotic rate of increase in eigenvalue 
concentration at $\lambda$ is given by $\frac {dJ}{d\lambda}=J'(\lambda)$; 
thus, the difference in eigenvalue concentration at $\lambda$ between 
times $t$ and zero is asymptotically  $tJ'(\lambda)$.  When 
$\lambda < inf~\sigma_{ess}(L)$, necessarily $J'(\lambda)\ge 0$ 
since there is no original eigenvalue concentration at $\lambda$. 
On the other hand $J'(\lambda)$ can take negative values when
$\lambda \in \sigma_{ess}(L)$.

We will now use the function $J(\lambda)$ and some assumptions based on 
numerical observations to derive the continuum limit of the eigenvalue 
evolution equations (\ref{s1.15}).  We begin by averaging the system 
(\ref{s1.15}) of equations $(j=1,\ldots N)$ over the time interval 
$(t, t+T)$ to obtain
\be
\D\frac{1}{2T} \ln \D\frac{\dot{\lambda}_{j}(t+T)}{\dot{\lambda}_{j}(t)}=
\label{s2.20}
\ee
\bd
\D\frac{1}{T}\int^{t+T}_{t} \lambda_{j} (t')dt' -
\D\frac{1}{T}\int^{t+T}_{t}a_{0}(t')dt' + 
\D\frac{1}{T}\sum\limits^{N}_
{i=1\atop i \ne j}\int^{t+T}_{t}
\D\frac{d\lambda_{i}(t')}{\lambda_{j}(t')-\lambda_{i}(t')},\quad j=1,\ldots N.
\ed
Let $\lambda$ satisfy $J'(\lambda)\neq 0$, and let $j=j(t)$ be such that 
in the asymptotic limit $1\ll T\ll t\ll N$ (note that we require $T\ll t$, 
not just $T<t$ as in (\ref{s2.10})) the following is true,
\be
\lambda_{j(t)}(t')\to\lambda \mbox{ for any } t'
\mbox{ that satisfies } t<t'<t+T.
\label{s2.25}
\ee
In practical terms, this means that we expect eigenvalues $\lambda_j$ 
to stay close to the value $\lambda$ throughout the time interval 
$[t, t+T]$.  This fact is clearly borne out in the results of numerical 
experiments as long as $J'(\lambda)\neq 0$. 

We make two more simplifying assumptions when $J'(\lambda)\neq 0$
that are again justified by numerical experiments:
\begin{itemize}
\begin{enumerate}

\item[(a)] The left hand side of (\ref{s2.20}) is negligible.
\item[(b)] The ``singular contribution'' in the sum of the right hand 
side corresponding to indices $i$ that are close to $j$ is also negligible.  
In practical terms we interpret this to mean that the limiting integral 
kernel is the Hilbert transform.

\end{enumerate}
\end{itemize}
Under these conditions we can take the limit $ 1\ll T\ll t \ll N$ 
in (\ref{s2.20})-(\ref{s2.25}).
\lfd
\begin{theorem}
\label{ts2.5}
(Continuum Limit of (\ref{s2.20})-(\ref{s2.25})).  
Under Ansatz (\ref{as2.5}) and under the further assumption described 
above we have:
\be
J'(\lambda) \neq 0 \Rightarrow \lambda - <a_{0}>-
\D\int\limits^{\infty}_{-\infty}\!\!\!\!\!\!\!\!= 
\frac{J(\mu)}{\lambda - \mu} d\mu = 0,
\label{s2.30}
\ee
where $<a_{0}>$  is the mean value of the periodic driver $a_{0}(t)$, and as 
usual the double bars on the integral indicate that the principal value is 
taken.
\end{theorem}
\begin{proof}
By the assumptions, and by (\ref{s2.25}) the only thing to be shown is 
that the 
sum in (\ref{s2.20}) tends to the integral in (\ref{s2.30}).  
If we partition the 
eigenvalue axis into a set of infinitesimal intervals and if 
$(\mu - d\mu, \mu)$ is such an interval then the contribution 
$\D\frac{1}{T}\left(\frac{-d\mu}{\lambda-\mu}\right)$ should 
arise in as many terms of the sum in (\ref{s2.25}) as there are eigenvalues 
that cross the value $\mu$ during the time interval $(t,t+T)$.  
This number is asymptotically $TJ(\mu)$.  The sum in (\ref{s2.25}) therefore 
tends to $-\D\int\limits^{\infty}_
{-\infty}\!\!\!\!\!\!\!\!=\frac{J(\mu)d\mu}{\lambda-\mu}$.
\end{proof}
 
\section{ The asymptotic spectral density of $L_{N}$}
\label{s3}

We now consider the problem of determining $J(\lambda)$.  Our solution 
is partial in the sense that we can calculate the function 
$J(\lambda)$ if we are given the set 
$\{\lambda:  J'(\lambda) \neq 0\}$.  Numerical calculations 
(see Figures \ref{fe5}, \ref{fe6}, \ref{fe7}) 
show that this set is a finite 
union of intervals.  Thus, we are assuming 
knowledge of a finite set of numbers that are in 
principle determined by the fluctuating part of 
the periodic driver $a_{0}(t)$.  Determining these 
numbers is, unfortunately, the part of the problem 
that we have not yet been able to solve.  

We proceed to give some basic definitions.
\lfd
\begin{definition}
\label{ds3.5}
Let the points  $p_{0}<q_{1}<p_{1}<q_{2}<p_{2}<\ldots <p_{g}<q_{g+1}$ 
be given.  These are $2g+2$ points in all.  We define the set of 
{\bf bands} $B$, where we have  $J'(\lambda)\neq 0$,
\be
B=[p_{0},q_{1}]\cup[p_{1},q_{2}]\cup \cdots \cup[p_{g},q_{g+1}],
\label{s3.5}
\ee
and the set of {\bf gaps} $G$, where $J'(\lambda)=0$, by
\be
G_{k}=(q_{k}, p_{k}),~ k=1,2,\ldots g;\qquad G=\cup^{g}_{k=1} G_{k} 
\label{s3.10}
\ee
(cf Figure \ref{fs1} below). We then define the hyperelliptic curve
\be
R(\lambda)=\{\Pi^{g}_{k=0}(\lambda-p_{k})
(\lambda-q_{k+1})\}^{1/2}
\label{s3.15}
\ee
with branch cuts along the set B and sign determination such that 
$R(\lambda)>0$ when $\lambda \to +~\infty$.  Finally, we define the polynomial
\be
P(\lambda)=\D\sum\limits^{g+1}_{i=1}(\lambda-\sigma_{i}) 
\quad,\sigma_{i}\in {\bf R},
\label{s3.20}
\ee
where the $\sigma_{i}\mbox{'s}, ~i=1,\ldots g+1$, are uniquely 
determined by the $g+1$ relations.
\be
\int\limits_{G_{k}}\D\frac{P(\lambda)}
{R(\lambda)}d{\lambda}=0\quad,k=1,2,\ldots g.
\label{s3.25}
\ee
\be
\int\limits_{B}\D\frac{P(\lambda)}{R(\lambda)}d{\lambda}=0.
\label{s3.30}
\ee
\end{definition}

\bef
\leavevmode \epsfysize=2.5cm
\epsfbox{fs1.ps}
\caption{{\em The band -- gap structure of the spectrum}}
\label{fs1}
\eef

We observe that the integrals in (\ref{s3.25}) and (\ref{s3.30}) can be 
easily understood as contour integrals on the Riemann surface 
associated with $R$.  The contour in (\ref{s3.30}) can be replaced by 
a circle of (large) radius.  We then obtain 
$\frac{P(\lambda)}{R(\lambda)}=1 + O(\frac{1}{\lambda^{2}}) 
\quad\mbox{ as }\quad \lambda\to\infty$, which implies through 
an easy asymptotic calculation that
\be
2\D\sum\limits^{g+1}_{i=1}\sigma_{i}=\D\sum\limits^{g}_{i=0}(p_{i}+q_{i+1}).
\label{s3.35}
\ee
\lfd
\begin{theorem}
\label{ts3.5}
Let $J(\lambda)$ be a continuous function supported on the set $B\cup G$, 
differentiable at all points, except possibly the boundary points of $B$, 
satisfying
\be
\lambda - <a_{0}> - \D\int\limits^{\infty}_{-\infty}\!\!\!\!\!\!\!\!=
\D\frac {J(\mu)}{\lambda-\mu} d\mu =  0\qquad  \lambda\in B,
\label{s3.40}
\ee
\be
J(\lambda)  =  c_{k} = \mbox{ const., }\qquad  \lambda\in G_{k}, ~k=1\ldots, g,
\label{s3.45}
\ee
\be
J(\lambda) =  0 \qquad \lambda \not\in B\cup G,
\label{s3.50}
\ee
where $<a_{0}>$ is a constant. Then $J(\lambda) + iHJ(\lambda)$ 
is the limiting value as $z\to \lambda+i0$, of the analytic function
\be
f(z)=\frac{-i}{\pi}\D\int\limits^{\infty}_{z}1-
\D\frac{P(z')}{R(z')}dz',\quad {\bf Im}~z\neq 0.
\label{s3.55}
\ee
Precisely:
\be
J(\lambda)=\D\frac{1}{\pi}{\bf Im}
\D\int\limits^{\infty}_{\lambda}\D\left(1-\frac{P(\mu)}{R(\mu)}\right)d\mu,
\label{s3.60}
\ee
\be
HJ(\lambda)=\frac{1}{\pi}\D\int\limits^{\infty}_{-\infty}
\D\frac{J(\mu)}{\lambda-\mu}d\mu = 
-\D\frac{1}{\pi}{\bf Re} \D\int\limits^{\infty}_{\lambda}
\left( 1-\D\frac{P(\mu)}{R(\mu)}\right)d\mu.
\label{s3.65}
\ee
The endpoints of $B$ satisfy the compatibility condition
\be
q_{g+1} + \D\int\limits^{\infty}_{q_{g+1}}
\D\left(1-\D\frac{P(\lambda)}{R(\lambda)}d\lambda\right ) = <a_{0}>.
\label{s3.70}
\ee
\end{theorem}
\lfd
\begin{remark}
\label{rs3.5}
on condition (\ref{s3.45}). \newline
In the above Theorem \ref{ts3.5} we have not specified the value
$c_k$ which the function $J(\lambda)$ obtains in the $k$-th gap. In fact,
$c_k$ will be determined by all the other conditions of 
Theorem \ref{ts3.5}.
However, as remarked in the Introduction (\ref{e.170}), one 
observes in numerical experiments that $c_k$ should equal
$k \gamma / 2 \pi$. Figures \ref{fc7} and \ref{fc8}
in Appendix \ref{c} demonstrate that the solution $J(\lambda)$
of the integral equation (\ref{s3.40}) -- (\ref{s3.50}) indeed
satisfies this additional relation.
\end{remark}
\begin{proof}
The differentiability properties of $J(\lambda)$ in the interior 
of $B$ as well as its constancy on the $G_{k}$'s is immediately 
obvious as soon as one sees that $R$ is pure imaginary in the 
interior of $B$, and pure real elsewhere.  The function 
$J(\lambda)$ is clearly continuous at the endpoints of 
each $G_{k}$ and at $q_{g+1}$.  Also $J(p_{0})=0$ by (\ref{s3.30}) 
and consequently $J(\lambda)=0$ when $\lambda<p_0$.  
Condition (\ref{s3.40}) follows from (\ref{s3.25}), (\ref{s3.65}) 
and (\ref{s3.70}) in a straightforward way, using once again 
the pure real/pure imaginary structure of $R$.

The function $J(\lambda)$ constructed is unique. Indeed, if by 
$\delta(\lambda)$ we denote the difference of two solutions of 
(\ref{s3.40})-(\ref{s3.50}), then $\delta(\lambda)$ 
satisfies the equations
\be
\left\{
\begin{array}{rlll}
H\delta(\lambda) & =  0,            &\mbox{ when } &\lambda\in B,\\
 \delta(\lambda) & =  c^{'}_{k}=const.,  & \mbox{ when } 
&\lambda\in G_{k}, ~ n=1,\ldots g,\\
\delta(\lambda)  & =  0,             & \mbox{ when }&\lambda\not\in B\cup G,
\end{array}
\right.
\label{s3.75}
\ee
and its derivative $\delta'(\lambda)$ satisfies
\be
\left\{ 
\begin{array}{rlll}
H\delta'(\lambda)  & = 0, &\mbox{ when } & \lambda\in B,\\
  \delta'(\lambda) & = 0, &\mbox{ when } & \lambda\not\in B.
\end{array}
\right.
\label{s3.80}
\ee
 By (\ref{s3.80}) and the third relation in (\ref{s3.75}) the 
derivative $\delta'(\lambda)$ of $\delta(\lambda)$ is 
identically zero; since $\delta(\lambda)$ is compactly 
supported we also have $\delta(\lambda)=0.$
\end{proof}

The above derivation of the equation for $J(\lambda)$ has been obtained 
under a variety of assumptions.  A fully rigorous derivation still 
eludes us.  At this stage however the justification lies in the
comparison with experiments as described e.g. in Figures
\ref{fc7} and \ref{fc8} of Appendix \ref{c}.

\chapter{A boundary matching technique}
\label{m}

\section{Introduction}
\label{i}

In this chapter we construct $\frac{2 \pi}{\gamma}$ -- time periodic
solutions of
\be
\label{i2.25}
\ddot{x} _n(t) =
F(x_{n-1}(t) - x_{n}(t)) - F(x_n(t) - x_{n+1}(t)),\quad n \geq 1,
\ee
with
\be
\label{i2.26}
x_0(t) := \epsilon \sum_{m \in \ganz} b_m e^{i \gamma m t}, \quad
\sum_{m \in \ganz} |b_m| = 1,
\ee
satisfying
\be
x_n (t) = c n + O(\epsilon).
\label{i1.1}
\ee
Throughout this chapter we assume that $F$ is real analytic and 
monoton increasing on an open interval. We will construct solutions 
satisfying (\ref{i2.25})-(\ref{i1.1}) for any $c \in \relz$,
such that $-c$ lies in the interval. Not all values of $c$,
however, can be observed as the spacing of an asymptotic state
of the driven lattice, 
described by the initial boundary value problem \ref{e.90}. 
To see this we look at the Toda lattice.
For $a < 1$, in the case where $h = 0$, the solution $x_n(t)$
converges as $t \rightarrow \infty$ to $x_{{\rm asymp},\,n}(t) = c n$, with
spacing $c = -2 \ln (1 + a)$. Thus the values of the spacing $c$
that can be observed by driving the Toda lattice with constant velocity $2a$,
lie between $-\ln 4$ and $0$.
For $a > 1$, as we know, the solution of the driven lattice does
\underline{not} converge to a quiescent state $x_{{\rm asymp},\,n}(t) = c n$
and, in particular, the values $c = -2 \ln (1 + a)$, for $a > 1$
cannot be observed in this experiment.

Expand $x_n(t)$ in a Fourier series,
\be
\label{i2.37}
x_n(t) = c n + \sum_{m \in \ganz} a(n,m) e^{i \gamma m t}, \quad
\mbox{ for } n \geq 0.
\ee
For $n=0$ we have $a(0,m) = 
\epsilon b_m$. Expanding $F$ in a power series at $-c$, we obtain
\bea
&&F(x_{n-1} - x_n) - F(x_{n} - x_{n+1}) \\
&=&
F'(-c) \sum_{m \in \ganz} [a(n-1,m) -2 a(n,m) + a(n+1,m)] e^{i \gamma m t}
+ \mbox{ higher order terms }.
\eea
Equation (\ref{i2.25}) with (\ref{i2.26}) is equivalent to 
\be
\label{i2.40}
(L_m a(\cdot,m))(n) + W(a)(n,m) + a(0,m) \delta _{1,n} = 0,\quad
\mbox{ for all } m \in \ganz, n \geq 1,
\ee
where for $m \in \ganz$ the linear operators $L_m$ acting 
on the $n$-variable are given by
\be
\label{i2.41}
L_m =
\left(
\begin{array}{cccc}
\delta _m & 1 && 0 \\
1& \delta _m & 1 & \\
& 1 & \delta _m & \ddots  \\
0 && \ddots & \ddots   \\
\end{array}
\right),\quad
\delta _m := -2 + \frac{(m \gamma)^2}{F'(-c)}.
\ee
$W(a)$ contains all terms of higher order and $\delta _{1,n}$ denotes
the Kronecker symbol.

We note that for $\epsilon = 0$ equation (\ref{i2.40}) is solved by
$a=0$. However, we cannot apply the implicit function theorem to 
obtain solutions of equation (\ref{i2.40}) for $\epsilon \neq 0$,
because the linearized operator, 
\be
\label{i1.17}
L = \bigoplus_{m \in \ganz} L_m,
\ee
is not invertible.
Indeed, the spectrum of the operator $L_m$
acting on $\ell _2$ sequences is given by 
$\sigma (L_m) = [\delta _m -2, \delta _m +2]$. This implies that 
$0 \in \sigma (L_m)$ for all $m \in \ganz$ satisfying $0 \leq 
\frac{(\gamma m)^2}{F'(-c)} \leq 4$. Denote
\be
\label{i1.20}
m_0 := \max \{ m \in \ganz : 0 \leq \frac{(\gamma m)^2}{F'(-c)} \leq 4 \},
\ee
then the multiplicity of $0$ in the 
$\ell _2$ spectrum of $L$ is $2m_0 + 1$.
However, due to the simple form of the operator $L_m$ we are able to
`` invert '' the operator explicitly.
In fact, consider the important case 
$0 < \frac{(\gamma m)^2}{F'(-c)} < 4$.
For a given vector $(y_n)_{n \geq 1}$ and given
$u_1$, the vector $(u_n)_{n \geq 2}$
solves the equation
\bd
L_m u = y,
\ed
if and only if 
\bea
u_n
&=&
\frac{1}{\sin \beta _m}
\left[
u_1 \sin n \beta _m  +  \sum_{k=1}^{n-1} y_k \sin (n-k)\beta _m
\right]  \nonumber \\
&=&
\frac{1}{\sin \beta _m}
\left[
\sin n \beta _m
\left(
u_1 + \sum_{k=1}^{n-1} y_k \cos k \beta _m
\right)
- \cos n \beta _m
\sum_{k=1}^{n-1} y_k \sin k \beta _m
\right], \nonumber 
\eea
\be
\beta _m := - \mbox{ sgn }(m) \arccos \left(- \frac{\delta _m}{2} \right).
\label{i1.25}
\ee
Note, that the value of $u_n$ is independent of the sign of $\beta _m$.
We will justify the particular choice we have made in Section \ref{m3} 
(see (\ref{m3.10})) below.

The following observation will prove to be useful. Suppose 
$(y_n)_{n \geq 1}$ decays exponentially, i.e. there exists a $\sigma > 0$,
such that $\sup _{n \geq 1} \left| y_n e^{\sigma n} \right| < \infty$,
then $\sup _{n \geq 1} \left| u_n e^{\sigma n} \right| < \infty$, provided
the following two relations hold.
\be
\label{i2.45}
u_1 + \sum_{k=1}^{\infty} y_k \cos k \beta _m = 0
\ee
and
\be
\label{i2.46}
\sum_{k=1}^{\infty} y_k \sin k \beta _m = 0.
\ee
Equation (\ref{i2.45}) can always be satisfied by an appropriate choice of
$u_1$, whereas equation (\ref{i2.46}) is a condition on the sequence
$(y_n)_{n \geq 1}$. Therefore the operator $L_m$ acts 1-1 on spaces
of exponentially decaying sequences and the range has codimension 1.
Furthermore a simple calculation shows that the inverse operator
acts on the range as a bounded operator with respect to the corresponding
exponentially weighted supremum norms.
Still we cannot apply a standard implicit function theorem
to obtain solutions of equation (\ref{i2.40}). Nevertheless, proceeding
formally, we transform equation (\ref{i2.40}) into a fixed-point equation.
\be
\label{i2.50}
a(n,m) = - L_m^{-1} 
[W(a)(n,m) + \epsilon b_m \delta _{1,n}].
\ee
As $W$ is of higher order we can in principle apply a Banach fixed-point
argument to obtain a solution of (\ref{i2.50})
as long as $L_m^{-1}$ is a bounded operator. We have seen
above that this can be achieved, if condition (\ref{i2.46}) is satisfied, i.e.
\be
\label{i2.51}
\epsilon b_m \sin \beta _m +
\sum_{k=1}^{\infty} W(a)(k,m) \sin k \beta _m = 0, \mbox{ for } 
0 < \frac{(\gamma m)^2}{F'(-c)} < 4 .
\ee
Equation (\ref{i2.51}) indicates that we will be  
able to solve equation (\ref{i2.50}) for sufficiently small $\epsilon$
in a space of sequences decaying exponentially in $n$,
only if the Fourier coefficients of the driver $b_m$
take on a special value for those $m \in \ganz$ 
satisfying $0 < \frac{(\gamma m)^2}{F'(-c)} <  4$.

This observation is consistent with the linear case where we conclude from
formulae (\ref{e.115}) -- (\ref{e.125})  
that the solutions decay exponentially
in $n$,
only if $b_m = 0$ for $0 \leq \frac{(\gamma m)^2}{F'(-c)} \leq 4$. 
The linear case 
also suggests that we should add multiphase waves in order to obtain
solutions of equations (\ref{i2.25}),(\ref{i2.26})
for general driving functions. This leads to the 
following ansatz for $a(n,m)$:
\be
\label{i3.5}
\hspace{0.5cm}
a(n,m) = u(n,m) + v(n,m) + ( \epsilon b_m - u(0,m) - v(0,m) ),\quad n \geq 0,
m \in \ganz,
\ee
where $u$ denotes the travelling wave part and $v$ corresponds to the 
exponentially decaying modes. Note that (\ref{i3.5}) implies
$a(0,m) = \epsilon b_m$, for all $m \in \ganz$.

\lfd
\begin{definition}
\label{ri1.1}
We will refer to the Fourier modes $m$ with 
$0 \leq \frac{(\gamma m)^2}{F'(-c)} \leq 4$, or equivalently
$|m| \leq m_0$, as {\bf resonant Fourier modes}.
On the other hand we say that a frequency $\gamma \in \relz ^+$ is
{\bf resonant} if $\frac{(\gamma m)^2}{F'(-c)} = 4$ for some $m \in \ganz$.
\end{definition}

The present chapter is organized as follows. We begin Section \ref{m2}
by deriving equations 
for the sequences of Fourier coefficients $u(n,m)$ and $v(n,m)$ (given by
(\ref{i2.37}) and (\ref{i3.5}) above), 
which are sufficient to prove that the 
corresponding functions $x_n(t)$ solve (\ref{i2.25}), (\ref{i2.26}).   
These equations, which
are given  in Lemma \ref{lb2.1} below, 
can be made rather explicit because of the assumption, that
the force function $F$ can locally be expanded in a power series 
and therefore 
we will obtain good estimates on the higher order terms
by carefully choosing the norm on the sequences of Fourier coefficients.
In the notation of Lemma \ref{lb2.1}, these equations can be described as 
follows.
\begin{itemize}
\item[(1)]
is an equation for $u$, which is satisfied by the Fourier coefficients
of solutions 
$x_n (t;u) := \sum _{m \in \ganz} u(n,m) e^{i \gamma m t}$
of the doubly infinite lattice.
\item[(2)]
is an equation for $v$, depending on $u$, which guarantees that
$u + v$ corresponds to a solution of the semi-infinite lattice.
\item[(3)]
represents the boundary condition by requiring 
$\epsilon b_m - u(0,m) - v(0,m) = 0$, for $m \neq 0$. The case 
$m = 0$ is special; we do not have to require
$\epsilon b_0 - u(0,0) - v(0,0) = 0$ for the reason that solutions of 
(\ref{i2.25}) are invariant under translations $x_n \rightarrow x_n + const$.
\end{itemize}

We then proceed in Section \ref{m3} to prove the basic result 
(Theorem \ref{tb3.1}) of this chapter.
Assume $\gamma$ is non resonant (see Definition \ref{ri1.1}, 
then for (small) $\epsilon$ and for given
(small) travelling wave solutions $u$ of the doubly infinite lattice, we
can construct sequences $v$ satisfying equation (2) above and solving 
equation (3) for all {\em non resonant Fourier modes $m$}
(compare with Definition \ref{ri1.1}), i.e. for those $m \in \ganz$ satisfying
$|m| > m_0$.
Furthermore suppose that $u$ is given as a $C^1$ function of a parameter $q$.
Then we will show that the resulting $v$ is a $C^1$ function of $q$
and $\epsilon$. Note that this statement is needed in order to ensure 
that the remaining equations of (3) (for {\em resonant Fourier modes})
can be solved by constructing a sufficiently large parameter family
of travelling wave solutions $u(q)$, and then applying a standard 
implicit function theorem.

The proof of this basic result (Theorem \ref{tb3.1}) rests on a Banach 
fixed-point argument.
The equation for $v$ takes the form 
\be 
\label{i4.5}
(Lv)(n,m) + W(u,v)(n,m) + v(0,m) \delta _{1,n} = 0, \mbox{ for } n \geq 1,
m \in \ganz.
\ee
$L = \bigoplus _{m \in \ganzz} L_m$ was defined in 
(\ref{i2.41}) and (\ref{i1.17}) and $W(u,v)$ denotes the 
higher order terms. We turn
(\ref{i4.5}) formally into a fixed-point equation,
\be
\label{i4.10}
v(\cdot,m) = - L_m^{-1} 
\left( W(u,v)(\cdot,m) + v(0,m) \delta _{1, \cdot} \right).
\ee
Denote
\bd
S_{\sigma} := \left\{ 
(y_n)_{n \geq 1} : \sup _{n \geq 1} \left| y_n e^{\sigma n} \right|
< \infty \right\}.
\ed
As indicated above we will be able to prove the 
following results on the invertibility of $L_m$ by explicit calculation
(see proof of Theorem \ref{tb3.1}).
\begin{itemize}
\item
For $0 < |m| \leq m_0$ and $\sigma > 0$, the linear operator $L_m$ maps
$S_{\sigma}$ onto $\{ y \in S_{\sigma}: \sum _{k \geq 1} y_k \sin (k 
\beta _m) = 0 \}$. (The  quantities $\beta _m$ were defined in (\ref{i1.25}). 
The inverse
operator acting on the range is bounded with respect to the 
corresponding norms.
\item
There exist weights $\sigma > 0$, such that the operators 
$L_m: S_{\sigma} \rightarrow S_{\sigma}$ are bijective and
have a bounded inverse for all $|m| > m_0$.
\item
Again the case $m=0$ is somewhat special as the Green's function of the
operator $L_0$ grows linearly and we will have to use the special structure
of $W(u,v)$ in order to define a bounded inverse. See the proof
of Theorem \ref{tb3.1}  below for more details.
\end{itemize}

Although $u(n,m)$ does not decay in $n$, we will nevertheless see by
explicit calculation 
that $W(u,v)$ decays exponentially in $n$. This makes it possible to prove
the exitence of a solution of (\ref{i4.10}) by a Banach fixed-point
argument. For $|m| > m_0$ we can choose $v(0,m) = \epsilon b_m - u(0,m)$ 
and hence
satisfy the boundary condition as described in equation (3) above,
whereas in the case $0 < |m| \leq m_0$ the choice of $v(0,m)$ is determined
by the condition that $W(u,v)(\cdot,m) + v(0,m) \delta _{1, \cdot}$ has 
to lie in the range of $L_m$.
A small technical problem arises when proving the smooth dependence of
$v$ on the parameters. It will turn out 
that the travelling wave solutions constructed in the subsequent 
chapters depend smoothly on $q$, but
$\frac{\partial u}{\partial q} (n,m)$ grows linearly in $n$.
Therefore $\frac{\partial u}{\partial q}$ does not lie in a
space which is suitable for our calculations. We will verify the 
smooth dependence of $v$ on $q$ and $\epsilon$ explicitly
by applying a Banach fixed-point
argument to the partial derivatives in the appropriate spaces.

In Section \ref{m4} we show that the results of Sections \ref{m2}
and \ref{m3} suffice to construct periodic solutions
of the driven lattice in the case that $m_0 = 0$, which
corresponds to high driving frequencies.

\section{The equation for the Fourier coefficients}
\label{m2}

In this section we introduce norms, which are suitable for the sequences 
of Fourier coefficients, and prove some of their basic properties.
Then the {\em general assumptions} on the force function $F$
and on the driver will be stated precisely. Using the assumptions on $F$
we derive estimates
on the nonlinear terms which allow us to give conditions on
the Fourier coefficients which are sufficient for proving that the 
correponding functions $x_n(t)$ given by (\ref{i2.37}) and 
(\ref{i3.5}) solve equations (\ref{i2.25}) and (\ref{i2.26}).

\subsection{Sequence spaces}
\label{m21}

The choice for the norms on the sequences of Fourier coefficients $u(n,m)$
and $v(n,m)$ (compare with equation (\ref{i3.5})) 
is motivated by the following observations. 
The nonlinear terms of the force function
make it necessary to take convolutions with respect to the $m$-variable
(see Section \ref{m23} below). Therefore
we choose an $\ell _1$-norm for $m$. In fact we use a
weighted $\ell _1$-norm in order to control the regularity of the solution.
Furthermore the weight function has to satisfy some additional conditions 
to insure that the norm is still compatible with respect to 
convolution (see Definition \ref{di2.1} of {\em admissible weight functions}). 
For the $n$-variable a supremum norm with an exponential weight is 
chosen which is suitable for inverting the linearized operators $L_m$.
\lfd
\begin{definition}
\label{di2.1}
A map $w:\ganz \rightarrow \relz$ is said to be an 
\underline{admissible weight function},
if
\be
\label{i2.30}
w(m) \geq 1 , \quad \mbox{ for all } m \in \ganz,
\ee
and
\be
\label{i2.31}
w(m) \leq w(m-n) w(n),\quad \mbox{ for all } m,n \in \ganz.
\ee
\end{definition}

\lfd
\begin{definition}
\label{di2.2}
Let $w$ be an {\em admissible weight function}. We denote 
\bd
\ell _{1,w} :=\left\{ u:\ganz \rightarrow \comz \left|
                  \sum_{m \in \ganz} w(m) |u(m)| < \infty \right. \right\},
\ed
with the corresponding norm
\bd
           \| u \| _{\ell _{1,w}} := \sum_{m \in \ganz} w(m) |u(m)|.
\ed
\end{definition}
Note that $\ell _{1,w}$ is a Banach space which lies in $\ell _1$ by
condition (\ref{i2.30}). The inequality (\ref{i2.31}) insures that
the $\ell _{1,w}$-norm is submultiplicative with respect to convolution
(see Proposition \ref{pb2.2} below). 
It also implies that the weight function
can not grow faster than exponentially. Indeed, it is easy to prove
that $w(m) \leq w(0) \left(\max (w(1),w(-1)) \right) ^{|m|}$.
We shall be interested in three types of weight functions
which will all satisfy the conditions specified in Definition \ref{di2.1}.
\begin{itemize}
  \item[(i)]   $ \forall m \in \ganz : w(m) := 1$ .
  \item[(ii)]  $ \forall m \in \ganz : w(m) := (1 + |m|)^{\beta} $, for 
    $ \beta \geq 0$ .
  \item[(iii)] $ \forall m \in \ganz : w(m) := e^{\beta |m|}$,
    for $ \beta \geq 0$.
\end{itemize}
Finally note that the product of two admissible weight functions is again an
admissible weight function. 
\lfd
\begin{definition}
\label{db2.1}
Let $w$ be an admissible weight function and 
let $ \sigma \in \relz$. Then
\bd
{\cal L}_{\sigma,w} :=
             \left\{u: \natz _{0} \times \ganz \rightarrow \comz  \left|
             \left( \sup_{n \geq 0} e^{\sigma n} |u(n, \cdot)| \right)_
                  {m \in \ganz}  \in \ell _{1,w}  \right. \right\} ,
\ed
with the corresponding norm 
\bd
     \| u \| _{\sigma,w} := \sum_{m \in \ganz} w(m)      
                                   \sup_{n \geq 0} e^{\sigma n} |u(n,m)| .
\ed
\end{definition}
It is easy to check that ${\cal L} _{\sigma ,w}$ are Banach spaces 
and that ${\cal L} _{\sigma _1,w} \subset {\cal L} _{\sigma _2,w}$ for
$\sigma _1 \geq \sigma _2$.
In the following proposition we recall some simple
properties of the convolution of sequences, which is defined by
$(u*v)(m):=\sum_{l \in \ganz} u(m-l) v(l)$. Furthermore we provide 
estimates on the convolution in terms of the norms defined above.

\lfd
\begin{proposition}
\label{pb2.2}
Let $c \in \comz , u_1, u_2, u_3, u_4 \in \ell _{1}$, then
  \begin {itemize}
  \item[(i)]  $u_1 * u_2 \in \ell _{1} 
  \mbox{ and } \| u_1 * u_2 \| _{\ell _{1}}
                      \leq \| u_1 \| _{\ell _{1}}  \| u_2 \| _{\ell _{1}}.$
  \item[(ii)] $u_1 * u_2 = u_2 * u_1 .$
  \item[(iii)] $u_1*(u_2*u_3) = (u_1*u_2)*u_3 .$
  \item[(iv)] $u_1*(u_2 + c u_3) = u_1 * u_2 + c (u_1 * u_3) .$
  \item[(v)] If $u_1$ and $u_2$ satisfy the reality condition (i.e.
    $ \forall m \in \ganz : u_i(m) = \overline{u_i(-m)} ; i=1,2$), then 
    $u_1 * u_2$ also satisfies the reality condition.
  \item[(vi)]  If for all $m \in \ganz : |u_1(m)| \leq u_2(m)$ and
     $|u_3(m)| \leq u_4(m)$, then $|(u_1*u_3)(m)| \leq (u_2*u_4)(m)$          
     for all $m \in \ganz$ .
  \item[(vii)] Convolution 
    respects the $\ell _{1,w}$ norm, i.e. let
    $v_1,v_2 \in \ell _{1,w}$ then
    $v_1 * v_2 \in \ell _{1,w}$ and $\| v_1*v_2 \| _{\ell _{1,w}} \leq 
    \| v_1 \| _{\ell _{1,w}}  \| v_2 \| _{\ell _{1,w}}$ .
  \item[(viii)]
    Let $\sigma _1, \sigma _2 \in \relz$ and $u \in {\cal L} 
    _{\sigma _1 ,w}, v \in {\cal L} _{\sigma _2 ,w}$ and define their $m$ -
    convolution $y$ by $y(n,m) := \sum_{l \in \ganz}
    u(n,m-l) v(n,l)$. Then
    $y \in {\cal L} _{\sigma _1 + \sigma _2,w}$
    and 
    $ \| y \| _{\sigma _1 + \sigma _2,w}  \leq
    \| u \| _{\sigma _1,w}  \| v \| _{\sigma _2,w}.$
  \end{itemize}
\end{proposition}

\begin{proof}
Properties (i)-(vi) are standard. In order to show (vii) we note that
the inequality (\ref{i2.31}) implies
\bea
\sum_{m} w(m) |(u_1 * u_2) (m)| &\leq& 
 \sum_{m,n} w(m-n) |u_1(m-n)| w(n)|u_2(n)| \\
&=&
 \| (w|u_1|)*(w|u_2|) \| _{\ell _1}, 
\eea
and by (i) this is all we need. 
Property (viii) is a consequence of (vi) and (vii). 
\end{proof}

\subsection{The general assumptions}
\label{m22}

Recall the notation which was introduced in equations (\ref{i2.25})
and (\ref{i2.26}). We now state the assumptions on the force function
$F$, the frequency $\gamma$ and the Fourier coefficients 
$(b_m)_{m \in \ganz}$ of the driver.

{\bf The general assumptions.}
\begin{itemize}
\item[(1)]
$F: \relz \rightarrow \relz$ is real analytic 
in a neighborhood of $-c, c \in \relz$, and
\be
\label{i2.35}
F'(-c) > 0.
\ee
\item[(2)]
$\gamma \in \relz ^+ \setminus 
\{\gamma: \frac{(m \gamma)^2}{F'(-c)} = 4$ for some $m \in \ganz \}$. 
\item[(3)]
$(b_m)_{m \in \ganz} \in \ell _{1,w}$ for some 
{\em admissible weight function} $w$
and $\| b_m \| _{\ell _{1,w}} = 1$.
\end{itemize}
{\bf Remark:}
\begin{itemize}
\item
We are looking for a solution of the type $x_n(t) = c n + O(\epsilon)$.
Therefore $F(x_{n-1} - x_n) = F(-c + O(\epsilon))$. Condition (1)
will allow us to expand $F(x_{n-1} - x_n) - F(x_{n} - x_{n+1})$
in a power series where the linear term does not vanish.
\item
It will be shown that the exceptional set of resonant frequencies (see 
Definition \ref{ri1.1})
$\{\gamma: \frac{(m \gamma)^2}{F'(-c)} = 4$ for some $m \in \ganz \}$
consists
of precisely those frequencies for which the number of phases in the 
travelling wave solution described above
changes. In the case of the Toda lattice these are also the frequencies
for which the number of gaps in the spectrum of the corresponding 
Lax operator at $t= \infty$ changes.
\item
It turns out that the weighted spaces $\ell _{1,w}$ are well suited 
to proving that the regularity of the solution is comparable to the 
regularity of the driver.
\end{itemize}

\subsection{The nonlinear terms}
\label{m23}

The force function $F$
is assumed to be a real analytic function at $-c$ (see {\em
general assumptions} above) and we can define
for all $k \geq 0$,
\be
\label{b2.10}
\alpha _{k} := \frac{\partial ^{k}}{\partial x^{k}} F(-c).
\ee
By $\rho _{F,c}$ we denote the minimum of $1$ and the radius of 
convergence of the power series 
$\sum_{k=0}^{\infty} \frac{\alpha _{k}}{k!} (x+c)^{k}$.
Recall that $\alpha _1 \neq 0$ by the {\em general assumptions}. 
Therefore  we obtain the
following estimates by standard arguments for power series.

\lfd
\begin{proposition}
\label{pb2.1}
There exists a constant $\tilde{C}_{F,c}$, such that for all $y$,
$|y| \leq \frac{\rhof}{2}$,
\bea
\frac{1}{|\alpha _{1}|} \sum_{k=2}^{\infty} \frac{|\alpha _{k}|}{k!} |y|^k 
&\leq& \tilde{C}_{F,c} |y|^2.  \\
\frac{1}{|\alpha _{1}|} \sum_{k=2}^{\infty} 
\frac{|\alpha _{k}|}{(k-1)!} |y|^{k-1} 
&\leq& \tilde{C}_{F,c} |y|.  \\
\frac{1}{|\alpha _{1}|} \sum_{k=2}^{\infty} 
\frac{|\alpha _{k}|}{(k-2)!} |y|^{k-2} 
&\leq& \tilde{C}_{F,c}.  
\eea
\end{proposition}

We now define the higher order terms of the equations for the
Fourier coefficients as formal power series.
The convergence of these series and various
differentiability properties will be discussed in
the subsequent proposition. In order to see that the following 
expressions indeed represent the higher order terms of the equation,
one may look at Lemma \ref{lb2.1} below.

If $u = u(n,m) \in {\cal L} _{\sigma,w}$, we use $u(n,\cdot)^{*k}$ to 
denote the $k$ -th $m$-convolution of $u$, that is
\be
\label{b2.13}
u(n,\cdot)^{*k}(m) = 
\sum_{l_1 + \cdots + l_k = m} u(n,l_1) \cdot \ldots \cdot u(n,l_k).
\ee
\lfd
\begin{definition}
\label{db2.2}
For $\sigma \geq 0, u,v \in {\cal L} _{\sigma,w}$ and $n \geq 1$, denote 
\begin{eqnarray*}
\triangle u(n,m) &:=&
    u(n-1,m) - u(n,m) . \\
W(u)(n,m) &:=& \frac{1}{\alpha _1} \summb \frac{\alpha _k}{k!}
    \left( (\triangle u)(n,\cdot)^{*k} 
    - (\triangle u)(n+1,\cdot)^{*k} \right)(m)          
    . \\
Y(u,v)(n,m) &:=& \frac{1}{\alpha _1} \summb \frac{\alpha _k}{k!}
    \left( (\triangle (u+v))(n,\cdot)^{*k} 
    - (\triangle u)(n,\cdot)^{*k} \right)(m)
    . \\
W(u,v)(n,m) &:=& Y(u,v)(n,m) - Y(u,v)(n+1,m) . 
\end{eqnarray*}
For $n = 0$ and for all $m \in \ganz$ set 
$\triangle u(0,m) := W(u)(0,m) := Y(u,v)(0,m) := W(u,v)(0,m) := 0$.
\end{definition}

\lfd
\begin{proposition}
\label{pb2.3}
There exists a constant $C_{F,c}$, such that
for all $0 \leq \sigma \leq 1,
u \in {\cal L} _{0,w}$ and $v \in {\cal L} _{\sigma ,w}$ with
$\| u\| _{0,w}, \| v \| _{\sigma ,w} < \frac{\rhof}{8}$ 
the following is true. 
The series in the definition of $W$ and $Y$ converge absolutely with 
$ W(u) \in {\cal L} _{0,w}, Y(u,v) \in {\cal L} _{\sigma ,w}$. Furthermore
\begin{itemize}
  \item[(i)] $\| W(u) \| _{0,w} \leq C_{F,c} \| u \| _{0,w}^2$.
  \item[(ii)] $\| Y(u,v) \| _{\sigma ,w}      \leq  C_{F,c}
    \| v \| _{\sigma ,w} \max ( \| u \| _{0,w}, \| v \| _{0,w} )$.
  \item[(iii)]
    The map $F_1: 
    \{ v \in {\cal L} _{\sigma ,w}: \| v \| _{\sigma ,w} < \frac{\rhof}{8} \}
    \rightarrow {\cal L} _{\sigma ,w}
    , v \mapsto Y(u,v) $ is $C^2$ and the derivatives satisfy the following
    estimates
    \bea
      \forall x \in {\cal L} _{\sigma ,w}  &:&
      \| (D_v Y)(u,v) x \| _{\sigma ,w}  \leq
      C_{F,c} \max ( \| u \| _{0,w}, \| v \| _{0,w} ) 
      \| x \| _{\sigma ,w}. \\
      \forall x_1, x_2 \in {\cal L} _{\sigma ,w}  &:&
      \| (D_v^2 Y)(u,v) [x_1, x_2 ] \| _{\sigma ,w}  \leq
      C_{F,c}  \| x_1 \| _{\sigma ,w} \| x_2 \| _{\sigma ,w} . 
    \eea
  \item[(iv)]
    The map $F_2: 
    \{ u \in {\cal L} _{0,w}: \| u \| _{0,w} < \frac{\rhof}{8} \}
    \rightarrow {\cal L} _{\sigma ,w}
    , u \mapsto Y(u,v) $ is $C^1$ with derivative
    \be
    \label{b2.30}
    \hspace{.4in}
      D_u Y (u,v) x =
      \frac{1}{\alpha _1} \summb \frac{\alpha _k}{k!} 
      \sum_{l=1}^{k-1} 
      \left(    \begin{array}{l} k \\ l \end{array}   \right)
      (k-l)
      (\triangle u)^{*(k-l-1)} *
      (\triangle v)^{*l} * \triangle x.
    \ee
    $D_u Y(u,v)$ as given in equation (\ref{b2.30}) can be regarded
    as a bounded linear operator from 
    ${\cal L} _{\sigma ',w}$ into 
    ${\cal L} _{\sigma + \sigma ',w}$ for $\sigma ' \in \relz$ and the
    corresponding operator norm is bounded by
    $(1 + e^{\sigma '}) C_{F,c}
    \max ( \| u \| _{0,w}, \| v \| _{\sigma, w} )$.
  \item[(v)]
    Let $\sigma ' \geq 0$ and
    fix $x \in {\cal L} _{- \sigma ',w}$. The map \newline
    $F_3 : {\cal L} _{\sigma ,w} \rightarrow 
    {\cal L} _{\sigma - \sigma ',w}, v \mapsto (D_u Y)(u,v) x$
    is $C^1$  and the derivative satisfies the estimate \newline
    $ \forall z \in {\cal L} _{\sigma ,w} : 
    \| D_v F_3 (v) z \| _{\sigma - \sigma ',w} \leq
    C_{F,c} \| x \| _{- \sigma ',w} \| z \| _{\sigma ,w}.$
\end{itemize}
\end{proposition}

{\bf Remark:}
The differentiability properties (iii)-(v) will not be used
in the present section, but they are needed in Section \ref{m3}
when we prove differentiability of the solution of the fixed-point
equation with respect to certain parameters (compare with the proof of
Theorem \ref{tb3.1} ).

\begin{proof}
(i) We begin by remarking that $u \in \raua$ implies 
$\triangle u \in \raua$ and 
$ \| \triangle u \nora  \leq  2 \| u \nora$.
By Proposition \ref{pb2.2} (viii) and Proposition \ref{pb2.1} it is easy
to see, that
\begin{eqnarray*}
\| W(u) \nora & \leq &
  2 \frac{1}{| \alpha _1 |} \summb \frac{| \alpha _k |}{k!}
  \| \triangle u  \| _{0,w}^k \\
  & \leq &
  2  \tilde{C}_{F,c} (2 \| u \nora)^2.
\end{eqnarray*}

(ii) In this case one has to evaluate
\bd
(\triangle (u+v))(n,\cdot)^{*k} - (\triangle u)(n,\cdot)^{*k}  =
  \sum_{l=1}^{k}
  \left(    \begin{array}{l} k \\ l \end{array}   \right)
  (\triangle v)(n,\cdot)^{*l} * (\triangle u)(n,\cdot)^{*(k-l)}.
\ed
Using again Proposition \ref{pb2.2} (viii) and Proposition \ref{pb2.1}
we obtain
\begin{eqnarray*}
  \| Y(u,v) \norb & \leq &
  \frac{1}{| \alpha _1 |} \summb \frac{| \alpha _k |}{k!}
  \sum_{l=1}^{k} 2^k
  \| \triangle v \nora ^{l-1}
  \| \triangle u \nora ^{k-l}
  \| \triangle v \norb  \\
  & \leq &
  \frac{1}{| \alpha _1 |} \summb \frac{| \alpha _k |}{(k-1)!}
  \left(2 \max ( \| \triangle u \nora ,\| \triangle v \nora) \right) ^{k-1}
  2 \| \triangle v \norb  \\
  & \leq &
  8 \tilde{C}_{F,c} \max ( \| u \nora ,\| v \nora)  \| \triangle v \norb.
\end{eqnarray*}
Observing that $\| \triangle v \norb \leq (1+e^{\sigma}) \| v \norb $ 
the claim follows.\newline

(iii) The proof of differentiability for $F_1$ (as well as $F_2$ and $F_3$)
uses the fact that these functions are sums over $l$ and $k, l \leq k$,
of monomials of the form $(\triangle v)(n, \cdot)^{*l} *
(\triangle u)(n, \cdot)^{*(k-l)}$. Therefore it
suffices to first prove the continuous 
differentiability of each term in the sum 
and to show secondly that the sum of the derivatives converges uniformly
in the corresponding norm.
\newline
Because of the simple algebraic rules for convolution (see Proposition
\ref{pb2.2}) it is straightforward to check that for $l \geq 1$ the map
\bd
F_4:v \mapsto (\triangle v)(n, \cdot)^{*l} * 
(\triangle u)(n, \cdot)^{*(k-l)} 
\ed
is a $C^1$ map from 
${\cal L}_{\sigma ,w}$ into ${\cal L}_{\sigma ,w}$ with derivative
\bd
DF_4(v) x = l (\triangle v)(n, \cdot)^{*(l-1)} *
(\triangle u)(n, \cdot)^{*(k-l)} * (\triangle x)(n, \cdot).
\ed
Proposition \ref{pb2.2} (viii) yields the estimate in the corresponding
operator norm
\bd
\| DF_4 (v) \| \leq
(1+e^{\sigma}) l \left(
\max (\| \triangle u \| _{0,w},\| \triangle v \| _{0,w}) \right)^{k-1}
\ed
and with Proposition \ref{pb2.1} we conclude the uniform convergence of
the sum, as
\bea
&&
  \frac{1}{| \alpha _1 |} \summb \frac{| \alpha _k |}{k!} 
  \sum_{l=1}^{k} 
  \left(    \begin{array}{l} k \\ l \end{array}   \right)
  (1+e^{\sigma})l  \left( 2
  \max (\| u \| _{0,w},\| v \| _{0,w}) \right)^{k-1}  \\
& \leq &
  \frac{1}{| \alpha _1 |} \summb \frac{| \alpha _k |}{(k-1)!} 
  (1+e^{\sigma}) 2  \left( 4
  \max (\| u \| _{0,w},\| v \| _{0,w}) \right)^{k-1}  \\
& \leq &  
  2 (1+e^{\sigma}) \tilde{C}_{F,c} 4
  \max (\| u \| _{0,w},\| v \| _{0,w}). 
\eea
This proves everything about the first derivative. For the second 
derivative we can proceed similarily. We get
\bd
D^2F_4(v)[x,y] = l (l-1) 
(\triangle v)(n, \cdot)^{*(l-2)} *
(\triangle u)(n, \cdot)^{*(k-l)} * (\triangle x)(n, \cdot)
* (\triangle y)(n, \cdot).
\ed
The convergence of the sum is guaranteed by
\bea
&&
  \frac{1}{| \alpha _1 |} \summb \frac{| \alpha _k |}{k!} 
  \sum_{l=1}^{k} 
  \left(    \begin{array}{l} k \\ l \end{array}   \right)
  (1+e^{\sigma})^2 l(l-1)  \left( 2
  \max (\| u \| _{0,w},\| v \| _{0,w}) \right)^{k-2}    \\
& \leq &
  \frac{1}{| \alpha _1 |} \summb \frac{| \alpha _k |}{(k-2)!} 
  (1+e^{\sigma})^2 4  \left( 4
  \max (\| u \| _{0,w},\| v \| _{0,w}) \right)^{k-2}      \\
& \leq &  
  4 (1+e^{\sigma})^2 \tilde{C}_{F,c}. 
\eea

(iv)  
The proof is rather similar to the one just given. For $k \geq 2, l \geq 1$
let
\bd
F_5: u \mapsto 
(\triangle u)^{*(k-l)} * (\triangle v) ^{*l}.
\ed
$F_5$ is a $C^1$ map from ${\cal L} _{0,w}$ into ${\cal L} _{\sigma,w}$
with
\bea
DF_5 (u) x &=&
(k-l) (\triangle u)^{*(k-l-1)} * (\triangle v) ^{*l}
* \triangle x . \\
\| DF_5 (u) x \| _{\sigma ,w}  &\leq&
(k-l) \left(
2 \max (\| u \| _{0,w} , \| v \| _{0,w} ) \right) ^{k-2}
(1 + e^{\sigma}) \| v \| _{\sigma ,w} 2 \| x \| _{0,w}.
\eea
Proposition \ref{pb2.1} gives the uniform convergence of the sum.
The remaining part of (iv) can be easily seen from 
Proposition \ref{pb2.1}, Proposition \ref{pb2.2} (viii) and the just given
formula.

(v) Applying the procedure again, we first convince ourselves that 
for $l \geq 1$ the function 
\bd
F_6: v \mapsto 
(\triangle u)^{*(k-l-1)} * (\triangle v) ^{*l}
* \triangle x
\ed
is a $C^1$ map from ${\cal L} _{\sigma,w}$ into
${\cal L} _{\sigma - \sigma ',w}$ with derivative
\bd
DF_6(v) z =
l (\triangle u)^{*(k-l-1)} * (\triangle v) ^{*(l-1)}
* \triangle x * \triangle z.
\ed
The sum of the operator norms of the derivatives can uniformly be 
estimated by
\bea
&&  
  \frac{1}{| \alpha _1 |} \summb \frac{| \alpha _k |}{(k-2)!} 
  2^k  \left( 2
  \max (\| u \| _{0,w},\| v \| _{0,w}) \right)^{k-2}      
  2 \| x \| _{- \sigma ',w} (1+e^{\sigma}) \\
& \leq &  
  8 (1+e^{\sigma}) \tilde{C}_{F,c} \| x \| _{- \sigma ',w}. 
\eea
This concludes the proof of the proposition.
\end{proof}

\subsection{The equations for the Fourier coefficients}
\label{m24}

Following the ansatz described in Section \ref{i}, we are now ready to 
give sufficient
conditions for the Fourier coefficients  $u(n,m)$ and $v(n,m)$ 
in order to obtain real
solutions for the driven nonlinear lattice described by equations
(\ref{i2.25}) and (\ref{i2.26}). Recall from equation (\ref{i2.41})
the definition $\delta _m = -2 + \frac{(m \gamma)^2}{F'(-c)}.$

\lfd
\begin{lemma}
\label{lb2.1}
Let $F,c,\gamma, (b_m)_{m \in \ganz}, w$ 
satisfy the general assumptions.
Suppose there exist $ u,v \in \raua$, for which the following conditions
hold.
\begin{itemize}
  \item[(1)] 
    \begin{eqnarray*}
      \forall n \geq 1,m \in \ganz: 
&u(n-1,m) + \delta _m u(n,m) + u(n+1,m) + W(u)(n,m) =0.& \\
      \forall n \geq 0,m \in \ganz:
&u(n,-m) = \overline{u(n,m)}.& \\
      \| u \nora < \frac{\rhof}{8}. &&
    \end{eqnarray*}
  \item[(2)] 
    \begin{eqnarray*}
      \forall n \geq 1,m \in \ganz:
&v(n-1,m) + \delta _m v(n,m) + v(n+1,m) + W(u,v)(n,m) =0.& \\
      \forall n \geq 0,m \in \ganz: 
&v(n,-m) = \overline{v(n,m)}.& \\
      \| v \nora < \frac{\rhof}{8} .&&
    \end{eqnarray*}
  \item[(3)]
    $\forall m \neq 0 : \epsilon b_m - v(0,m) - u(0,m) = 0.$
\end{itemize}
Then the family of real valued and periodic functions 
\bd
x_n(t) := cn + \sum_{m \in \ganz} a(n,m)  e^{i m \gamma t}
,\quad \mbox{ for } n \geq 1, 
\ed
with 
\bd
a(n,m) := u(n,m) + v(n,m) + \epsilon b_m - u(0,m) - v(0,m) , 
\quad \mbox{ for } n \geq 0,
\ed
solves the equations 
\bd
\ddot{x}_n  = F(x_{n-1} - x_n) - F(x_n - x_{n+1}),\quad n \geq 1,
\ed
where 
\bd
x_0(t) = \epsilon \sum_{m \in \ganz} b_m  e^{i m \gamma t} .
\ed
\end{lemma}

Note that the sequence $a$ is also defined for $n=0$ and that
$a(0,m) = \epsilon b_m$, the Fourier coefficients of the driver.

\begin{proof}
First we note that the $x_n$ are twice differentiable functions for $n \geq 1$. 
In fact, we know that $(\delta _m u(n,m))_{m \in \ganz}$ and   
$(\delta _m v(n,m))_{m \in \ganz}$ are sequences in $\ell _1$ , as they can be 
expressed in terms of the $\ell _1$ sequences 
$u(n-1,\cdot) ,u(n+1,\cdot), W(u)(n,\cdot)$ and 
$v(n-1,\cdot) ,v(n+1,\cdot), W(u,v)(n,\cdot)$ . But this implies
that $(m^2 a(n,m))_{m \in \ganz}$ is in $\ell _1$, which yields the $C^2$
regularity of $x_n$. Furthermore it is immediate that all functions $x_n$ 
are real valued and periodic with period $\frac{2 \pi}{\gamma}$.

Let us now turn to the main point of the proof, namely to verify
that the $x_n$ are solutions of the driven lattice. As $\ddot{x} _n$
and $ F(x_{n-1} - x_n) - F(x_n - x_{n+1})$ are both continuous functions
of the same period, it suffices to show that their Fourier coefficients
coincide. One checks from the definitions  that
\begin{eqnarray*}
\forall n \geq 1 :
x_{n-1} (t) - x_n (t) &=&
-c + \sum_{m \in \ganz} \triangle a(n,m) e^{i m \gamma t}   \\
&=& -c + \sum_{m \in \ganz} \triangle (u+v)(n,m) e^{i m \gamma t}.
\end{eqnarray*}
Substituting into the Taylor series for $F$ yields
\begin{eqnarray*}
F(x_{n-1} (t) - x_n (t)) & = &
\sum_{k \geq 0} \frac{\alpha _k}{k!}   \left(
  \sum_{m \in \ganz} \triangle (u+v)(n,m) e^{i m \gamma t}  \right) ^k   \\
&=&
\sum_{k \geq 0} \frac{\alpha _k}{k!}   
  \sum_{m \in \ganz} (\triangle (u+v))(n,\cdot)^{*k} (m) e^{i m \gamma t} \\
&=&
\sum_{m \in \ganz}   \left(
  \sum_{k \geq 0} \frac{\alpha _k}{k!}
  (\triangle (u+v))(n,\cdot)^{*k} (m)   \right)   e^{i m \gamma t},
\end{eqnarray*}
where all the manipulations are justified as the sums converge absolutely
(compare with proof of Proposition \ref{pb2.3}). We can now read off the 
Fourier coefficients.
\begin{eqnarray*}
&& \frac{ \gamma }{2 \pi} \int_{0}^{\frac{2 \pi }{ \gamma }}
\left[  F(x_{n-1}(t) - x_n(t)) - F(x_n(t) - x_{n+1}(t))  
\right] e^{- i m \gamma t} dt \\
&=&
\alpha _1  \left[
\triangle (u+v) (n,m) - \triangle (u+v) (n+1,m) + W(u)(n,m) + W(u,v)(n,m)
\right].
\end{eqnarray*}
On the other hand, using condition (3) of the hypothesis
\begin{eqnarray*}
  \frac{ \gamma }{2 \pi} \int_{0}^{\frac{2 \pi }{ \gamma }}
  \ddot{x} _n (t) e^{- i m \gamma t} dt  
  &=&
  -( \gamma  m)^2 a(n,m) \\
  &=&
  -( \gamma  m)^2  (u+v)(n,m) .
\end{eqnarray*}
The equality of the Fourier coefficients follows from (1) and (2) of
the hypothesis.
\end{proof}

It is a simple corollary of the proof of the Lemma, to see that 
condition (1) is satisfied if we have a ``small'' solution of the 
doubly infinite lattice equation. This is stated more precisely 
in the following remark.
 
\lfd
\begin{remark}
\label{rb2.1}
Suppose that $u : \ganz \times \ganz \rightarrow \comz$ satisfies
$\sum_{m \in \ganz} \sup_{n \in \ganz} |u(n,m)| < \frac{\rhof}{8}$.
For $ n \in \ganz$ set 
\bd
x_n^{(0)} (t) :=  c n + \sum_{m \in \ganz} u(n,m) e^{i m \gamma t}.
\ed
Then $u$ satisfies condition 
\begin{itemize}
  \item[(1')] 
    \bea
      \forall n \in \ganz, m \in \ganz & : &
      u(n-1,m) + \delta _m u(n,m) + u(n+1,m) + W(u)(n,m) =0, \\
      \forall n \in \ganz, m \in \ganz & : &
      u(n,-m) = \overline{u(n,m)}, 
    \eea
\end{itemize}
if and only if $x_n^{(0)} (t)$ is a real valued solution of the differential
equation
\bd
\ddot{x} _n (t) =
F(x_{n-1} (t) -x_{n} (t)) -
F(x_{n} (t) -x_{n+1} (t)), \mbox{ for } n \in \ganz.
\ed
\end{remark}

\section{Solving for the non resonant modes}
\label{m3}

The present section is devoted to the proof of our basic result, which was
explained and motivated in Chapter \ref{e} and in the introduction to this 
chapter.
Before the theorem can be stated we recall some notation of Section \ref{i}
and we add a few definitions. 

In equation (\ref{i2.41}) we have set
\bd 
\delta _m = -2 + \frac{(m \gamma)^2}{F'(-c)}.
\ed
Furthermore we denoted in equation (\ref{i1.20}),
\bd 
m_0 = \max \{m \in \ganz : 0 \leq \frac{(m \gamma)^2}{F'(-c)} \leq 4 \}.
\ed
By separation of variables one obtaines solutions of the free linearized
problem of the form
\bd
y_n(t) = z_m^n e^{i \gamma m t}, \quad m \in \ganz,
\ed
where
\be
z_m^2 + \delta _m z_m + 1 = 0. 
\label{m3.3}
\ee
The case $|m| > m_0$ corresponds to $\delta _m > 2$ and therefore we can
pick $z_m$ to be the solution of the above equation 
with $|z_m| < 1$, which is given by
\be
\forall |m| > m_0 : z_m  := 
       - \frac{ \delta _{m}}{2} +
             \sqrt{ \frac{ \delta _m ^2}{4} -1} , \quad
       \mbox{ for } \delta _m > 2.  
\label{m3.5}
\ee

In the case $0 < |m| \leq m_0$, the {\em general assumptions} on the
frequency $\gamma$ imply that $|\delta _m| < 2$. We choose for $z_m$
the solution of equation (\ref{m3.3}), 
which corresponds to an outgoing wave $y_n(t)$,
i.e. $z_m = e^{i \beta _m}$ with
\be
\beta _m = - \mbox{ sgn }(m) \arccos \left( - \frac{\delta _m}{2} \right).
\label{m3.10}
\ee
This explains the choice of the sign in equation (\ref{i1.25}).
Note, that for all $0 < |m| \leq m_0$ we have
\be
| \sin \beta _m | > 0,
\label{m3.15}
\ee
as $|\delta _m| < 2$.

\lfd
\begin{definition}
\label{db3.1}
Let $\gamma ,F,c, (b_m)_{m \in \ganz}, w$ satisfy the 
{\em general assumptions}.
\begin{eqnarray*}
\sigma _0& := &
\min (1, - \frac{1}{2} \ln (| z_{m_0 +1}|) ). \\
\sigma _1& := &
\frac{1}{4} \sigma _0 . \\
C_K & := &
\max \left[ \frac{1 + e^{-2 \sigma _0}}{ 1 - e^{-4 \sigma _0}}
  \left( \frac{1}{e^{\sigma _0 }-1} + 
  \frac{1}{1- e^{-2 \sigma _0}} \right),\right. \\
  &&\left. \frac{2}{1-e^{- \sigma _1}} \max \left\{ 
  \frac{1}{|\sin \beta _m|} : 0<m \leq m_0 \right\} \right]. \\
\rho _{u} & := &
\sum_{|m| > m_0} w(m)|u(0,m)|.
\end{eqnarray*}
\end{definition}
 
{\bf Remarks :}      
\begin{itemize}
\item
The constant $C_K$ is well defined, as $\sigma _0, \sigma _1 > 0$ and
$|sin \beta _m| > 0$ for $0 < |m| \leq m_0$. We will see below,
that $C_K$ is an upper bound on a linear operator $K$ which is related
to the to the inverse of the operator 
$\bigoplus _{m \in \ganz} L_m$ (see  
Proposition \ref{pb3.1} below).
\item
The use of $\rho _u$ will not become clear before Chapter \ref{l}.
There we will see that this quantity is the main tool for proving
that the sequence $v$ ( in the notation of Lemma \ref{lb2.1}),
which will be constructed in the next theorem,
is of higher order.
\end{itemize}

The following is the main result of this chapter.

\lfd
\begin{theorem}
\label{tb3.1}
Let $\gamma ,F,c, (b_m)_{m \in \ganz},w$ satisfy the general assumptions.
Furthermore we assume that there exists a choice of 
constants $N \in \natz, c_0 > 0, C_0 > 1$,
and a map $\{q \in \relz ^N: |q| \leq c_0 \} \rightarrow 
{\cal L} _{0,w}:
q \mapsto u(q)$,  such that  
\begin{itemize}
\item[(1)]
$\forall |q| \leq c_0 : \| u(q) \| _{0,w} \leq c_1$ and 
$\rho _{u(q)} \leq \frac{c_1}{4}$, where
\be
\label{b3.5}
c_1: = \min \left( \frac{1}{4 C_K C_{F,c} }, \frac{\rhof}{8} \right).
\ee
\item[(2)]
$q \mapsto u(q)$ is a $C^1$ map from $\{q \in \relz ^N: |q| \leq c_0 \}$ to 
${\cal L} _{-\sigma _1,w}$ and the following estimates hold. 
\bea
\forall |q| \leq c_0, 1 \leq j \leq N &:&  
\left\| \frac{\partial}{\partial q_j} u(q) \right\| _{-\sigma _1,w}
\leq C_0. 
\eea
\item[(3)]
For all $|q| \leq c_0, n \geq 0$ and $m \in \ganz :
u(q)(n,-m) = \overline{u(q)(n,m)}$.
\end{itemize}
Then for all $ (q, \epsilon)  \in \relz ^{N+1}$ with 
$|q|, |\epsilon | \leq \min (c_0, \frac{c_1}{4 C_0}) $
there exists a unique $v \in 
{\cal L} _{\sigma _0,w}$ with the following properties (i)-(iii). 
\begin{itemize}
\item[(i)]
  $ \forall m \in \ganz, n \geq 1:
  v(n-1,m) +\delta _m v(n,m) +v(n+1,m) +W(u(q),v)(n,m) = 0. $
\item[(ii)]
  $\| v\| _{ \sigma _0,w} \leq 2(\rho _{u(q)} +|\epsilon |)  \leq c_1
  \leq \frac{\rhof }{8}.$
\item[(iii)]
  $\forall |m| > m_0 : v(0,m) = \epsilon b_m - u(q)(0,m).$
\end{itemize}
Furthermore the following holds.
\begin{itemize}
\item[(iv)]
For all $n \geq 0$ and $m \in \ganz :
v(n,-m) = \overline{v(n,m)}$.
\item[(v)]
  The map $(q,\epsilon ) \longmapsto v$ is a $C^1$ map into 
  ${\cal L} _{\sigma _1,w}.$
\end{itemize}
\end{theorem}

\begin{proof}  
The proof proceeds via a Banach fixed-point argument for $v$ and the
derivatives of $v$ with respect to the parameters $q_j, \epsilon$. 
First we define a map $\tilde{T}$ (Step~1),
which we then show to be a contraction on a certain set
(Step~2). We conclude that the first component of the fixed-point
of this map is a
solution of the equations given in (i) and (iii) (Step~3). 
Step 4 settles the question of differentiability for $v$ and 
Step 5 deals with the remaining properties (ii) and (iv).
\newline

{\bf Step 1:} Definition of the contraction.

First we will turn equation (2) of Lemma \ref{lb2.1} 
into a fixed-point equation
$v = T(q,\epsilon,v)$, depending on the parameters $q$ and $\epsilon$,
by applying the inverse of $\bigoplus _{m \in \ganz} L_m$ on it. As it was 
pointed out in the introduction of the present chapter, 
the operators $L_m$ are not invertible 
for $|m| \leq m_0$. Nevertheless we will define a formal inverse for
$|m| \leq m_0$, acting on exponentially decaying sequences,
as motivated in the introduction. 
Note the special role of
$m = 0$, where the Green's function of $L_0$ grows linearly. Using the fact
that the nonlinear term is given by
$W(u,v)(n,0) = Y(u,v)(n,0) - Y(u,v)(n+1,0)$, we end up with a bounded
kernel acting on $Y(u,v)(n,0)$.
We then proceed to define the maps $\tilde{T} _{q, \epsilon, j}$,
which give rise to the fixed-point equation for $v$ in their first argument
and to the fixed-point equation for the partial derivative of $v$ with 
respect to the $j$-th component of the parameters in their second
argument. The map $\tilde{T}$ is introduced in order to show that
the fixed-point $v$ of the map $T$ depends smoothly on the parameters
$q, \epsilon$ in Step~4.

Let $\nu :=(q,\epsilon)$ denote the
parameters in the construction. Furthermore it is convenient to
scale these parameters by a factor 
\be
\label{b3.8}
\eta := \min \left( c_0, \frac{c_1}{4 C_0} \right). 
\ee
Hence we can choose
$\nu \in U$, with
\be
\label{b3.9}
U := \{(q, \epsilon) \in \relz ^{N+1}:
|q|, |\epsilon| < 1 \}.
\ee
Let us further define two sets on which the map will act.
\be
\label{b3.10}
B := \{ v \in {\cal L} _{\sigma _0 ,w} :
\| v \| _{\sigma _0,w} \leq c_1 \}.
\ee
\be
\label{b3.11}
B' := \{ y \in {\cal L} _{\sigma _0 -2 \sigma _1,w} :
\| y \| _{\sigma _0 -2 \sigma _1,w} \leq c_1 \}.
\ee

We now define the map $T(\nu, v)$ explicitly.
For $\nu \in U, v \in B$ let
\begin{itemize}
\item
  $m=0 , n \geq 0 :$
  \bd
  T(\nu,v)(n,0) :=-\sum_{k=n+1}^{\infty} Y(u(\eta q),v)(k,0).
  \ed
\item
  $0< |m| \leq m_0, n \geq 0 :$
  \bd
  T(\nu,v)(n,m) :=  \frac{1}{\sin \beta _m}
    \sum_{k=n+1}^{\infty} W(u(\eta q),v)(k,m) \sin (n-k)\beta _m.
  \ed
\item
  $|m| > m_0 , n=0 :$
  \bd
  T(\nu,v)(0,m) :=  \eta \epsilon b_m-u(\eta q)(0,m).
  \ed
\item
  $ |m| > m_0, n \geq 1 :$
  \bea
  T(\nu,v)(n,m) &:=& 
  \sum_{k=1}^{\infty}  \frac{1-z_m^{2 \min (n,k)}}{ 1-z_m^2} z_m^{|n-k|+1}
     W(u(\eta q),v)(k,m) \\ 
     &+&  z_m^n ( \eta \epsilon b_m - u(\eta q)(0,m) ).
  \eea
\end{itemize}
It is useful to rewrite $T$ in the following way.
Define the linear map $K$, which acts on spaces
${\cal L} _{\sigma ,w}.$
\be
\label{b3.15}
(Ky)(n,m) := \sum_{k=1}^{\infty} K(k,n,m)y(k,m), 
\ee
with kernel $K(k,n,m)$.
\begin{itemize}
\item
  $m=0 :$
\bd
  K(k,n,0) :=  \left\{
  \begin{array}{ccc}
    0, & \mbox{ for } & k \leq n  \\
    -1, & \mbox{ for } & k \geq n+1 .
  \end{array}
  \right.
\ed
\item
  $0 < |m| \leq m_0 :$
\bd
  K(k,n,m) :=  \left\{
  \begin{array}{ccc}
    0, & \mbox{ for } & k \leq n  \\
    \frac{1}{\sin \beta _m} 
    [\sin (n-k)\beta _m - \sin (n+1-k)\beta _m],
    & \mbox{ for } & k > n. 
  \end{array}
  \right.  
\ed
\item
  $|m| > m_0 :$
\bd
  K(k,n,m) :=
  \frac{1-z_m^{2 \min (n,k)}}{ 1-z_m^2} z_m^{|n-k|+1} -
  \frac{1-z_m^{2 \min (n,k-1)}}{ 1-z_m^2} z_m^{|n+1-k|+1}.
\ed
\end{itemize}
It is straightforward to check that 
\be
\label{b3.16}
\hspace{.3in}
T(\nu,v)(n,m) = 
[K Y(u(\eta q),v)](n,m) + 
(\eta \epsilon b_m - u(\eta q)(0,m) ) z_m^n 
{\bf 1}_{ \{ |m| > m_0 \} }.
\ee

Now we can define the map $\tilde{T}$, which
acts on the complete metric space $B \times B'$, equipped
with the norm
\be
\label{b3.18}
\| (v,y) \| _{B \times B'} := 
\max \left(
\| v \| _{\sigma _0,w}, \| y \| _{\sigma _0 - 2 \sigma _1}  \right).
\ee
Fix $\nu \in U, 1 \leq j \leq N+1$.
\newline
$ \tilde{T} _{\nu,j} : B \times B' \longrightarrow B \times B' $,
\be
\label{b3.20}
\tilde{T} _{\nu,j}
\left(
\begin{array}{c} v \\ y  
\end{array}
\right)
:=
\left(
\begin{array}{c}  
  T(\nu,v) \\
  h_j(\nu , v) + K (D_v Y)(u(\eta q),v) y
\end{array}
\right),
\ee
where we now give explicit expressions for $h_j(\nu , v)$. There are
two cases.
\newline
Case 1:  $\nu _j = \epsilon$   \newline
\be
\label{b3.21}
h_j(\nu , v)(n,m) :=
\eta b_m  z_m^n {\bf 1}_{ \{ |m| > m_0 \} }.
\ee
Case 2:  $\nu _j = q_j$.  \newline
\be
\label{b3.22}
h_j(\nu , v)(n,m) :=
\eta [K (D_u Y) (D_j u) ](n,m)
- \eta (D_j u) (0,m) z_m^n {\bf 1}_{ \{ |m| > m_0 \} }.
\ee
This completes the definition of $\tilde{T}$.
\newline

{\bf Step 2:} 
$\tilde{T}_{\nu,j} : B \times B' \longrightarrow B \times B' $ 
is a contraction.
\newline
\newline
We first obtain a bound on the norm of the linear operator $K$.
\lfd
\begin{proposition}
\label{pb3.1}
For all $\sigma_1 \leq \sigma \leq \sigma _0$, the linear operator $K$ maps 
${\cal L} _{\sigma,w}$ into ${\cal L} _{\sigma,w}$ and the corresponding 
operator norms of $K$ are bounded by $C_K$. (See Definition \ref{db3.1}). 
\end{proposition}

\begin{proof}
The proof is a consequence of the following estimates.
\begin{itemize}
\item
  $m=0, n \geq 0:$
  \begin{eqnarray*}
  |(Ky)(n,0)|e^{ \sigma n}  & \leq & 
  e^{ \sigma n} \sum_{k=n+1}^{\infty} e^{- \sigma k}
  \sup_{j \geq 0} |e^{ \sigma j} y(j,0)| \\
  & \leq &
  \frac{1}{1-e^{- \sigma _1}}
  \sup_{j \geq 0} |e^{ \sigma j} y(j,0)|. \\
  \end{eqnarray*}
\item
  $0 < |m| \leq m_0, n \geq 0:$
  \begin{eqnarray*}
    |(Ky)(n,m)|e^{ \sigma n}  & \leq & 
    \frac{2}{\sin \beta _m}
    \sum_{k=n+1}^{\infty} e^{- \sigma k}
    \sup_{j \geq 0} |e^{ \sigma j} y(j,m)| \\
    & \leq &
    \frac{1}{1-e^{- \sigma _1}}
    \frac{2}{\sin \beta _m}
    \sup_{j \geq 0} |e^{ \sigma j} y(j,m)|. \\
  \end{eqnarray*}
\item
  $|m| > m_0, n = 0: K(k,0,m) = 0$.
\item
  $|m| > m_0, n \geq 1:$
  By definition of $\sigma _0$, $|z_m| \leq e^{ -2 \sigma _0 }$
  (see Definition \ref{db3.1}),
  which yields the estimate for the Greens function 
  $\left| \frac{1-z_m^{2 \min (n,k)}}{ 1-z_m^2} z_m^{|n-k|+1}    
  \right|  \leq
  \frac{ e^{-2 \sigma _0}}{ 1 - e^{-4 \sigma _0}}
  e^{-2 \sigma _0 |n-k|}$. Convolution of this bound with an
  exponentially decaying sequence gives
  \begin{eqnarray*}
    \sum_{k=1}^{\infty} e^{-2 \sigma _0 |n-k|} e^{- \sigma k}
    & \leq &
    \sum_{k=0}^{n-1} e^{-2 \sigma _0 (n-k) - \sigma k} +
    \sum_{k=n}^{\infty} e^{-2 \sigma _0 (k-n) - \sigma k}   \\
    & \leq &
    e^{-2 \sigma _0 n} \frac{e^{ (2\sigma _0 - \sigma ) n} -1}
                            {e^{ (2\sigma _0 - \sigma )} -1} +
    e^{- \sigma n} \frac{1}{1 - e^{-2 \sigma _0 - \sigma }} \\
    & \leq &
    e^{- \sigma n} \left( 
    \frac{1}{e^{ \sigma _0} -1} + \frac{1}{1 - e^{-2 \sigma _0}}  \right), \\
  \end{eqnarray*}
  and we arrive at
  \begin{eqnarray*}
    |(Ky)(n,m)|e^{ \sigma n}  & \leq &
    \frac{ 1 + e^{-2 \sigma _0}}{ 1 - e^{-4 \sigma _0}}  \left(
    \frac{1}{e^{ \sigma _0} -1} + \frac{1}{1 - e^{-2 \sigma _0}}  \right)
    \sup_{j \geq 0} |e^{ \sigma j} y(j,m)|. \\
  \end{eqnarray*}
\end{itemize}
\end{proof}

Now we are ready to prove that $\tilde{T} _{\nu,j}$ is a contraction. 
We use the estimates of 
Proposition \ref{pb2.3} and equations (\ref{b3.16})-(\ref{b3.22}).
\begin{itemize}
\item
\bd
\| T(\nu,v) \| _{\sigma _0, w}   \leq 
C_K C_{F,c} c_1^2 + \frac{c_1}{4} + \frac{c_1}{4} \leq 
\frac{3}{4} c_1.
\ed
Furthermore it was shown in Proposition \ref{pb2.3} (iii)  
that, for the $u,v$ in question, the map 
$v \mapsto Y(u,v)$ is $C^1$ from 
${\cal L} _{\sigma _0, w}$ into ${\cal L} _{\sigma _0, w}$ and the 
ususal operator norm of the derivative is bounded by $C_{F,c} c_1$.
Therefore
\bd
\| T(\nu,v') - T(\nu,v) \| _{\sigma _0, w}   \leq
C_K C_{F,c} c_1 \| v' - v \| _{\sigma _0, w} \leq
\frac{1}{4} \| v' - v \| _{\sigma _0, w}.
\ed
\item
We have for $\nu \in U,v \in B, y \in B'$, that
\bd
\| K (D_v Y)(u(\eta q),v) y \| _{\sigma _0 - 2 \sigma _1, w} \leq
C_K C_{F,c} c_1^2 \leq \frac{1}{4} c_1.
\ed
The bound on the second derivative of $Y$ with respect to $v$ 
(see Proposition \ref{pb2.3} (iii))
allows us to make the following estimate.
\bea
&&
\| K (D_v Y)(u,v') y' - K (D_v Y)(u,v) y \| 
_{\sigma _0 - 2 \sigma _1, w}  \\
&\leq&
\| K [(D_v Y)(u,v')  - (D_v Y)(u,v)] y' \| 
_{\sigma _0 - 2 \sigma _1, w}  \\
&& +
\| K (D_v Y)(u,v) ( y' - y) \| 
_{\sigma _0 - 2 \sigma _1, w}   \\
&\leq&
C_K C_{F,c} c_1 \| v' - v \|_{\sigma _0, w}  +
C_K C_{F,c} c_1 \| y' - y \|_{\sigma _0 - 2 \sigma _1, w}   \\
&\leq&
\frac{1}{2} \max \left(
\| v' - v \|_{\sigma _0, w},  
\| y' - y \|_{\sigma _0 - 2 \sigma _1, w}
\right).
\eea
\item
Finally we are dealing with the estimates for $h_j(\nu,v)$.
Corresponding to the definition (\ref{b3.21}),
(\ref{b3.22}) we have to distinguish two cases.
\newline
Case 1:  $\nu _j = \epsilon$   \newline
\bd
\| \eta b_m  z_m^n {\bf 1}_{ \{ |m| > m_0 \} } \| 
_{\sigma _0 - 2 \sigma _1, w}
\leq \eta \leq \frac{1}{4} c_1.
\ed
Case 2:  $\nu _j = q_j$.  \newline
Using hypothesis (2) of Theorem \ref{tb3.1} we see immediately that
\bd
\| \eta (D_j u) (0,m) z_m^n {\bf 1}_{\{ |m| > m_0 \}} \| 
_{\sigma _0 - 2 \sigma _1, w}
\leq  \eta C_0 \leq \frac{1}{4} c_1.
\ed
From Proposition \ref{pb2.3} (iv), we obtain
\bd
\| \eta K (D_u Y)(u,v) (D_j u) \| _{\sigma _0 - 2 \sigma _1, w}   \leq
2 \eta C_K C_{F,c} c_1 C_0 \leq \frac{1}{8} c_1.
\ed
Finally we use that for fixed $x \in {\cal L} _{-\sigma _1, w}$, the map
$v \mapsto (D_u Y)(u,v) x$ is $C^1$ from ${\cal L} _{\sigma _0, w}$ into
${\cal L} _{\sigma _0 - 2 \sigma _1, w}$ 
with a bound on the derivative as given in
Proposition \ref{pb2.3} (v).
\bea
\| \eta K [(D_u Y)(u,v') - (D_u Y)(u,v)] (D_j u) \| 
_{\sigma _0 - 2 \sigma _1, w}
&\leq& \eta C_K C_{F,c} C_0 \| v' - v \|_{\sigma _0, w}  \\
&\leq& \frac{1}{4} \| v' - v \|_{\sigma _0, w}.
\eea
\end{itemize}

The claim of the present step is a consequence of all these estimates. 
Thus we have proven the existence of a fixed-point of 
$\tilde{T} _{\nu,j}$ in $B \times B'$.\newline \newline

{\bf Step 3 :} For $v \in B$ the following eqivalence holds. \newline
v satisfies properties (i) and (iii) of Theorem \ref{tb3.1}
$ \Longleftrightarrow  T(\nu,v) = v$.
\newline

{\bf $\Longleftarrow$ :}  \newline
Recall that we have scaled $\nu$ by the factor $\eta$ in the beginning 
of the proof. Then
property (iii) is evident as $T(\nu,v)$ satisfies it by definition. 
Therefore it suffices to prove that for all $v \in B,
n \geq 1, m \in \ganz$:
\bd
T(\nu,v)(n-1,m) + \delta _m T(\nu,v)(n,m) + T(\nu,v)(n+1,m) = -W(u,v)(n,m).
\ed
We will verify this by evaluating the lefthandside for all different cases.
\begin{itemize}
\item
  $m=0, n\geq 1:$
  \bea
     \mbox{ LHS } &=& 
     - \sum_{k=n}^{\infty} Y(u,v)(k,0)
     +2 \sum_{k=n+1}^{\infty} Y(u,v)(k,0)
     - \sum_{k=n+2}^{\infty} Y(u,v)(k,0)    \\
     &=&
     -Y(u,v)(n,0) + Y(u,v)(n+1,0).  \\
  \eea
\item
  $0 < |m| \leq m_0, n \geq 1:$
  \bea
\hspace*{-.12in}
     \mbox{ LHS } &=& 
     \frac{1}{\sin \beta _m}  \left(
     \sum_{k=n+1}^{\infty} W(u,v)(k,m) G(k,n,m)
     + W(u,v)(n,m) \sin (-\beta _m)  \right),   \\
   \eea
   where by the definition of $\beta _m$ (see (\ref{i1.25}))
   \bd
   G(k,n,m) := 
   \sin (n-1-k) \beta _m  -2 \cos \beta _m  \sin (n-k) \beta _m  
   + \sin (n+1-k) \beta _m.  
   \ed
  As $\sin (n-1-k) \beta _m + \sin (n+1-k) \beta _m  =
  2 \cos \beta _m  \sin (n-k) \beta _m $, we conclude $G(k,n,m) = 0.$
\item
  $|m| > m_0 , n=1:$
  \bea
     \mbox{ LHS } &=& 
     ( 1+ \delta _m z_m + z_m^2 ) 
     ( \eta \epsilon b_m - u(\eta q)(0,m) )\\
     && + \sum_{k=1}^{\infty} W(u,v)(k,m) \left(
     \delta _m z_m^k + \frac{1 - z_m^{2 \min (2,k)}}{1-z_m^2} z_m^{|2-k|+1}
     \right). \\
  \eea
  The identity $1+ \delta _m z_m + z_m^2 =0$ implies that all terms 
  in the sum vanish with the exception of $k=1$. In fact
  \begin{itemize}
  \item[$\bullet$]
    $k=1 :$
    \bea
      \delta _m z_m^k + \frac{1 - z_m^{2 \min (2,k)}}{1-z_m^2} z_m^{|2-k|+1}
      & = &
      \delta _m z_m + z_m^2 = -1. \\
    \eea
  \item[$\bullet$]
    $ k \geq 2:$
    \bea
      \delta _m z_m^k + \frac{1 - z_m^{2 \min (2,k)}}{1-z_m^2} z_m^{|2-k|+1}
      & = &
      \delta _m z_m^k + z_m^{k-1} (1 + z_m^2) = 0. \\
    \eea
  \end{itemize}
\item
  $|m| > m_0, n \geq 2:$
  \bea
    \mbox{ LHS } &=&
    \left( z_m^{n-1} + \delta _m z_m^{n} + z_m^{n+1} \right) 
    ( \eta \epsilon b_m - u(\eta q)(0,m) )  \\
    &&
    + \sum_{k=1}^{\infty} W(u,v)(k,m) G(k,n,m) , \\
  \eea
  with
  \bea
    G(k,n,m) & := &
    \frac{1 - z_m^{2 \min (n-1,k)}}{1-z_m^2} z_m^{|n-1-k|+1} +
    \delta _m \frac{1 - z_m^{2 \min (n,k)}}{1-z_m^2} z_m^{|n-k|+1}  \\
    &&
    + \frac{1 - z_m^{2 \min (n+1,k)}}{1-z_m^2} z_m^{|n+1-k|+1}. \\
  \eea
   We show that $G(k,n,m) = - \delta _{k,n}$. To that end we
  evaluate $(1 - z_m^2) G(k,n,m) =$
  \begin{itemize}
  \item[$\bullet$]
    for $k < n$: \newline
    $ = (1 - z_m^{2k}) z_m^{n-k} (1 + \delta _m z_m +z_m^2) = 0.$
  \item[$\bullet$]
    for $k = n$: \newline
    $=  (z_m^2 + \delta _m z_m + z_m^2)
    - z_m^{2(n-1)} ( z_m^2 + \delta _m z_m^3 + z_m^4)
    = z_m^2 - 1.$
  \item[$\bullet$]
    for $k > n$: \newline
    $= z_m^{k-n} (1 + \delta _m z_m +z_m^2)
    - z_m^{k-n} z_m^{2(n-1)} ( z_m^2 + \delta _m z_m^3 + z_m^4)  = 0.$
  \end{itemize}
\end{itemize}
{\bf $\Longrightarrow$ :}  \newline
Let $v \in B$ satisfy (i) and (iii). Denote $d := v - T(\nu, v).$ We want 
to show that $d = 0$. 
The following list of properties for $d$ is immediate. 

\begin{itemize}
\item[(I)]
$d \in {\cal L} _{\sigma _0, w}.$ 
\item[(II)]
$d(0,m) = 0$ for $|m| > m_0$ (as $v$ satisfies (iii)). 
\item[(III)]
$\forall m \in \ganz, n \geq 1:
d(n-1,m) + \delta _m d(n,m) + d(n+1,m) = 0.$
\end{itemize}
The last relation enables us to express $d(n,m)$ in terms of $d(0,m)$
and $d(1,m).$
\begin{itemize}
\item 
  $m=0 :
  d(n,0) = n d(1,0) - (n-1) d(0,0)$ . As $d$ has to decay exponentially with
  $n \rightarrow \infty$ there is no other choice than $\forall n \geq 0:
  d(n,0)=0.$
\item
  $0<|m|\leq m_0:$
  \bea
    d(n,m) &=& \frac{ \sin  n\beta _m}{ \sin \beta _m} d(1,m)
    - \frac{ \sin  (n-1)\beta _m}{ \sin \beta _m} d(0,m)  \\
    &=&
    \frac{\sin  n\beta _m }{ \sin \beta _m} 
    \left( d(1,m) - d(0,m) \cos \beta _m  \right)
    + \cos n\beta _m  d(0,m). \\
  \eea
   $ \beta _m \in (0,\pi)$ and again the exponential decay of $d(\cdot,m)$
   force $d(0,m)$ and $d(1,m)$ to be zero and hence $d(n,m) = 0$ 
   for all $n \geq 0.$
\item
  $|m| > m_0 :$ Recall that we know already 
  by (II) that $d(0,m) = 0.$ Then
  $d(n,m) =  d(1,m) \frac{z_m^n - z_m^{-n}}{z_m - z_m^{-1}}$. Again we 
  conclude that $d(1,m) = 0.$
\end{itemize}

This concludes the proof of Step 3.
\newline

{\bf Step 4:} Differentiability of $v$ with respect to the parameters $\nu$.
\newline \newline
It is a standard and well known problem to prove smooth dependence
of the solution of a contraction problem on the parameters.
There is a small technical problem as $u$ is a $C^1$-function of the 
parameter $q$ only in the space ${\cal L} _{- \sigma _1,w}$, and not in the
space ${\cal L} _{0,w}$. But, at least formally, 
differentiating $v = T(\nu, v)$ with respect to $\nu _j$ gives
\bd
\frac{\partial v}{\partial \nu _j} =
(1 - D_v T)^{-1} \frac{\partial}{\partial \nu _j} T (\nu ,v).
\ed
The lack of differentiability of $u(q)$ in ${\cal L} _{0,w}$ prevents
that $\frac{\partial}{\partial \nu _j} T$ lies in
${\cal L} _{\sigma _0,w}$ , but it lies in 
${\cal L} _{\sigma _0 - 2 \sigma _1,w}$. Fortunately, however,
$D_v T$ maps any ${\cal L} _{\sigma ,w}$, $\sigma _1 \leq 
\sigma \leq \sigma _0$, into itself with small norm. Hence $v$ is
differentiable in an appropriate norm.
We find it convenient to proceed as follows.

The fixed-point $(v,y)(\nu)$ of $\tilde{T} _{\nu,j}$ can be constructed as 
the limit of the iteratives of the map, i.e. let 
\bd
(v_0(\nu),y_0(\nu)) := (0,0)
\ed
and define inductively 
\bd
(v_{s+1}(\nu),y_{s+1}(\nu)) := \tilde{T} _{\nu,j} (v_s(\nu),y_s(\nu)),
\ed 
then $(v(\nu), y(\nu)) = \lim_{s \rightarrow \infty} (v_s(\nu),y_s(\nu))$.
The limit is in the norm of $B \times B'$ and uniform in $\nu$.

We will prove inductively for each choice of $\nu _j$ and for all 
$s \geq 0$ that the following properties hold.
\begin{itemize}
\item[(I)]
$U \rightarrow {\cal L} _{\sigma _0 - \sigma _1,w}: \nu
\mapsto v_s (\nu)$ is continuous.
\item[(II)]
For each variable $\nu _j$, the map 
$\nu _j \mapsto v_s (\nu)$ is differentiable as a map into 
${\cal L} _{\sigma _0 - 2 \sigma _1,w}$ and
\bd
\frac{\partial}{\partial \nu _j} v_s (\nu) = y_s (\nu).
\ed
\item[(III)]
$U \rightarrow {\cal L} _{\sigma _0 - 3 \sigma _1,w}: \nu
\mapsto y_s (\nu)$ is continuous.
\end{itemize}
Once we have established (I)-(III), the proof of the claim is immediate.
In fact, we can deduce that for all $s \in \natz _0:
\nu \mapsto v_s (\nu)$ is a $C^1$ map from $U$ into
${\cal L} _{\sigma _0 - 3 \sigma _1,w} =
{\cal L} _{\sigma _1,w}$ (see Definition \ref{db3.1}), with partial derivatives
$\frac{\partial v_s}{\partial \nu _j} = y_s$. The uniform convergence
of $(v_s, y_s)(\nu)$ to the fixed-point of $\tilde{T} _{\nu,j}$,
$(v,y)(\nu)$, in the 
norm of $B \times B'$ yields the desired information, that 
$\nu \mapsto v(\nu)$ is a $C^1$ map from $U$ into
${\cal L} _{\sigma _1,w}$ and that $\frac{\partial v}{\partial \nu _j} = y$.

Let us therefore return to the statements (I)-(III). 
They are trivially satisfied for $s = 0$. We will now prove the induction
step  $ s \rightarrow s + 1$.
\newline \newline
{\bf (I) Continuity:}
\newline \newline
Let $ \nu, \nu ' \in U$. We have to show that 
\be
\label{b3.30}
\| T(\nu ', v_s(\nu ')) - T(\nu , v_s(\nu )) \|
_{\sigma _0 - \sigma _1,w} \longrightarrow 0, \mbox{ as } 
\nu ' \rightarrow \nu.
\ee
To simplify the
notation, denote $v := v_s (\nu), v' := v_s (\nu '),
u := u(\eta q), u' := u(\eta q')$.
It is not hard to verify that
\be
\label{b3.31}
\| [(\eta \epsilon ' b_m - u'(0,m) ) - (\eta \epsilon  b_m - u(0,m) )]
z_m^n {\bf 1}_{ \{ |m| > m_0 \} } \| _{\sigma _0 - \sigma _1,w}
\leq \eta | \epsilon ' - \epsilon | + 
\| u' - u \| _{- \sigma _1,w}.
\ee
Expressing the difference by telescoping sums, one obtains for $l \geq 1$
\bea
&&\| (\triangle u')^{*(k-l)} * (\triangle v')^{*l}  -
   (\triangle u)^{*(k-l)} * (\triangle v)^{*l}
\| _{\sigma _0 - \sigma _1,w}  \\
&\leq&
(1 + e^{\sigma _0}) k (2 c_1)^{k-1} 
\left( \| u' - u \| _{- \sigma _1,w} +
\| v' - v \| _{\sigma _0 - \sigma _1,w}  \right).
\eea

Proposition \ref{pb2.1} and Proposition \ref{pb3.1} imply
\lfd
\begin{eqnarray}
&& \| K Y(u',v') - K Y(u,v) \| _{\sigma _0 - \sigma _1,w} \nonumber \\
&\leq&
8 (1 + e^{\sigma _0}) C_K \tilde{C} _{F,c} c_1
\left( \| u' - u \| _{- \sigma _1,w} +
\| v' - v \| _{\sigma _0 - \sigma _1,w}  \right). \label{b3.33}
\end{eqnarray}
Equations (\ref{b3.31}), (\ref{b3.33}), assumption (2) of the Theorem
and the induction hypothesis suffice to prove equation (\ref{b3.30}).
\newline \newline
{\bf (II) Existence of partial derivatives:}
\newline \newline
We consider only the more difficult case $\nu _j = q_j$.
The proof for the case $\nu _j = \epsilon$ requires only
a proper subset of the arguments given below.
Denote
$\nu ' := (q + h q_j, \epsilon)$ and let $u',v',u,v, D_j u',D_j u,  
D_j v',D_j v$ have the obvious meaning. We have to prove
\be
\label{b3.40}
\lim_{h \rightarrow 0}
\left\| \frac{1}{h} \left(
T(\nu ', v') - T(\nu ,v) \right) - y_{s+1} (\nu) \right\| 
_{\sigma _0 - 2 \sigma _1,w}  = 0.
\ee
We break this statement up into several estimates.
\lfd
\begin{eqnarray}
\label{b3.41}
&&
\left\| \left[ \frac{1}{h} \left( u'(0,m) - u(0,m) \right) -
\eta (D_j u)(0,m) \right]
z_m^n {\bf 1}_{ \{ |m| > m_0 \} } \right\| _{\sigma _0 - 2 \sigma _1,w} 
\nonumber \\
&\leq& 
\left\| \frac{1}{h}(u' - u) - \eta D_j u \right\| _{- \sigma _1,w},
\end{eqnarray}
which tends to $0$ as $h \rightarrow 0$ by assumption (2).
Next we look at the monomials, from which $Y, D_u Y$ and
$D_v Y$ are built up. We use the induction hypothesis, which says
that $y_s = D_j v_s$.
\bea
&&
\frac{1}{h} \left[
(\triangle u')^{*(k-l)} * (\triangle v')^{*l}  -
(\triangle u)^{*(k-l)} * (\triangle v)^{*l}
\right]  
\\
&&
- (k-l)(\triangle u)^{*(k-l-1)} * (\triangle v)^{*l} 
* \eta \triangle(D_j u)
\\
&&
- l(\triangle u)^{*(k-l)} * (\triangle v)^{*(l-1)} 
* \triangle(D_j v)
\\
&=&
\left[
\frac{1}{h} \left( 
(\triangle u')^{*(k-l)}  - (\triangle u)^{*(k-l)}
\right) -
(k-l)(\triangle u)^{*(k-l-1)} * \eta \triangle(D_j u)
\right] 
* (\triangle v')^{*l} 
\\
&&
+ (k-l)(\triangle u)^{*(k-l-1)} * \eta \triangle(D_j u)
* \left( (\triangle v')^{*l} - (\triangle v)^{*l}  \right) 
\\
&&
+ (\triangle u)^{*(k-l)} *
\left[
\frac{1}{h} \left( 
(\triangle v')^{*l}  - (\triangle v)^{*l}
\right) -
l (\triangle v)^{*(l-1)} * \triangle(D_j v)
\right]
\\
&=&
(a) + (b) + (c).
\eea
The three terms are now investigated seperately. Assume $l \geq 1$.
\begin{itemize}
\item[(a)]
Telescoping differences twice we obtain
\bea
&&
\frac{1}{h} \left( 
(\triangle u')^{*(k-l)}  - (\triangle u)^{*(k-l)}
\right) -
(k-l)(\triangle u)^{*(k-l-1)} * \eta \triangle(D_j u)
\\
&=&
\sum_{j=1}^{k-l-1} \sum_{i=0}^{j-1}
(\triangle u')^{*i} * (\triangle u)^{*(k-l-i-2)}
* \triangle (u' - u) * \triangle \left(\frac{u' - u}{h} \right)
\\
&&
+ (k-l)(\triangle u)^{*(k-l-1)} 
* \triangle \left(\frac{u' - u}{h} - \eta D_j u \right).
\eea
This implies, that
\bea
\| (a) \| _{\sigma _0 - 2 \sigma _1,w}
&\leq&
2 (1 + e^{\sigma _0}) k (k - 1) (2 c_1)^{k-2} \\
&&
\hspace*{-.15in}
\times
\left(
\| u' - u \| _{- \sigma _1,w}
\left\| \frac{u' - u}{h} \right\| _{- \sigma _1,w} +
c_1 \left\|
\frac{u' - u}{h} - \eta D_j u
\right\| _{- \sigma _1,w}
\right).
\eea
\item[(b)]
Recalling the definition of $C_0$ in the statement of the theorem, it 
is easy to see that
\bd
\| (b) \| _{\sigma _0 - 2 \sigma _1,w}
\leq
2 (1 + e^{\sigma _0}) k (k - 1) (2 c_1)^{k-2} 
\eta C_0 \| v' - v \| _{\sigma _0 - \sigma _1,w}.
\ed
\item[(c)]
Proceeding as in (a) we obtain
\bea
\| (c) \| _{\sigma _0 - 2 \sigma _1,w}
&\leq&
2 (1 + e^{\sigma _0}) k (k - 1) (2 c_1)^{k-2} \\
&&
\hspace*{-.3in}
\times
\left(
\| v' - v \| _{0,w}
\left\| \frac{v' - v}{h} \right\| _{\sigma _0 - 2 \sigma _1,w} +
c_1 \left\|
\frac{v' - v}{h} - D_j v
\right\| _{\sigma _0 - 2 \sigma _1,w}
\right).
\eea
\end{itemize}
Substituting all these estimates in the power series and using
Proposition \ref{pb2.1}, 
the induction hypothesis, assumption (2)
and equation (\ref{b3.41}) we arrive at the assertion of equation 
(\ref{b3.40}).
\newline \newline
{\bf (III) Continuity of partial derivatives:}
\newline \newline
Using the notation and the methods of the last two proofs, we see
that the following estimates are enough to prove the claim.
\begin{itemize}
\item
\bea
&&
\| \eta ( (D_j u')(0,m) - (D_j u)(0,m)) 
z_m^n {\bf 1}_{ \{ |m| > m_0 \} } \| _{\sigma _0 - 3 \sigma _1,w} \\
&\leq&
\eta
\| D_j u' - D_j u \| _{- \sigma _1,w}.
\eea
\item
The usual telescoping technique readily yields for $l \geq 1$ 
the following estimate.
\bea
&&
\hspace*{-.48in}
\left\|
(\triangle u')^{*(k-l-1)} * (\triangle v')^{*l}  
* \triangle (D_j u') -
(\triangle u)^{*(k-l-1)} * (\triangle v)^{*l}
* \triangle (D_j u)
\right\| _{\sigma _0 - 3 \sigma _1,w}
\\
&\leq&
2 (1 + e^{\sigma _0}) k (2 c_1)^{k-2} \\
&&
\times
\left(
C_0 \| u' - u \| _{- \sigma _1,w} +
C_0 \| v' - v \| _{\sigma _0 - \sigma _1,w} +
c_1 \| D_j u' - D_j u \| _{- \sigma _1,w}
\right).
\eea
\item
We know from Step 2 that 
$\| D_j v' \| _{\sigma _0 - 2\sigma _1,w},
\| D_j v \| _{\sigma _0 - 2\sigma _1,w} \leq c_1$. Proceeding as above we 
obtain
\bea
&&
\hspace*{-.57in}
\left\|
(\triangle u')^{*(k-l)} * (\triangle v')^{*(l-1)}  
* \triangle (D_j v') -
(\triangle u)^{*(k-l)} * (\triangle v)^{*(l-1)}
* \triangle (D_j v)
\right\| _{\sigma _0 - 3 \sigma _1,w}
\\
&\leq&
(1 + e^{\sigma _0}) k (2 c_1)^{k-1} \\
&&
\times
\left(
\| u' - u \| _{- \sigma _1,w} +
\| v' - v \| _{\sigma _0 - \sigma _1,w} +
\| D_j v' - D_j v \| _{\sigma _0 - 3 \sigma _1,w}
\right).
\eea
\end{itemize}
The proof of Step 4 is completed.
\newline

{\bf Step 5:} The remaining properties (ii) and (iv).
\newline

(ii) is a consequence of the fact that for all $\nu \in U$ the 
map $v \mapsto T(\nu, v)$ sends 
$\{ v \in {\cal L}_{\sigma _0, w}: \| v \| _{\sigma _0, w}
\leq 2(\eta |\epsilon| + \rho _{u(\eta q)}) \}$ into itself.
In fact, from Proposition \ref{pb2.3} and equation (\ref{b3.16}) 
it follows, that
\bea
\| T(\nu,v) \| _{\sigma _0, w}     & \leq &
C_K C_{F,c} c_1 \| v \| _{\sigma _0, w}  +
(\eta |\epsilon| + \rho _{u(\eta q)})  \\
& \leq &
(\eta |\epsilon| + \rho _{u(\eta q)})
(2 C_K C_{F,c} c_1 + 1).
\eea
Note that the appearance of $\eta$ in the proof, which is not present  
in the formulation of property (ii) in the theorem, comes from the 
scaling of the parameters which we performed at the beginning of the proof.

(iv) can be shown inductively, following the iterative construction of
the fixed-point which we have already employed in Step 4.
$v_0 = 0$ clearly satisfies the reality condition, 
i.e. $v_0(n,-m) = \overline{v_0(n,m)}$, for all $n \geq 0, m \in \ganz$. 
Using Propositon \ref{pb2.2} (v) it follows that 
$v_{s+1}$ satisfies the reality condition as 
$(b_m)_{m \in \ganz}, u(\eta q)$ and $v_s$ do. 
This property is  
preserved as we pass to the limit $s \longrightarrow \infty$.

\end{proof}

\section{An immediate application: high frequency driving}
\label{m4}

Suppose that $m_0 = 0$, or equivalently $\gamma ^2 > 4 F'(-c)$.
In this case we can choose $u(q)$ identically equal to zero in 
Theorem \ref{tb3.1} and we obtain $v(\epsilon)$ with the 
corresponding properties. If we now take these choices for $u$
and $v$ we see immediately that all the conditions in Lemma 
\ref{lb2.1} are satisfied and hence we have constructed a
periodic solution of the differential equation given by (\ref{i2.25}),
(\ref{i2.26}). 
This proves the following result.
\lfd
\begin{theorem}
\label{tb4.1}
Let $F, c, (b_m) _{m \in \ganz}, w$ satisfy the general assumptions. 
Furthermore let $\gamma ^2 > 4 F'(-c)$. Then there exists a neighborhood
$D$ of $0$, such that for all $\epsilon \in D$ there exists a
sequence $v(\epsilon) \in {\cal L} _{\sigma _0,w}$ with
$\sigma _0 > 0$ as defined in Definition \ref{db3.1} such that
\bd 
x_n(t) := cn + \epsilon b_0 - v(\epsilon)(0,0) + \sum_{m \in \ganz}
v(\epsilon) (n,m) e^{i m \gamma t} ,\quad \mbox { for } n \geq 1,
\ed 
is a time periodic, real valued solution of the the differential equation,
given  by (\ref{i2.25}) and
(\ref{i2.26}).
\end{theorem} 

\lfd
\begin{remark}
\label{rm4.1}
We have restricted our attention to the case where $F'(-c) > 0$.
If $F'(-c)$ is negative, however, then $\gamma ^2 > 4 F'(-c)$ for
all $\gamma \in \relz ^+$ and it is easy to see that 
Theorem \ref{tb4.1} holds without
any restrictions on the driving frequency.
\end{remark}

\chapter{General Lattices}
\label{l}

In the last section of the last chapter we have seen that in the 
case of $m_0 = 0$, it is possible to construct periodic solutions 
for arbitrary lattices, where the interacting forces between 
neighboring particles satisfy the general assumptions. The goal of this
chapter is to obtain the same result in the case of $m_0 = 1$, i.e.
\be
\label{l1.6}
F'(-c) < \gamma ^2 < 4 F'(-c).
\ee
The difference from the case $m_0 = 0$ is that now Theorem \ref{tb3.1}
no longer yields a sequence $v$ which solves 
$\epsilon b_m = u(0,m) + v(0,m)$, for all $m \neq 0$, but only for
$|m| > 1$. Therefore we have a resonance equation
for $m = 1$. In order to be able to solve this complex valued
equation ( for $m = -1$ the equation is the complex conjugate of the 
equation for $m = 1$), we must obtain a family of sequences $u(q)$ 
depending on two real parameters $q_1$ and $q_2$. The idea is to 
obtain $u(q)$ by constructing travelling wave solutions for the
doubly infinite lattice. Physically one may view these as the waves the driver
excites and which travel through the lattice. They can be observed
almost unperturbed away from the boundary some time after  
we let the driver act on the system (see Figures \ref{fc2} and \ref{fc5}).

More precisely we make the ansatz
\be
\label{l1.10}
x_n^{(0)} (t) := c n + \sum_{m \in \ganz}
r(m) e^{i m (\beta n + \gamma t)}.
\ee
In Section \ref{l2} we will give an equation for $r$ which we then solve
by a Lyapunov-Schmidt decomposition. The idea is as follows.
It is easy to see that 
the linearized  equation at $r = 0$ is given by a diagonal
operator $\Lambda (\beta)$ with entries 
\be
\label{l1.11}
\Lambda (\beta) (m,m) := \Lambda (\beta, m) := 2 \cos \beta m + \delta _m.
\ee
If we choose $\beta = \beta_1$ (see (\ref{i1.25}), $\beta _1$ exists as
$|\delta _1| < 2$) , then $\Lambda (\beta)$ has a nontrivial
kernel and we can apply the decomposition procedure.
The assumptions on $\gamma$ as given in (\ref{l1.6})
imply that $\delta _m > 2$ for $|m| > 1$ and therefore $\Lambda (\beta) (m,m)$
is bounded away from zero for 
arbitrary $\beta \in \relz$. Hence the 
infinite dimensional part of the decomposition
poses no problems. The finite dimensional
part needs additional consideration. First, one has to use 
various symmetries  
to obtain the correct count of variables. Then by expanding the degenerate
equation to second order one proves that the finite dimensional part
can be solved by choosing the spatial frequency $\beta$ as a function 
of the two remaining real parameters $q_1$ and $q_2$ with
\bd
|\beta (q) - \beta _1 | \leq C |q|^2.
\ed
\newline
In Section \ref{l3} we then take these solutions $r(q)$ and show that 
\be
\label{l1.15}
u(q)(n,m) := r(q)(m) e^{i m n \beta (q)} 
\ee
satisfies all the conditions of Theorem \ref{tb3.1} in Chapter \ref{m}.
Note, that from equation (\ref{l1.15}) it is already obvious that the
map $ q \mapsto u(q)$ cannot even be continuous in 
${\cal L} _{0,w}$, if $\beta(q)$ is not identically constant. This is
why we had to introduce the weaker space 
${\cal L} _{- \sigma _1,w}$ in the previous chapter, 
which led to some additional difficulties.
Furthermore, as we will see below, 
it will be necessary to relax the weight function $w$ by
a factor $m^2$. We will anticipate this loss and construct
$r(m)$ in $\ell _{1,\tilde{w}}$, with
\be
\label{l1.20}
\tilde{w} (m) := w(m) (1 + |m|) ^2.
\ee
Note that $\tilde{w}$ is again an {\em admissible weight function}.

\section{Construction of travelling waves}
\label{l2}

We make the following definitions.
\be
\label{l2.11}
[\Lambda (\beta) r](m) := \Lambda (\beta,m) r(m) , 
\mbox{ with } \Lambda (\beta,m) =
2 \cos (\beta m) + \delta_m.
\ee
\be
\label{l2.12}
\triangle r(m) := (e^{-i \beta m} - 1) r(m).
\ee
\be
\label{l2.13}
\tilde{W}(r)(m) :=
(1 - e^{i \beta m})
\frac{1}{\alpha _1} \sum_{k=2}^{\infty} \frac{\alpha _k}{k!}
(\triangle r)^{*k} (m).
\ee
Finally, define an operator $T_{\xi}$ acting on $\ell _{1, \tilde{w}}$
in the following way.
\be
\label{l2.14}
(T_{\xi} r)(m) := e^{i \xi m}  r(m).
\ee
This multiplication operator on a Fourier sequence corresponds to 
a translation of the function, an operation under which the autonomous 
differential equation is invariant.

The next lemma follows easily from Remark \ref{rb2.1}.
\lfd
\begin{lemma}
\label{ll2.51}
Suppose that $\beta \in \relz$ and $r : \ganz \rightarrow \comz$ satisfies
$\| r \| _{\ell _1} < \frac{\rhof}{8}$.
For $ n \in \ganz$ set 
\bd
x_n^{(0)} (t) :=  c n + \sum_{m \in \ganz} r(m) e^{i m (\beta n + \gamma t)}.
\ed
Then $r$ satisfies conditions (A) and(B), with 
\begin{itemize}
  \item[(A)]
      $\Lambda (\beta) r + \tilde{W} (r) = 0$, 
  \item[(B)] 
     $ \forall m \in \ganz  : 
      r(-m) = \overline{r(m)}$, 
\end{itemize}
if and only if $x_n^{(0)} (t)$ is a real valued solution of the differential
equation
\bd
\ddot{x_n} (t) =
F(x_{n-1} (t) -x_{n} (t)) -
F(x_{n} (t) -x_{n+1} (t)), \mbox{ for } n \in \ganz.
\ed
\end{lemma}

Before starting the construction of solutions of equation
(A), we have to investigate the smoothness and 
symmetry properties of the nonlinearity $\tilde{W}$.

\lfd
\begin{proposition}
\label{pl2.1}
There exists a constant $C_{F,c}$ (which can be taken to be the same as in
Proposition \ref{pb2.3}) 
such that for all $r \in \ell _{1,\tilde{w}}$ with 
$\| r \| _{\ell _{1,\tilde{w}}} \leq \frac{\rhof}{8}$ the following holds.
The series in the definition of $\tilde{W}(r)$ converges absolutely to 
an element in $\ell _{1,\tilde{w}}$. Furthermore
\begin{itemize}
\item[(i)]
  $\tilde{W}: \left\{ r \in \ell _{1,\tilde{w}}: 
  \| r \| _{\ell _{1,\tilde{w}}} \leq
  \frac{\rhof}{8} \right\} \longrightarrow \ell _{1,\tilde{w}}$ 
  is a smooth map and
  the following estimates hold.
  \bea
    &&
    \| \tilde{W} (r) \| _{\ell _{1,\tilde{w}}} \leq 
    C_{F,c} \| r \| _{\ell _{1,\tilde{w}}}^2.    \\
    \forall x \in \ell _{1,\tilde{w}} &:&
    \| D\tilde{W} (r) x \| _{\ell _{1,\tilde{w}}} \leq
    C_{F,c} \| r \| _{\ell _{1,\tilde{w}}} 
    \| x \| _{\ell _{1,\tilde{w}}} . \\
    \forall x_1, x_2 \in \ell _{1,\tilde{w}} &:&
    \| D^2\tilde{W} (r) [x_1,x_2] \| _{\ell _{1,\tilde{w}}} \leq
    C_{F,c} \| x_1 \| _{\ell _{1,\tilde{w}}} 
    \| x_2 \| _{\ell _{1,\tilde{w}}} . 
  \eea
\item[(ii)]
  If for all $m \in \ganz, r(-m) = \overline{r(m)}$, then we have
  $\tilde{W}(r)(-m) = \overline{\tilde{W}(r)(m)}$, for all $m \in \ganz$.
\item[(iii)]
  $\tilde{W}( T_{\xi} r ) = T_{\xi} \tilde{W} (r).$
\item[(iv)]
  $r(m) \in i \relz$ for all $m \in \ganz$ implies
  $\tilde{W}(r)(m) \in i \relz$ for all $m \in \ganz$.
\end{itemize}  
\end{proposition}

\lfd
\begin{remark}
\label{rl2.1}
Property (iii) reflects that the underlying equation is autonomous.
Property (iv) corresponds in the original space to the fact that if
$x_n(t)$ is a solution, then $-x_{-n} (-t)$ is again a solution.
\end{remark}

\begin{proof}
The methods for proving (i) form a proper subset of what has already 
been done in Section \ref{m23}, Proposition \ref{pb2.3}, 
and therefore these arguments are not repeated here.
The properties (ii) and (iii) can be deduced from the corresponding
properties of the convolution (see Proposition \ref{pb2.2}).
Finally we observe that $r(m) \in i \relz$ implies  
$(\triangle r)(m) \in e^{-\frac{i \beta}{2} m} \relz$. Arguing as in (iii)
we conclude for all $k \geq 1$, that
$(\triangle r)^{*k}(m) \in e^{-\frac{i \beta}{2} m} \relz$, which yields
(iv).
\end{proof}

For the construction of the solution we will use a Lyapunov-Schmidt
decomposition. We introduce the projections
\lfd
\begin{eqnarray}
P: \ell _{1,\tilde{w}} & \rightarrow & \ell _{1,\tilde{w}}  \nonumber \\
r & \mapsto & (Pr)(m) :=
\left\{
\begin{array}{ccc}
r(m), &\mbox{ for }& |m| > 1 \\
0, &\mbox{ else }&
\end{array}
\right. . \label{l2.17}
\end{eqnarray}
\be
\label{l2.20}
Q := I - P.
\ee
Denoting $\varphi := Q r, \mu := P r$, we obtain two equations:
\lfd
\begin{eqnarray}
\label{l2.25}
\Lambda (\beta) \mu + P \tilde{W} (\varphi + \mu) &=& 0  \\
\lfd \label{l2.26}
\Lambda (\beta) \varphi + Q \tilde{W} (\varphi + \mu) &=& 0 
\end{eqnarray}
For given $\varphi$ and $\beta$ we will 
be able to solve the first, infinite dimensional equation for $\mu$
(see Lemma \ref{ll2.1} below). 
Substituting this 
solution into the second equation and using various symmetries 
we end up with one equation which we can solve by choosing $\beta$
as a function of $\varphi$, with $\beta$ close to 
$\beta _1$ and $\varphi$ close to $0$. As mentioned in the introduction
of the present chapter
the spatial frequency $\beta _1$ was chosen, because 
$\Lambda (\beta _1,1) = 0$. It is obviously true that
$\Lambda (\beta _{-1},1) = 0$ (recall from (\ref{i1.25}), that 
$\beta _1 = - \beta _{-1}$) and therefore one could replace
$\beta _1$ by $\beta _{-1}$ in the construction which follows. 
However, in all the  
numerical experiments that we have performed (see e.g. Figures \ref{fc2}
and \ref{fc5})
we only observe solutions corresponding to 
$\beta _1$. This is related to the direction
in which the energy flows in the lattice and is explained in more 
detail in Remark \ref{rl2.2} at the end of this section.
\newline
Let us now parameterize $\varphi$. Observe that the 
second equation (\ref{l2.26}) is 
automatically satisfied for $m=0$
(use (\ref{l2.11}) and (\ref{l2.13})). On the other hand $(\triangle r)(0)
= 0$, independent of the value of $r(0)$, i.e. the value of $r(0)$ has 
no influence and we can normalize it to be zero. 
This freedom reflects an additional symmetry, namely if $x_n^{(0)} (t)$
is a solution of the differential equation of the doubly infinite lattice,
then $x_n^{(0)} (t)$ + constant is again a solution for any constant.
Consequently we
parameterize only a subspace of the range of $Q$, namely
\be
\label{l2.21}
\varphi (q_1,q_2) (m) :=
\left\{
\begin{array}{ccc}
\frac{1}{2} (q_1 + i q_2) , &\mbox{ for }& m = 1 \\
\frac{1}{2} (q_1 - i q_2) , &\mbox{ for }& m = -1 \\
0, &\mbox{ else }&
\end{array}
\right. \quad .
\ee
Under a slight abuse of notation we define $T_{\xi}: \relz ^2 
\rightarrow \relz ^2$, 
\be
\label{l2.22}
T_{\xi} 
\left(
\begin{array}{c}
q_1 \\ q_2
\end{array}
\right) 
:=
\left(
\begin{array}{cc}
\cos \xi & - \sin \xi \\ \sin \xi & \cos \xi
\end{array}
\right) 
\left(
\begin{array}{c}
q_1 \\ q_2
\end{array}
\right). 
\ee
We then have
\bd
\varphi (T_{\xi} q) = T_{\xi} (\varphi (q)).
\ed
Let us now solve the first equation (\ref{l2.25}).

\lfd
\begin{lemma}
\label{ll2.1}
There exists a neighborhood $U \subset \relz ^3$ of 
$(\beta _1,0,0)$, a $C^2$ map
$\mu : U \rightarrow \ell _{1,\tilde{w}} ,
(\beta,q_1,q_2) \mapsto \mu (\beta ,q)$ and a constant $C$, such that
\begin{itemize}
\item[(i)]
$\Lambda (\beta) \mu (\beta, q) + 
P \tilde{W} (\varphi (q) + \mu (\beta, q) ) = 0$.
\item[(ii)]
$Q \mu (\beta, q) = 0$.
\item[(iii)]
$\forall \beta: \mu (\beta, 0) = 0$.
\item[(iv)]
$\| \mu (\beta, q) \| _{\ell _{1,\tilde{w}}}  \leq C |q|^2$.
\item[(v)]
$\| \varphi (q) + \mu (\beta, q) \| _{\ell _{1,\tilde{w}}}  \leq 
\frac{\rhof}{8}$.
\item[(vi)]
$\mu (\beta, q) (-m) = \overline{\mu (\beta, q) (m)}$, for all $m \in \ganz$.
\item[(vii)]
$\mu (\beta, T_{\xi} q ) = T_{\xi} \mu (\beta ,q)$. 
\item[(viii)]
$\mu (\beta, 0 , q_2) \in i \relz$.
\end{itemize}
\end{lemma}
\begin{proof}
We show first that the existence of the function 
$\mu$ can be deduced from the implicit
function theorem applied to the map
\bea
G: \tilde{U} \times 
\left\{ \mu \in \mbox{ Ran }(P) : \| \mu \| _{\ell _{1,\tilde{w}}} 
\leq \frac{\rhof}{16} \right\}   &\rightarrow& \mbox{ Ran }(P),  \\
(\beta, q, \mu) &\mapsto& 
\mu + \Lambda (\beta)^{-1} P \tilde{W} (\varphi (q) + \mu).
\eea
Indeed, note that $\Lambda (\beta)^{-1}$ is well defined on Ran$(P)$, as 
$\Lambda(\beta,m) \neq 0$ for $|m| > 1$ (follows from condition
(\ref{l1.6}) in the introduction of the present chapter). 
Furthermore it is easy to check
that $\Lambda ^{-1}$ is a $C^2$ function of $\beta$ in the operator norm,
which shows together with Proposition \ref{pl2.1} that $G$
is a $C^2$ map, if $\tilde{U}$ is chosen as a suitably small neighborhood
of $(\beta _1,0)$ in $\relz ^3$.
Using Proposition \ref{pl2.1} again we see, that
\bd
G(\beta _1,0,0) = 0,
\ed
and
\bd
D_{\mu} G (\beta _1,0,0) = \mbox{ id }_{\mbox{Ran}(P)}.
\ed
which allows us to apply the implicit function theorem and to obtain 
properties (i) and (ii) immediately. (iii) follows from the uniqueness
of the implicit function theorem and the observation, that $\forall \beta:
G(\beta,0,0) =0$. In order to show the remaining properties we will 
investigate the family of maps
\bd
H_{\beta ,q}: \mu \mapsto 
-\Lambda (\beta)^{-1} P \tilde{W}(\varphi (q) + \mu).
\ed
Note that a fixed-point of $H_{\beta ,q}$ corresponds to a zero of $G$.
That $H_{\beta ,q}$ is a contraction for suitably small values of
$q$ and $\mu$ is a consequence of Proposition \ref{pl2.1}.
In order to prove (iv) it suffices to show that $H_{\beta ,q}$ maps
$\left\{ \mu \in \mbox{ Ran }(P) : \| \mu \| _{\ell _{1,\tilde{w}}}
\leq C |q|^2 \right\}$ into itself for a constant $C$, which will 
be chosen below. In order to verify this claim,
let us denote 
\bd
C_{\Lambda} := \sup_{|m| > 1} \frac{1}{-2 + \delta _m} =
\frac{1}{-2 + \delta _2}.
\ed
Proposition \ref{pl2.1} implies
\bea
\| H_{\beta ,q} (\mu) \| _{\ell _{1,\tilde{w}}}
& \leq & C_{\Lambda} C_{F,c} 
2 \left(
\left( \frac{\tilde{w}(1)+\tilde{w}(-1)}{2} \right) ^2 
|q|^2 + \| \mu \| _{\ell _{1,\tilde{w}}} ^2 \right) \\
&\leq &
C |q|^2,
\eea
if we choose
\bd
C := 4 C_{\Lambda} C_{F,c} 
\left( \frac{\tilde{w}(1)+\tilde{w}(-1)}{2} \right) ^2 ,
\ed
and we only allow $q \in \relz ^2$ such that 
\bd
2 C_{\Lambda} C_{F,c} C |q|^2 \leq \frac{1}{2}.
\ed
Property (v) follows trivially from (iv) by making $|q|$ sufficiently
small.\newline
The last three properties can be proved in the following way.
$\mu(\beta, q)$ is the fixed-point of $H_{\beta ,q}$ and can be constructed
as the limit of the sequence $\mu _k (\beta, q)$ as $k
\rightarrow \infty$. We define inductively 
\bea
\mu _0 (\beta, q) &:=& 0,    \\
\mu _{k+1} (\beta, q) &:=&
H_{\beta ,q} (\mu _k (\beta, q)).
\eea
(vi),(vii),(viii) are trivially satisfied for $\mu _0$ as well as for the  
sequence $\varphi (q)$. Proposition \ref{pl2.1} permits us to conclude
inductively that all $\mu _k$ possess the three properties, which then is 
preserved under the limit $k \rightarrow \infty$. 
\end{proof}

We turn now to the finite dimensional equation (\ref{l2.26}).
\bd
\Lambda (\beta) \varphi (q)+ 
Q \tilde{W} (\varphi (q)+ \mu (\beta ,q)) = 0.
\ed
Let us first reduce the dimensions of these equations by factoring out 
all the symmetries. We recall that the
equation is satisfied for $m=0$ by (\ref{l2.11}) and (\ref{l2.13}).
Furthermore it suffices to solve equation (\ref{l2.26}) for $m=1$ as
the equation for $m = -1$
only yields the complex conjugate of the same equation.
This leaves us with one complex valued or 
equivalently with two real valued equations. Finally we can make one further 
reduction using the $T _{\xi}$ invariance. It implies that a choice
of $(\beta ,q_1 ,q_2)$ is a solution if and only if 
$(\beta , 0, |q|)$ is a solution
($|q| := \sqrt{q_1^2 +q_2^2}$). But we have seen that for  
arguments of the form $(\beta , 0, |q|)$, 
all the terms in equation (\ref{l2.26}) are purely imaginary
and hence we have reduced 
equation (\ref{l2.26}) to one equation in two unknowns.
Accordingly define
\bea
g(\beta, p) &:=& \frac{1}{i} \left(
\Lambda (\beta ,1) \varphi (0,p)(1) +
\tilde{W}(\varphi (0,p) + \mu (\beta , 0, p)) (1) \right) \\
&=&
(2 \cos \beta + \delta _1) \frac{p}{2} +
\frac{1}{i} \tilde{W}(\varphi (0,p) + \mu (\beta , 0, p)) (1).
\eea
We want to identify the set of zeros of the $C^2$ function $g$.
Property (iii) of the last lemma implies that $g(\beta,0) = 0$
for all $\beta$. We introduce the associated function
\bd
\tilde{g}(\beta , p) := \frac{g(\beta ,p)}{p}.
\ed 
Note that this function has a $C^1$ extension for $p=0$. Furthermore
\bea
\tilde{g}(\beta _1 , 0) &=&  \partial _p g(\beta _1 ,0)  \\
&=&
(2 \cos \beta _1+ \delta _1) \frac{1}{2} +
\frac{1}{i} D \tilde{W}(0)[\partial _p (\varphi + \mu)](1) = 0.
\eea
We compute
\bea
\partial _{\beta} \tilde{g}(\beta _1 , 0)  &=&
\partial _{\beta} \partial _p g(\beta _1 ,0)  \\
&=&
- \sin \beta_1 +
\frac{1}{i} D \tilde{W}(0)[\partial _{\beta , p} \mu](1) +
\frac{1}{i} D^2 \tilde{W}(0)
[\partial _p (\varphi + \mu),\partial _{\beta} \mu](1)  \\
&=&
- \sin \beta_1  \neq 0.
\eea
In the above calculation we used that $D \tilde{W} (0) = 0$ and  that
$\partial _{\beta} \mu (\beta _1 ,0) = 0$, as 
$\mu (\beta ,0) = 0$ for all $\beta$.
Hence we can again apply the implicit function theorem and obtain in this
way an even $C^1$ function $\beta (p)$, defind as map from a neighborhood
of zero to a neighborhood of $\beta _1$, such that
\bd
g(\beta (p), p) = 0.
\ed
The evenness of $\beta$ is a consequence of the uniqueness in the implicit
function theorem and the evenness of $\tilde{g}$ in $p$.
We can summarize our considerations.

\lfd
\begin{lemma}
\label{ll2.2}
Given the function $\mu (\beta ,q)$ from Lemma \ref{ll2.1},
there exists an even $C^1$ function $\beta (q)$, mapping from
a neighborhood of zero to a neighborhood of $\beta _1$ such that 
the set of zeros of
\bd
\Lambda (\beta) \varphi (q)+ 
Q \tilde{W} (\varphi (q)+ \mu (\beta ,q)) = 0
\ed
in a neighborhood $U_0$ of $(\beta _1,0,0)$ is given by
\bd
\{(\beta ,q) \in U_0 : q = 0 \} \cup 
\{(\beta ,q) \in U_0 : \beta = \beta(\sqrt{q_1^2 + q_2^2}) \}.
\ed
\end{lemma}

Combining the last two lemmas we have proved the following theorem.
\lfd
\begin{theorem}
\label{tl2.1}
There exists a neighborhood $D_0$ of $0$ in $\relz ^2$, a 
$C^1$ map 
\bd
r:D_0 \rightarrow \ell _{1,\tilde{w}}, (q_1 , q_2)
\mapsto r(q_1, q_2) 
\ed
and a $C^1$ map
\bd
\tilde{\beta}: D_0 \rightarrow \relz, 
(q_1, q_2) \mapsto \tilde{\beta} (q_1, q_2),
\ed
such that 
\begin{itemize}
\item[(i)]
$\forall q \in D_0: 
\Lambda (\tilde{\beta} (q)) r(q) + \tilde{W} (r(q)) = 0.$
\item[(ii)]
$r(q)(-m) = \overline{r(q)(m)}$ for all $m \in \ganz, q \in D_0$.
\item[(iii)]
$\forall q \in D_0:
\| r(q) \| _{\ell _{1,\tilde{w}}} \leq \frac{\rhof}{8}.$
\item[(iv)]
$ \exists C > 0 \forall q \in D_0 :
\sum_{|m|>1} \tilde{w}(m) |r(q)(m)| \leq C |q|^2.$
\end{itemize}
\end{theorem}

\begin{proof}
Let 
\bea
\tilde{\beta} (q) &:=& \beta (|q|),   \\
r(q) &:=& \varphi(q) + \mu(\tilde{\beta}(q),q). 
\eea
Then all the properties follow immediately from the last two lemmas.
The differentiability of $\tilde{\beta}$ at $0$ follows from the fact 
that $\beta (\cdot)$
is an even $C^1$ function.
\end{proof}
As we have proved Lemma \ref{ll2.1} and Lemma \ref{ll2.2} using implicit
function theorems, we can prove a uniqueness result for the 
constructed sequences $r(q)$.
\lfd
\begin{corollary}
\label{cl2.1}
There exists a $\delta > 0$, such that the following holds.
Let $s \in \ell _{1, \tilde{w}}$ and $\beta \in \relz$, satisfying
\begin{itemize}
\item[(I)]
$0 < \| s \| _{1, \tilde{w}} < \delta.$
\item[(II)]
$s(0) = 0.$
\item[(III)]
$s(-m) = \overline{s(m)}$, for all $m \in \ganz$.
\item[(IV)]
$|\beta - \beta _1| < \delta$.
\end{itemize}
If 
\bd
y_n (t) := c n + \sum_{m \in \ganz} s(m) \exp (i m (\beta n + \gamma t)),\quad
 n \in \ganz,
\ed
is a solution of
\bd
\ddot{y} _n = F(y_{n-1} - y_n) - F(y_n - y_{n+1}),\quad n \in \ganz,
\ed
then $s = r(q)$ and $\beta = \tilde{\beta} (q)$, with
\bd
q = (q_1, q_2) = 2 (\mbox{ Re}(s(1)), \mbox{ Im}(s(1)) ),
\ed
and $r, \tilde{\beta}$ are the functions defined in Theorem \ref{tl2.1} above.
\end{corollary}
\begin{proof}
If $\delta$ is chosen suitably small, it follows from Lemma
\ref{ll2.51} above
that
\bd
\Lambda(\beta) s + \tilde{W} (s) = 0.
\ed 
Defining $q$ as above, we see that the
projection of $s$, $Ps$,
satisfies assumptions (i) and (ii) of Lemma \ref{ll2.1}. Consequently
we have $Ps = \mu (\beta, q)$. 
Furthermore the spatial frequency
$\beta$ satisfies the equation in Lemma \ref{ll2.2} with $q \neq 0$
(otherwise $\mu = s = 0$, which contradicts assumption (I)). Lemma
\ref{ll2.2} yields that $\beta = \tilde{\beta} (q)$  
which finally implies $s = r(q)$ and the 
Corollary is proved.
\end{proof}
\lfd
\begin{remark}
\label{rl2.2}
Heuristic explanation of the choice of the sign of $\beta$.
\end{remark}
We investigate the transport of energy of the travelling wave solutions
constructed in Theorem \ref{tl2.1}. The amount of energy exchanged
between the particles $x_0$ and $x_1$ is given by
\bd
E = - \int_{0}^{\frac{2 \pi}{\gamma}} 
F(x_0(t) - x_1(t)) \dot{x} _0(t) dt.
\ed
Note that $E < 0$ means that on the average energy flows
from $x_0$ to $x_1$. Recall from equation (\ref{l1.10}) in
the introduction of this chapter, that 
\bd
x_n^{(0)}(t) = c n + \sum_{m \in \ganz} 
r(m) e^{i \beta m n} e^{i m \gamma t}.
\ed
Let
\bea
\tilde{r}(m) &:=& (1 - e^{i \beta m}) r(m).  \\
r' (m) &:=& i \gamma m r(m).
\eea
Then
\bd
E = - \frac{2 \pi}{\gamma}
\left[
\sum_{k=0}^{\infty} \frac{\alpha _k}{k!} \tilde{r} ^{*k} * r'
\right] (0).
\ed
Note from Lemma \ref{ll2.1} (iv) that only 
$r(1)$ and $r(-1)$ are of first order in $q$, all others
are of higher order. We now evaluate the formula above.
\begin{itemize}
\item
$1^{\mbox{ st }}$ order:
$ - \frac{2 \pi}{\gamma} \alpha _0 r'(0) = 0.$
\item
$2^{\mbox{ nd }}$ order:
\bea
- \frac{2 \pi}{\gamma} \alpha _1 (\tilde{r} * r')(0) &\simeq&
- \frac{2 \pi}{\gamma} \alpha _1
(\tilde{r}(1) r'(-1) + \tilde{r}(-1) r'(1)) \\
&=&
- \frac{4 \pi}{\gamma} \alpha _1 \mbox{ Re } \left[
(1 - e^{i \beta}) |r(1)|^2 (-i \gamma) \right] \\
&=&
4 \pi \alpha _1 |r(1)|^2 \sin \beta.
\eea
\item
$3^{\mbox{ rd }}$ order:
no terms $\neq 0$.
\end{itemize}
In the system we are investigating we expect, that the driver excites
outgoing waves, hence the energy should be transported in direction of
increasing $n$. The above calculation shows up to third order in $q$,
that this this is achieved for $\beta$ with $\sin \beta < 0$. Therefore
we have chosen solutions with $\beta$ close to $\beta _1$ rather than
close to $\beta _{-1} = - \beta _1$.

\section{Construction of the periodic solution via Chapter~3}
\label{l3}

The following theorem is the main result of the present chapter.
\lfd
\begin{theorem}
\label{tl3.1}
Let $F,c, (b_m) _{m \in \ganz}, w$ satisfy the general assumptions and 
let $F'(-c) < \gamma ^2 < 4 F'(-c)$. Then there exists a neighborhood
$D$ of $0$ such that for all $\epsilon \in D$ there exist sequences
$u(\epsilon) \in {\cal L} _{0,w}; v(\epsilon) \in 
{\cal L} _{\sigma _0, w}$ ($\sigma _0 > 0$ defined in Definition 
\ref{db3.1})
such that
\bd
x_n(t) := c n + \epsilon b_0 - (u + v)(\epsilon)(0,0)
+ \sum_{m \in \ganz} (u+v)(\epsilon)(n,m) e^{i m \gamma t},\quad
\mbox{ for } n \geq 1,
\ed
is a time periodic solution of the differential equation given by 
(\ref{i2.25}) and (\ref{i2.26}).
\end{theorem}

\begin{proof}
We start the construction with $r(q) \in \ell _{1, \tilde{w}}$, which was 
obtained in Theorem \ref{tl2.1}. Recall the definition in equation
(\ref{l1.20}), namely
$\tilde{w} (m) := w(m) (1 + |m|)^2$, which again is an admissible 
weight function. According to the 
ansatz which was explained in the beginning of the present chapter, we
define
\be
\label{l3.10}
u(q)(n,m) := r(q)(m) \exp (i \tilde{\beta} (q) m n),\quad \mbox{ for } q \in D_0.
\ee
From Theorem \ref{tl2.1} it follows that all of the conditions of 
Theorem \ref{tb3.1} on the function $u(q)$ are trivially satisfied,
with the exception of condition (2), which states the differentiability
with respect to the parameter $q$ in a certain norm. 
It is obvious from Theorem \ref{tl2.1} that each component of $u$ is 
a $C^1$ function of $q$ with
\bd
\frac{\partial}{\partial q_j} u(q) (n,m) =
\left(\frac{\partial}{\partial q_j} r(q) (m)+ 
i r(q)(m) \frac{\partial}{\partial q_j} \tilde{\beta} (q) m n  \right)
\exp (i \tilde{\beta} (q) m n).
\ed
In order to show the $C^1$ dependence of $u$ on $q$ in the
${\cal L} _{-\sigma _1, w}$ norm, we only have to verify that
$q \mapsto \frac{\partial}{\partial q_j} u(q)$ is a continuous map
from $D_0$ into ${\cal L} _{-\sigma _1, w}$. It is straightforward
to check this by hand, using the following observations.
\begin{itemize}
\item
 \bea
 &&|\exp (i \tilde{\beta} (q') m n) - \exp (i \tilde{\beta} (q) m n)| \\
 &\leq& m n |q' - q| \sup_{s \in [0,1]}  \left|
 \frac{\partial}{\partial q_j} \tilde{\beta} (q + s(q' - q))
 \right| .
 \eea
\item
 $w(m) m^2 \leq \tilde{w}(m).$
\item
 $ \sup_{n \geq 0} n^2 e^{- \sigma _1 n} < \infty$.
\end{itemize}

Therefore we are in a position to apply Theorem \ref{tb3.1} of
Chapter \ref{m} and obtain this way the function $v(q, \epsilon)$.
Equipped with both, $u$ and $v$, we turn now to Lemma \ref{lb2.1}
of Chapter \ref{m}.
Using Lemma \ref{ll2.51} it only remains to solve
\be
\epsilon b_m - v(0,m) - u(0,m) = 0, \mbox{ for } |m| = 1,
\label{l3.105}
\ee
as by the conditions on $\gamma$ we have $m_0 = 1$.
It suffices to solve equation (\ref{l3.105}) for $m=1$,
because for $m=-1$ we obtain the complex conjugate 
of the equation. Hence there are two real equations 
to be solved (real part
and imaginary part). Let
\bd
g(q, \epsilon) := \epsilon b_1 - v(q, \epsilon)(0,1) - u(q)(0,1)
\in \relz ^2.
\ed
Note that $v(0,1)$ is a $C^1$ function of $(q, \epsilon)$ (see 
Theorem \ref{tb3.1} (v))
with $|v(q, \epsilon)(0,1)| \leq 2(|\epsilon| + C |q|^2)$ 
(see Theorem \ref{tb3.1} (ii), Lemma \ref{ll2.1} (iv)) and therefore 
$D_q v(0,0) (0,1) = 0$. Furthermore $u(q)(0,1) = \frac{1}{2}
(q_1 , q_2)$. We conclude
\bea
g(0,0) &=& 0, \\
(D_q g) (0,0) &=& \frac{1}{2} 
\left(
\begin{array}{cc} 1 & 0 \\ 0 & 1 \end{array}
\right).
\eea
The implicit function theorem then guarantees the 
existence of a neighborhood $D$ of
$0$ and a $C^1$ function $q(\epsilon)$ such that 
$g(q(\epsilon),\epsilon ) = 0$. This proves that for
\bea
u(\epsilon) &:=& u(q(\epsilon)), \\
v(\epsilon) &:=& v(q(\epsilon), \epsilon),
\eea
the hypothesis of Lemma \ref{lb2.1} is satisfied and we have 
successfully constructed the solution to the differential equation
given by (\ref{i2.25}) and (\ref{i2.26}).
\end{proof}

\chapter{The Toda Lattice}
\label{g}

In this chapter we will use the complete integrability of the doubly
infinite Toda lattice ($F(x) = e^x$) and show how the well known
$g$ -gap solutions contain a sufficiently large family of travelling waves
to construct solutions of equations (\ref{i2.25}) and (\ref{i2.26})
for any number of resonances $m_0$.

$G$-gap solutions were first constructed for the continuous analog
of the Toda lattice, the KdV equation. Combining the complete integrability
of the system as given in (\cite{Lax1}) with methods,
developed in algebraic geometry (see \cite{Akh}), these solutions 
can be expressed in terms of ratios of theta functions (see \cite{Dub},
\cite{IM}). Although these solutions have been
studied and used in many different contexts (see e.g.\cite{DT}, \cite{Kri1}, 
\cite{Kri2},\cite{McK}, etc.), we repeat their construction in the Toda case, 
because we will need some specific details about 
these solutions, which are not readily available in the literature.
The resulting formulae will be evaluated in the case of small gaps
and we derive a $C^1$ parameterization of time-periodic 
$g$-gap solutions. This will be done in Section \ref{g1}.
In order to keep the presentation self contained we provide
the details of the construction in the Appendix \ref{g2}, where we will
use only some general facts about hyperelliptic curves.
The resulting parameters are 
closely related to those introduced first by T. Kappeler (\cite{Kap1}) for
spatially periodic potentials in the case of KdV, then by T. Kappeler et al.
\cite{BGGK} in the Toda case.
It turns out that the condition of time periodicity determines
the position of the midpoints of the gaps. In fact, there are two possible
positions for each gap, corresponding to outgoing and incoming waves.
We then proceed to 
prove in Section \ref{g4} that the basic result of Chapter \ref{m}
can be applied to obtain periodic solutions of the driven lattice 
for an arbitrary number of resonances. 
In the case of $m_0 \leq 1$ these solutions are shown to coincide with 
those constructed for general lattices in the previous two chapters.
Finally the essential spectrum of the corresponding
Lax operator is investigated. Clearly it has a band -- gap structure and 
we will determine the
width of the gaps and their dependence on the 
Fourier coefficients of the driver to first order in $\epsilon$.

\section{The $g$-gap solutions}
\label{g1}

First we briefly describe the construction of $g$-gap solutions
via {\em Baker-Akhiezer functions}. We follow the construction in 
\cite{Kri1} (see also \cite{Akh},\cite{Dub}, \cite{IM}).

Let $g \in \natz _0 = \{ n \in \ganz: n \geq 0 \}$ and denote by 
$R_g$ the hyperelliptic curve of
genus g, which is constructed by pasting together two copies of the
Riemann sphere $\comz \cup \{ \infty \}$ along the slits 
$[E_0,E_1], [E_2,E_3], \ldots, [E_{2g},E_{2g+1}]$ ,where 
$(E_0 < E_1 < \cdots < E_{2g+1} ).$
Denote by $\pi$ the canonical projection of $R_g$ on
the Riemann sphere. We fix $g$ points $P_j \in R_g, 1 \leq j \leq g$,
satisfying $ \pi (P_j) \in [E_{2j-1},E_{2j}]$.

Then for each $n \in \ganz, t \in \relz$, there exists an unique
{\em Baker-Akhiezer function $\psi _n(t,\cdot)$} 
(up to multiplication by a constant)
meromorphic on
$R_g \setminus \{ P_{\infty} , P_{\infty}^* \}$, with at most simple poles
at $P_j, 1 \leq j \leq g$ and a certain prescribed 
behavior at the essential singularities $ P_{\infty}$ and 
$P_{\infty}^*$, depending on $n$ and $t$. See Theorem \ref{tg2.1} in
Appendix \ref{g2}.
for the precise statement. The existence of these functions is proved 
by explicitly constructing them in terms of theta functions and
uniqueness is a consequence of the Riemann-Roch theorem applied 
to $R_g$.
Using the defining properties of the {\em Baker-Akhiezer functions}
one is able to obtain functions $x_n(t)$ such that the corresponding
operators $\tilde{L}$ and $\tilde{B}$ (see (\ref{e.65}) and (\ref{e.70})) 
satisfy for all $t \in \relz$
and $P \in R_g$ the following equations.
\be
\label{g2.5}
(\tilde{L} \psi)_n (t,P) = \frac{1}{2} \pi(P) \psi _n (t,P).
\ee
\be
\label{g2.6}
\frac{\partial}{\partial t} \psi _n (t,P) =(\tilde{B} \psi)_n (t,P).
\ee
From equations (\ref{g2.5}), (\ref{g2.6}) we conclude for $t \in \relz$, that
\be
\label{g2.4}
\tilde{L}_t = [\tilde{B},\tilde{L}].
\ee
Hence $x_n(t)$ is a solution of the doubly infinite Toda chain. Furthermore
$x_n(t)$ can be expressed in terms of theta functions, namely
\be
\label{g1.5}
x_n(t) = n I + t R + f_{\tau} \left(
Un + Vt - Z \right).
\ee
The function $f_{\tau} : \relz ^g / \ganz ^g \longrightarrow \relz$
is essentially given by the logarithm of a ratio of theta functions
with period matrix $\tau$. In general these solutions are quasiperiodic
in $n$ and $t$. The parameters $I, R, \tau, U, V, Z$
depend on the spectral data chosen in the construction, i.e.
they depend on the $3g + 2$ real parameters $E_i, 0 \leq i \leq 2g + 1$ and 
$P_j, 1 \leq j \leq g$, or equivalently on $(a, b, \lambda _j, p_j, P_j,
1 \leq j \leq g)$, where $\lambda _j$ denotes
the midpoint and $p_j$ the half-width of the $j$-th gap.
We are interested in obtaining the following choice for the parameters.
\lfd
\begin{eqnarray}
I &=& c.
\label{g1.15}  \\
R &=& 0.
\lfd \label{g1.20}  \\
V &=& 
-\frac{\gamma}{2 \pi}
\scriptsize
\left(
\begin{array}{c}
1 \\
\vdots \\
g \\
\end{array}
\right).
\normalsize
\lfd \label{g1.25} 
\end{eqnarray}
The next theorem will show that we can choose the parameters $a,b,
\lambda _j, 1 \leq j \leq g$ as $C^1$ functions of the remaining $2g$
parameters $p_j, P_j, 1 \leq j \leq g$ such that equations 
(\ref{g1.15}), (\ref{g1.20}), (\ref{g1.25}) are satisfied.

\lfd
\begin{theorem}
\label{tg1.5}

We assume that $g \in \natz _0, \gamma > 0, c \in \relz$ satisfy 
\be
g \gamma < 2 e^{- \frac{c}{2}} < (g+ 1) \gamma.
\label{g1.30}
\ee 
Then there exists a positive number $\delta$ and $C^1$ functions
$a, b, \lambda _j, U _j, \tau^{(reg)}_{i,j},
1 \leq i,j \leq g$, mapping
$\{(p_1, \ldots, p_g) \in \relz ^g: |p_k| < \delta, 1 \leq k \leq g \}$
into $\relz$ and which are even in each argument 
$p_k, 1 \leq k \leq g$ such that 
the following holds.

For $0 < p_1, \ldots, p_g < \delta$ the $g$ -gap solution corresponding
to the choice of parameters 
\bd
p_j, P_j, a(p_1, \ldots, p_g), b(p_1, \ldots, p_g), 
\lambda _j (p_1, \ldots, p_g) ,1 \leq j \leq g,
\ed
is given by
\be
\label{g3.110}
x_n (t) =
c n + \ln
\frac{\vartheta (\frac{1}{2} U - Z |\tau)
      \vartheta ((n-\frac{1}{2}) U + tV - Z |\tau)}
     {\vartheta (- \frac{1}{2} U - Z |\tau)
      \vartheta ((n+\frac{1}{2}) U + tV - Z |\tau)}
\ee
with
\begin{itemize}
\item[(i)]
$V = 
-\frac{\gamma}{2 \pi}
\scriptsize
\left(
\begin{array}{c}
1 \\
\vdots \\
g \\
\end{array}
\right).
\normalsize$
\item[(ii)]
The map $(P_1,\ldots,P_g) \longmapsto Z$ is a surjective map
from $(S^1)^g$ to $\relz ^g / \ganz ^g$
for all choices of parameters $0 < p_1, \ldots, p_g < \delta$.
\item[(iii)]
\be
\label{g3.62}
\pi i \tau =
\mbox{ diag } (\ln p_k) + \tau^{(reg)} .
\ee
\end{itemize}
Furthermore formula (\ref{g3.110}) is also well defined for
$p_j \geq 0$ and the function we obtain by letting 
some or all of the $p_j$
converge to $0$ agrees with the corresponding lower gap solution.
\end{theorem}

\lfd
\begin{remark}
\label{rg3.2}
\begin{em}
\begin{itemize}
\item[(1)]
The proof of Theorem \ref{tg1.5} is given in the Appendix \ref{g3}.
\item[(2)]
In fact there are two possible choices for each function $\lambda _j$
The special choice made in the proof of Theorem
\ref{tg1.5} corresponds to the numerical observation,
that gaps open up only in the lower half of the band. 
We will see in 
Remark \ref{rg3.3} at the end of the next section, that the physical
reason for this lies in the direction in which energy is transported in
the corresponding $g$ -gap solution. The reader may recall that
exactly the same situation occurred in Chapter \ref{l} with the choice
of the spatial frequency $\beta$ (see Remark \ref{rl2.2}).
\item[(3)]
The functions in Theorem \ref{tg1.5} are also defined for negative
values of the $p_k$'s. This extension is purely formal and is used
to simplify regularity proofs at $p_k = 0$.
\end{itemize}
\end{em}
\end{remark}

\section{Construction of the periodic solution via Chapter~3}
\label{g4}

In this section we use the $g$-gap solutions of the last section 
to construct periodic solutions of the driven lattice by means of the 
procedure in Chapter \ref{m}.
Our goal is to prove the following
theorem.
\lfd
\begin{theorem}
\label{tg4.1}
Let $c,\gamma ,(b_m) _{m \in \ganz}, w$ satisfy the general assumptions. Then
there exists a neighborhood $D$ of $0$, such that for all 
$\epsilon \in D$, there exist sequences 
$u(\epsilon) \in {\cal L} _{0,w}, v(\epsilon) \in
{\cal L} _{\sigma _0,w}$  (where $\sigma _0 > 0$ was introduced in Definition
\ref{db3.1}) and
\be
\hspace{.7in}
x_n(t) := c n + \epsilon b_0 - (u + v)(\epsilon)(0,0) +
\sum_{m \in \ganz} (u + v)(\epsilon)(n,m) e^{i m \gamma t} , \mbox{ for }
n \geq 1,
\label{g4.3}
\ee
is a time periodic solution of the differential equation
given by (\ref{i2.25}) and (\ref{i2.26}).
\end{theorem}

\begin{proof}
By the {\it general assumptions on $\gamma$},
there exists a $m_0 \in \natz _0$,
such that 
\bd
\frac{(m_0 \gamma) ^2}{e^{-c}} < 4 <
\frac{((m_0 +1) \gamma) ^2}{e^{-c}}.
\ed
We choose $g := m_0$ and thus satisfy the assumptions (\ref{g1.30})
of Theorem \ref{tg1.5}.
To put the results of the preceeding sections in a form that is 
suitable for the procedure of Chapter \ref{m}, we still have to
make some technical definitions and remarks.
Let us begin with the definition of the parameters. Denote
\be
\label{g4.5}
\tilde{p} _j :=  p_j \exp (2 \pi i Z _j) \in \comz,\quad  1 \leq j \leq g.
\ee
The $2g$ real variables $q_j$ are then defined by 
\be
\label{g4.6}
\tilde{p} _j = q_{2j-1} + i q_{2j}.
\ee

Note that any choice of $q$ in a sufficiently small neighborhood
of $0$ in $\relz ^{2g}$ corresponds to a choice of spectral data.
In fact, let
$p_j := |q_{2j - 1} + i q_{2j}|$. Furthermore Lemma \ref{lg3.6} in 
Appendix \ref{g3} 
shows that for
any given choice of phase $Z_j \in [0,1)$, there is
a choice of points $P_j$ which corresponds to the phase .
Using (\ref{g3.115}) we can now write equation (\ref{g3.110}) 
in the following form,
\be
\label{g4.10}
x_n^{(0)} (t) = c n + \ln \frac{1 + \Gamma _n(t,q)}
{1 + \Gamma _{n+1}(t,q)} + \ln \frac{1 + \Gamma _1(0,q)}
{1 + \Gamma _0(0,q)} , 
\ee
with
\lfd
\begin{eqnarray}
\label{g4.11}
\Gamma _n(t,q)  &:=&
\sum_{l \in {\bf z} ^g \setminus \{0\}} 
r(q)(n,l) \exp \left(-i \left(
l_1 + 2 l_2 + \cdots + g l_g \right)  \gamma t \right), \\
\lfd \label{g4.12}
&=& \sum_{m \in {\bf Z}} s(q)(n,m) e^{i \gamma m t}, \mbox{ where } \\
r(q)(n,l) &:=&
\left( \prod_{j=1}^g 
\overline{\tilde{p}_j^{l_j} p_j^{l_j(l_j - 1)}} \right) \nonumber \\
&&
\exp \left( 2 \pi i <l, U(p)>(n - \frac{1}{2}) 
+ <l, \tau^{(reg)}(p) l> \right),\lfd \label{g4.13} \\
\lfd \label{g4.14}
s(q)(n,m) &:=&
\sum_{
\scriptsize
\begin{array}{c} 
l \in \ganz ^g \setminus \{0\} \\
l_1 + \cdots g l_g = -m
\end{array}
} 
\normalsize
r(q)(n,l).
\end{eqnarray}

The Fourier series of the $g$-gap solution is now given by 
\lfd
\begin{eqnarray}
\label{g4.15}
x_n^{(0)} (t) &=&
c n + \sum_{m \in \ganz} u(q)(n,m) e^{i m \gamma t}, \mbox{ with } \\
\lfd \label{g4.16}
u(q)(n,m) &:=&
u_1(q)(n,m) - u_1(q)(n+1,m) + u_2(q) {\bf 1}_{ \{ m= 0 \} },  \\
\lfd \label{g4.17}
u_1(q) (n,m) &:=&
\sum_{k=1}^{\infty} \frac{(-1)^{k-1}}{k} s(q)(n,\cdot) ^{*k} (m), \\
\lfd \label{g4.18}
u_2 (q) &:=&
\ln ( 1 + 
\sum_{l \in {\bf z} ^g \setminus \{0\}} r(q)(1,l) ) -
\ln ( 1 + 
\sum_{l \in {\bf z} ^g \setminus \{0\}} r(q)(0,l) ).
\end{eqnarray}
The convergence of the above series follows from the smallness of 
$p_1, \ldots, p_g$ and will be verified below.
\newline
{\bf Claim:}
\newline
There exists a neighborhood $D_0$ of $0$ in $\relz ^{2g}$ and a constant
$C > 0$, such that
\begin{itemize}
\item[(I)]
for all $q \in D_0: u(q) \in {\cal L} _{0,w}$ and 
$\| u(q) \| _{0,w} \leq C |q|$.
\item[(II)]
for all $q \in D_0: \sum_{|m| > m_0} w(m)|u(q)(0,m)| \leq C|q|^2.$
\item[(III)]
$u : D_0 \rightarrow {\cal L} _{-\sigma _1,w} , q \mapsto u(q)$, is a 
$C^1$ map, where $\sigma _1$ was introduced in Definition \ref{db3.1}.
\item[(IV)]
$f:D_0 \rightarrow \relz ^{2g}, q \mapsto (u(q)(0,m))_{m=1}^g$ is  
continuously differentiable and 
\bd
\det (D_q f(0)) \neq 0.
\ed
\item[(V)]
$x_n^{(0)} (t), n\in \ganz$ is a smooth, real valued solution of the 
doubly infinite Toda lattice.
\end{itemize}

Before we proceed to check all these properties, let us show that 
they suffice to prove the theorem. 
First we have to verify that $u(q)$ satisfies the conditions (1)-(3)
in Theorem \ref{tb3.1}, but this is an immediate consequence
of (I),(III) and (V). Hence we obtain a $v(q,\epsilon) \in 
{\cal L} _{\sigma _0,w}$, satisfying conditions (i) to (v) of the
Theorem \ref{tb3.1}. Using (II) we obtain in addition, that
\be
\label{g4.20}
\| v(q,\epsilon) \| _{\sigma _0,w} \leq 2(|\epsilon| + C|q|^2).
\ee
We now show that it is possible to satisfy
(1), (2) and (3) in Lemma \ref{lb2.1}.\newline
(1) follows for $u(q)$ from (V) and the Remark  \ref{rb2.1}. \newline
(2) follows for $v(q,\epsilon)$ from Theorem \ref{tb3.1}.
\newline
By Theorem \ref{tb3.1} it suffices to
solve (3) for $1 \leq m \leq g$. This gives $2g$ real equations, 
for which we have the $2g$ real variables $q_j$ available.
Define 
\bd 
g(q,\epsilon) := (\epsilon b_m - u(q)(0,m) - v(q,\epsilon)(0,m))_{m=1}^g
\in \relz ^{2g}.
\ed
Using Theorem \ref{tb3.1} (v), equation (\ref{g4.20}) and 
what we know about function $f$ from claim (IV), we see that $g$ is a $C^1$
function and
\bea
g(0,0) &=& 0. \\
D_q g(0,0) &=& - D_q f(0),
\eea
which is an invertible matrix by (IV). 
The implicit
function theorem yields a function $q(\epsilon)$, such that 
$g(q(\epsilon),\epsilon) = 0$. Thus
\bea
u(\epsilon) &=& u(q(\epsilon)). \\
v(\epsilon) &=& v(q(\epsilon),\epsilon). 
\eea
satisfies all the conditions of Lemma \ref{lb2.1}.

It remains to prove properties (I)-(V). \newline \newline
{\bf (I):} \newline
It is easy to read off the following estimates from the  above definitions.
Denote for $l \in \ganz ^g: |l| := \sqrt{l_1^2 + \cdots + l_g^2}$.
\bd
\exists C, \delta > 0 : \forall |q| \leq \delta, 
n \in \natz _0, l \in \ganz ^g: |r(q)(n,l)| \leq 
C^{|l|^2} |q|^{|l|^2}.
\ed
$l_1 +2l_2 + \cdots + gl_g = -m$ implies that $|l| \geq \frac{|m|}{g^{1.5}}$.
Furthermore
\bd
\exists C, \delta > 0  : \forall |q| \leq \delta, 
n \in \natz _0, k_0 \in \natz: 
\sum_{|l| \geq k_0}|r(q)(n,l)| \leq 
C^{|k_0|^2} |q|^{|k_0|^2}.
\ed
Together with the observation, that the sum 
defining $s(q)$ does not contain the term where $l=0$, it follows,
that
\bd
\exists C, \delta > 0  : \forall |q| \leq \delta, 
n \in \natz _0, m \in \ganz: |s(q)(n,m)| \leq 
(C|q|)^{\max(1,\frac{m^2}{g^3})}.
\ed
Because of the estimate on the weightfunction $w(m) \leq Ce^{\sigma m}$,
for some constants $C, \sigma > 0$ (see Definition \ref{di2.2} and below),
we conclude that $s(q) \in {\cal L} _{0,w}$
for $|q|$ small enough and 
\bd
\exists C, \delta > 0 : \forall |q| \leq \delta: 
\| s(q) \| _{0,w} \leq C|q|.
\ed
Equation (\ref{g4.17}) shows, that $u_1(q) \in {\cal L} _{0,w}$ and
$\| u_1(q) \| _{0,w} \leq C|q|$ for $|q|$ small enough. Applying
the same kind of estimates again we obtain 
\bd
\exists C, \delta > 0 : \forall |q| \leq \delta: 
| u_2(q) | \leq C|q|.
\ed
This completes the proof of (I). \newline
\newline
{\bf (II):} \newline
The definitions at the beginning of the proof yield
\bd
\forall |m| > m_0: u(0,m) =
s(0,m) - s(1,m) +
\sum_{k=2}^{\infty} \frac{(-1)^{k-1}}{k} \left(
s(0,\cdot) ^{*k} (m) - s(1,\cdot) ^{*k} (m) \right).
\ed
Reworking the corresponding estimates in the proof of (I) one obtains
\bd
\exists C, \delta > 0 : \forall |q| \leq \delta, 
n \in \natz _0, |m| > m_0: |s(q)(n,m)| \leq 
(C|q|)^{\max(2,\frac{m^2}{g^3})},
\ed
from which we easily conclude that
\bd
\exists C, \delta > 0 : \forall |q| \leq \delta: 
\| s(q) {\bf 1}_{ \{|m|>m_0 \} } \| _{0,w} \leq C|q|^2.
\ed
On the other hand the estimate above on $ \| s(q) \| _{0,w}$ implies,
\bd
\exists C, \delta > 0 : \forall |q| \leq \delta, 
n \in \natz _0: 
\| \sum_{k=2}^{\infty} \frac{(-1)^{k-1}}{k} 
s(q)(n,\cdot) ^{*k} \| _{0,w} \leq C|q|^2.
\ed
This proves (II).\newline
\newline
{\bf (III):} \newline
Let us start with the differentiability of each $r(q)(n,l)$ with 
respect to $q_j$. From equation (\ref{g4.13}) and Theorem \ref{tg1.5} 
we learn that the only possible
problem lies at points where one of the $p_j = 0$, as
$p_j = \sqrt{q_{2j-1}^2 + q_{2j}^2}$ is not a differentiable function
of $q$ in those points. However, the functions $U$
and $\tau ^{(reg)}$ are even in each variable $p_j$ 
(see Theorem \ref{tg1.5}, $l_j (l_j - 1)$ is even and $\tilde{p} _j =
q_{2j-1} + i q_{2j}$. Hence $r(q)(n,l)$ is $C^1$.
Simple estimates show that
\bd
\exists C,\tilde{C}, \delta > 0 : \forall |q| \leq \delta, 
n \in \natz _0, l \in \ganz ^g \setminus \{ 0 \}: 
|\frac{\partial}{\partial q_j} r(q)(n,l)| \leq 
\tilde{C} C^{|l|^2} |q|^{|l|^2 - 1}(n + 1).
\ed
$
\exists C,\tilde{C}, \delta > 0 : \forall |q|, |q'| \leq \delta, 
n \in \natz _0, l \in \ganz ^g \setminus \{ 0 \} : 
$
\bd
|\frac{\partial}{\partial q_j} r(q')(n,l) - 
\frac{\partial}{\partial q_j} r(q)(n,l) | \\
\leq 
\tilde{C} C^{|l|^2} \max (|q'|,|q|)^{\max (|l|^2 - 2,0)}(n + 1)^2
\triangle (q',q),
\ed
with
\bd
\triangle (q',q) :=
\max \left(
|q' - q|, 
|\frac{\partial}{\partial q_j} U(q') - 
\frac{\partial}{\partial q_j} U(q) |,
|\frac{\partial}{\partial q_j} \tau^{(reg)}(q') - 
\frac{\partial}{\partial q_j} \tau^{(reg)}(q) | 
\right).
\ed
Note that the powers of $|l|$, which are produced by the 
differentiation have been subsumed into the $C^{|l|^2}$ term, simply by
increasing the constant $C$.

Our next goal is to prove that $q \mapsto s(q)$ is a $C^1$ map
into ${\cal L} _{- \sigma _2,w}$, with 
$\sigma _2 := \frac{\sigma _1}{2}$. The differentiability of
$r(q)(n,l)$ with respect to $q$ and the absolute and uniform
convergence of the sum
\bd
\sum_{
\scriptsize
\begin{array}{c} 
l \in \ganz ^g \setminus \{ 0 \} \\
l_1 + \cdots g l_g = -m
\end{array}
} 
\normalsize
\frac{\partial}{\partial q_j} r(q)(n,l),
\ed
shows that for all
$n \in \natz _0, m \in \ganz : s(q)(n,m)$ is a $C^1$ function of $q$. Therefore 
we have only to prove that the map $q \mapsto 
\frac{\partial}{\partial q_j} s(q)(n,m)$ is a continuous map into
${\cal L} _{- \sigma _2,w}$. This, however, follows by 
applying the arguments given in the proof of (I) to the estimates for 
$|\frac{\partial}{\partial q_j} r(q)(n,l)|$ and
$|\frac{\partial}{\partial q_j} r(q')(n,l) - 
\frac{\partial}{\partial q_j} r(q)(n,l) |$ given above and
from the observations that $\triangle (q',q) \rightarrow 0$, as 
$q' \rightarrow q$ and that 
$\sup_{n \geq 0} (n + 1)^2 e^{- \sigma _2 n} \leq \infty$.

Let us now investigate the differentiability of $u_1(q)$.
Employing the calculations we made in the proof of Theorem \ref{tb3.1},
step 4,
where we were in a similar situation, we obtain 
\bd
\frac{\partial}{\partial q_j} u_1(q) =
\sum_{k=0}^{\infty} (-1)^{k}  s(q)(n,\cdot) ^{*k} 
* \frac{\partial}{\partial q_j} s(q). 
\ed
Standard arguments yield the continuity of the map
$q \mapsto \frac{\partial}{\partial q_j} u_1(q)$ in
${\cal L} _{-2 \sigma _2,w} = {\cal L} _{- \sigma _1,w}$. 

The differentiability of the function $u_2(q)$ can be established
from the differentiability of $r(q)(n,l)$, the absolute and uniform
convergence of the corresponding sums 
(see equation (\ref{g4.18})) and the differentiability
of the logarithm away from $0$.
The proof of property (III) is completed. \newline \newline
{\bf (IV):} \newline
The differentiability was already proven in (III). We only 
have to investigate the terms of first order in $q$ of $u(q)(0,m)$ for
$1 \leq m \leq g$. Following the reasoning in the proof of (II),
we see that we only have to consider the first order terms of
$s(q)(0,m) - s(q)(1,m)$ which are given by
$r(q)(0,e_m) - r(q)(1,e_m)$, where $e_m := (0,\ldots ,-1, \ldots,0)$ 
denotes the unit vector in $\relz ^g$ with $-1$ at the $m$ -th entry.
This implies
\be
\label{g4.55}
\forall 1\leq m \leq g:
u(q)(0,m) =
\tilde{p}_m 2i \sin (\pi U_m (0)) 
\exp(\tau_{m,m}^{(reg)}(0)) + O(|q|^2).
\ee
We are done if we can show that 
$\sin (\pi U_m (0)) \neq 0$. It turns out that we can compute this
quantity explicitly. By Lemmas \ref{lg3.2}, \ref{lg3.4} 
and Proposition \ref{pg3.1} we conclude,
that
\be
U_m(0) =
- \frac{2}{\pi} \arctan \left(
\sqrt{\frac{b-\lambda _m}{\lambda _m -a}} \right). \label{g4.60}
\ee
Therefore
\lfd
\begin{eqnarray}
\sin (\pi U_m (0)) &=&
- 2 \tan \left(  \arctan \left(
\sqrt{\frac{b-\lambda _m}{\lambda _m -a}} \right) \right)
\cos ^2 \left(  \arctan \left(
\sqrt{\frac{b-\lambda _m}{\lambda _m -a}} \right) \right) \nonumber \\
&=&
- \frac{2}{b-a} \sqrt{(\lambda _m -a)(b - \lambda _m)} \nonumber \\
&=&
- \frac{2}{b-a} m \gamma \neq 0, \label{g4.65}
\end{eqnarray}
by Lemma \ref{lg3.5}.
\newline \newline
{\bf (V):}  See Appendix \ref{g2}, Theorem \ref{tg2.2}
\end{proof}
\lfd
\begin{remark}
\label{rg3.3}
The choice of $\lambda _j$ (compare with Remark \ref{rg3.2} (2), see also
(\ref{g3.3})).
\end{remark}
In Appendix \ref{g3} it is shown 
that the condition of time periodicity of the $g$-gap
solution leads to equation (\ref{g3.73}) in Lemma \ref{lg3.5}. 
This equation was 
solved by choosing the $\lambda _j$ 's as functions of the remaining
parameters. As it was stated in Remark \ref{rg3.2} (2), there exist two
choices for each $\lambda _j$, one in the lower half of the spectrum
and one in the upper half of the spectrum. We will now discuss
the difference of these two choices.
Proceeding as in Remark \ref{rl2.2}, we investigate
in which direction the energy is transported in the $g$-gap solution
up to second order in $q$.
The arguments given in the proof of property (IV) in the proof
of Theorem \ref{tg4.1} above yield
\bea
x_0^{(0)} (t) &=&
\sum_{m=1}^g u(q)(0,m) e^{i \gamma m t} + O(|q|^2), \\
x_1^{(0)} (t) &=& c +
\sum_{m=1}^g u(q)(0,m) \exp (-2 \pi i U _m (0) )
e^{i \gamma m t} + O(|q|^2),
\eea
where $u(q)(0,m)$ was determined in (\ref{g4.55}).
Now we compute the energy which is exchanged between particles $x_0^{(0)}$
and $x_1^{(0)}$ during one period. Repeating the calculations in 
Remark \ref{rl2.2} we arrive at
\bea
E &=& - \int_0^{\frac{2 \pi}{\gamma}} 
F(x_0^{(0)} (t) - x_1^{(0)} (t) ) \dot{x}_0^{(0)} (t) dt \\
&=& 4 \pi \exp (-c) \sum_{m=1}^g m |u(q)(0,m)|^2 \sin (- 2 \pi U _m (0) )
+ O(|q|^3).
\eea
Each term in the sum corresponds to the energy transported by one 
phase of the multiphase solution. It is transported in the direction
of  increasing $n$, if 
\bd
\sin (- 2 \pi U _m (0) ) < 0.
\ed
Equation (\ref{g4.60}) implies that
\bd
- 2 \pi U _m (0) = 4  \arctan \left(
\sqrt{\frac{b-\lambda _m}{\lambda _m -a}} \right).
\ed
Therefore $- 2 \pi U _m (0) \in (\pi, 2 \pi)$, if $\lambda _m < 
\frac{a + b}{2}$ and
$- 2 \pi U _m (0) \in (0, \pi)$, if $\lambda _m > \frac{a + b}{2}$.
Thus the choice we make in (\ref{g3.3}), 
corresponds to a solution where the energy is transported outwards
in the direction of increasing $n$. As described above our choice
implies that all gaps open up only in the lower half of the spectrum.

\lfd
\begin{remark}
\label{rg4.5}
Comparison with solutions of general lattices for $g = 0, 1$.
\end{remark}
Recall, that for $m_0 =0$ and $m_0 =1$, we were able to construct
periodic solutions for general lattices (see Sections \ref{m4} and
\ref{l3} respectively). We will verify that these solutions 
are the same as we have constructed in the present section
in the Toda case, by
showing that they produce the same families of sequences 
of Fourier coefficients $u(q)$, which are used in the construction
described in Chapter \ref{m}.

The case $g=0$ is trivial, as the $0$-gap solution is simply
given by $x_n (t) = c n$, i.e. $u(q) = 0$, which is the choice
we made in Section \ref{m4}.

In the case $g = 1$ equation (\ref{g3.110}) and the periodicity
of the theta function in the real direction shows that the 
one-gap solution is of a form to which Corollary
\ref{cl2.1} applies. We have to check the assumptions (I)-(IV) in the 
Corollary. (I) and (III) are trivially satisfied.

Property (IV) is satisfied, if we can prove, that $-2 \pi U(0)
 = \beta _1$, where $\beta _1$ was defined in (\ref{i1.25}). By the Remark 
\ref{rg3.3} above
we know already that both quantities lie in $(\pi, 2 \pi)$ 
(mod $2 \pi$) and therefore it
suffices to show that $\cos( - 2 \pi U(0)) = \cos ( \beta _1)$.
Equation (\ref{i1.25}) yields
$\cos ( \beta _1) = - \frac{\delta _1}{2}$.
Using (\ref{g4.65}) we compute
\bea
\cos( - 2 \pi U(0)) + \frac{\delta _1}{2}
&=&
\left( 1 - \frac{8 \gamma ^2}{(b-a)^2} \right) +
\left( -1 + \frac{ \gamma ^2 }{2 \exp (-c)} \right) \\
&=&
0,
\eea
as by the choice of $b^{(0)}$ 
(see (\ref{g3.98})) in the proof of Lemma \ref{lg3.7},
we have $b-a = 4 \exp (-\frac{c}{2} ) $.

Property (II) of Corollary \ref{cl2.1}
is not satisfied as the zero Fourier coefficient 
\bd
s(0) = 
\ln
\frac{\vartheta (\frac{1}{2} U - Z |\tau)}
     {\vartheta (- \frac{1}{2} U - Z |\tau)}
\ed
does not necessarily vanish.
But this corresponds only to adding a constant to $u(q)(n,0)$,
which does not change the equation for $v$ in Lemma \ref{lb2.1}
(note that $W(u,v)$ only contains $\triangle u$, see Definition
\ref{db2.2}) and which does not change the value of $a(n,m)$ in 
Lemma \ref{lb2.1}.

\lfd
\begin{remark}
\label{rg5.2}
The opening of the gaps to first order in $\epsilon$.
\end{remark}
We will now investigate the size of the instability regions and 
their dependence on $\epsilon$ and the Fourier coefficients of the
driver. This question of basic interest was well studied in the
continuous case for periodic potentials (see \cite{MW} and references
therein). In our situation
recall that $L$ denotes the semi-infinite Lax operator corresponding to the
solution constructed in Theorem \ref{tg4.1} above (see also \ref{e.135}).
As $L$ is obtained from the doubly infinite operator $\tilde{L}$ 
by restricting it to $n \geq 1$ and adding an operator decaying exponentially
in $n$ (which corresponds to the $v$-term in \ref{g4.3}) standard
arguments of spectral theory imply that
\bd
\sigma _{\mbox{ess}} (L) = \sigma _{\mbox{ess}} (\tilde{L}),
\ed
and hence
\bd
\sigma _{\mbox{ess}} (L(t)) =
\left[ \frac{a}{2}, \frac{\lambda_1 - p_1}{2} \right] \cup 
\left[ \frac{\lambda _1 + p_1}{2}, \frac{\lambda _2 - p_2}{2} \right] \cup
\cdots \cup \left[ \frac{\lambda _g + p_g}{2}, \frac{b}{2} \right].
\ed
(See equation (\ref{g2.5})).
Therefore the width of the $j$-gap in the spectrum of $L(t)$ 
is given by $|p_j|, 1 \leq j \leq g$. They were
determined as functions of $\epsilon$, when we solved the 
resonance equations $\epsilon b_m = u(q)(0,m) +v(q, \epsilon)(0,m),
1 \leq m \leq g$ in Theorem \ref{tg4.1}.
As $v(q,\epsilon)$ is a fixed-point of the operator $T(q, \epsilon, \cdot)$
given by  (\ref{b3.16}) in Chapter \ref{m}, we conclude from 
Theorem \ref{tb3.1} (ii), Proposition \ref{pb2.3} (ii) and 
the fact that $q(\epsilon) = O(\epsilon)$ that
$v(q, \epsilon)(0,m) = O(\epsilon ^2),$ for $1 \leq m \leq g$. 
Equation (\ref{g4.55}) then yields that
the following relation holds.
\be
\label{g5.80}
|\epsilon b_m| = 2 |\sin (\pi U_m (0)) |
\exp(\tau_{m,m}^{(reg)}(0)) |p_m| + O(\epsilon ^2).
\ee
Below we will show that 
\be
\label{g5.85}
\tau_{m,m}^{(reg)}(0) = - \ln \frac{8 (m \gamma)^2}{b-a}.
\ee
Using (\ref{g4.65}) and (\ref{g5.80}) we can express $|p_m|$.
\be
\label{g5.90}
|p_m| = 2 |\epsilon b_m| m \gamma + O(\epsilon ^2).
\ee
We compare this formula with numerical results shown in 
Figures \ref{fe5} and \ref{fe6}. We recall the choice of parameters,
$|\epsilon b_1| = 0.1, |\epsilon b_2| = 0.05$. Up to first order,
equation (\ref{g5.90}) yields
\begin{itemize}
\item
$\gamma = 1.8 : p_1 \approx 0.36 $.
\item
$\gamma = 1.1 : p_1 \approx 0.22, p_2 \approx 0.22 $.
\end{itemize}
These values are in good agreement with the numerical experiments. In fact,
\begin{itemize}
\item
$\gamma = 1.8, \mbox{ Figure \ref{fe5} }: p_1 = 0.34 \pm 0.007 $.
\item
$\gamma = 1.1, \mbox{ Figure \ref{fe6} }: p_1 = 0.223 \pm 0.005, 
p_2 = 0.215 \pm 0.15 $.
\end{itemize}
We will now sketch the derivation of formula (\ref{g5.85}).
By Theorem \ref{tg1.5} and Proposition \ref{pg3.1} (vi), we may consider
the one-gap situation, where only the $m$ -th gap is open and all other
gaps are closed. Using equations (\ref{g3.60}),(\ref{g3.61}),(\ref{g3.63})
(\ref{g3.31}),(\ref{g3.20}) and (\ref{g3.6}) and Lemmas 
\ref{lg3.1}, \ref{lg3.2} we arrive at
\be
\label{g5.100}
 \tau_{m,m}^{(reg)}(0) =
- \frac{1}{h_m}  \lim_{p_m \rightarrow 0} \left(
\int_a^{\lambda _m - p_m}
\frac{1}
{\sqrt{(E-a)(b-E)[(E-\lambda _m)^2 - p_m^2]}} dE
- h_m \ln \frac{1}{p_m} 
\right)\hspace{-4pt}.
\ee
The quantity $h_m$ was defined in (\ref{g3.15}) and is given
in this case by
\be
\label{g5.105}
h_m = \frac{1}{\sqrt{(\lambda _m - a)(b - \lambda _m)}}.
\ee
We define the auxiliary functions
\lfd
\begin{eqnarray}
f(E) &:=& \frac{1}{\sqrt{(E - a)(b - E)}} \label{g5.110}, \\
g(E) &:=& \frac{f(E) - f(\lambda _m)}{E - \lambda _m}.\lfd \label{g5.111}
\end{eqnarray}
Then
\lfd
\begin{eqnarray}
&&
\lim_{p_m \rightarrow 0} \left(
\int_a^{\lambda _m - p_m}
\frac{1}
{\sqrt{(E-a)(b-E)[(E-\lambda _m)^2 - p_m^2]}} dE
- h_m \ln \frac{1}{p_m} \right)  \nonumber \\
&=&
\lim_{p_m \rightarrow 0} \left(
\int_a^{\lambda _m - p_m}
\frac{f(\lambda _m)}
{\sqrt{(E-\lambda _m)^2 - p_m^2}} dE
- h_m \ln \frac{1}{p_m} \right) \nonumber \\
&& +
\lim_{p_m \rightarrow 0} \left(
\int_a^{\lambda _m - p_m}
\frac{(E - \lambda _m) g(E)}
{\sqrt{(E-\lambda _m)^2 - p_m^2}} dE \right) \nonumber \\
&=&
(I) + (II). \label{g5.112}
\end{eqnarray}
Using the appropriate changes of variables we evaluate
\lfd
\begin{eqnarray}
(I) &=& h_m \ln 2(\lambda _m - a). \label{g5.115} \\
(II) &=& \int_a^{\lambda _m} g(E) dE \nonumber \\
&=& \lim_{\epsilon \rightarrow 0} \left(
\int_a^{\lambda _m - \epsilon} \frac{f(E)}{E - \lambda _m} dE -
h_m \ln \frac{\lambda _m - a}{\epsilon} \right). \lfd \label{g5.120}
\end{eqnarray}
The remaining integral can be integrated and we obtain
\be
\label{g5.125}
\int_a^{\lambda _m - \epsilon} \frac{f(E)}{E - \lambda _m} dE =
\frac{1}{m \gamma} \left(
\ln \frac{1}{\epsilon} + \ln \frac{4 (m \gamma)^2}{b-a} + \ln 
(1 + O(\epsilon)) \right).
\ee
Note that for $\lambda _m = \lambda _m (p_m = 0)$ equation (\ref{g3.74})
yields $h_m = \frac{1}{m \gamma}$.
Therefore we can determine 
\lfd
\begin{eqnarray}
(II) &=& h_m \left(\ln \frac{4 (m \gamma)^2}{b-a} -
\ln (\lambda _m - a) \right).
\label{g5.130} \\
(I) + (II) &=& h_m \ln \frac{8 (m \gamma)^2}{b-a}.
\lfd \label{g5.135}
\end{eqnarray}
It follows from (\ref{g5.100}) and (\ref{g5.112}) that
$\tau_{m,m}^{(reg)}(0) = -\ln \frac{8 (m \gamma)^2}{b-a}$,
which proves (\ref{g5.85}).

\appendix
\chapter{Definition of the g-gap solution}
\label{g2}

It is our goal to derive a formula for the well known $g$-gap
solutions of the Toda lattice.

Let $g \in \natz _0 = \{ n \in \ganz: n \geq 0 \}$ and denote by 
$R_g$ the hyperelliptic curve of
genus g, which is constructed by pasting together two copies of the
Riemann sphere $\comz \cup \{ \infty \}$ along the slits 
$[E_0,E_1], [E_2,E_3], \ldots, [E_{2g},E_{2g+1}]$ ,where 
$(E_0 < E_1 < \cdots < E_{2g+1} ).$

A point on $R_g$ is denoted by $P$, and the canonical projection of $R_g$ on
the Riemann sphere is given by $\pi$. We write $E = \pi (P).$
\begin{itemize}
\item
For $ 1 \leq k \leq g$ let $P_j \in R_g$ be a point in the j-th gap, i.e
$ \pi (P_j) \in [E_{2j-1},E_{2j}].$
\item
$\alpha _k, \beta _k ,
1 \leq k \leq g$ denote the canonical homology basis for $R_g$ 
(see \cite[III.1]{FK} for a definition of a canonical homology basis)
as shown in the Figure \ref{fa1} below.
\end{itemize}
\bef
\leavevmode \epsfysize=2.5cm
\epsfbox{fa1.ps}
\caption{The canonical homology basis of the Riemann surface}
\label{fa1}
\eef
The existence and uniqueness of the following differentials can be 
deduced from the Riemann Roch Theorem (see \cite[III.4.8]{FK})
applied to $R_g$
(see e.g. \cite[VI.4, Satz 32]{BS}, \cite[Remark 16.17]{Ges}).
\begin{itemize}
\item
$\omega := (\omega _l )_{1 \leq l \leq g}$ is defined as the unique basis 
of holomorphic differentials on $R_g$, which is normalized by the 
condition $\int_{\alpha _k} \omega _l = \delta _{k,l}$. Note that all 
the holomorphic differentials on $R_g$ can be written in the following form.
\lfd
\beqa
\nu &=& \frac{p(E)}{\sqrt{R(E)}} dE, \label{aa.3} \\
R(E) &:=&  \prod_{j=0}^{2g+1} (E-E_j) \lfd \label{g3.4},
\eeqa 
where $p$ is any polynomial in $E$ of degree $\leq g-1$
(see \cite[III.7.5, Corollary 1]{FK}).
\item
$\omega ^{(1)}$ is the unique meromorphic differential on $R_g$ with
simple poles at $P_{\infty}$ and $P_{\infty}^* $(residues: $1$ /resp. $-1$),
holomorphic everywhere else and normalized by 
$\int_{\alpha _k} \omega ^{(1)} = 0$, for $ 1 \leq k \leq g$.
Here the points $P_{\infty}$ and $P_{\infty}^* $ denote the point at 
infinity on the upper and lower sheet respectively. We can write
\be
\omega ^{(1)} = \frac{p^{(1)}(E)}{\sqrt{R(E)}} dE,
\label{aa.5}
\ee
where $p^{(1)}$ denotes an uniquely determined polynomial of degree
$g$.
\item
$\omega ^{(2)}$ is the unique meromorphic differential on $R_g$ with
second order poles at $P_{\infty}$ and $P_{\infty}^* $
(principal parts: $\frac{1}{2\xi ^2}d\xi$ /resp.
$\frac{-1}{2\xi ^2}d\xi$, where $\xi$ is the coordinate around 
$\infty : \xi := \frac{1}{\pi (P)}),$
holomorphic everywhere else and normalized by 
$\int_{\alpha _k} \omega ^{(2)} = 0$, for $ 1 \leq k \leq g$.
Again this differential can be written in the following way.
\be
\omega ^{(2)} = \frac{p^{(2)}(E)}{\sqrt{R(E)}} dE,
\label{aa.10}
\ee
where $p^{(2)}$ denotes an uniquely determined polynomial of degree
$g+1$.
\item
We define the following $\beta$ periods.\newline
\be
\tau := (\tau _{k,l})_{1 \leq k,l \leq g},\quad \mbox{ with }
\tau _{k,l} := \int_{\beta _k} \omega _l.  
\label{aa.15}
\ee
The Riemann bilinear relations imply that $\tau$ is a symmetric matrix
with positive definite imaginary part (see \cite[III.3.1, III.3.2]{FK}).
\lfd
\beqa
\label{g2n.1}
U_k &:=& \frac{1}{2\pi i} \int_{\beta _k} \omega ^{(1)}. \\
\lfd \label{g2n.2}
V_k &:= & \frac{1}{2\pi i} \int_{\beta _k} \omega ^{(2)}.
\eeqa
\item
We denote K to be the vector of Riemann constants with respect to
the basepoint $E_{2g + 1}$ (compare \cite[VI.2.4]{FK}). 
\end{itemize}

The following three integral functions are multivalued and 
depend on the path of integration.\newline 
$A(P) :=  \int_{E_{2g+1}}^P \omega.$ \newline
$\Omega ^{(1)}(P) :=  \int_{E_{2g+1}}^P \omega ^{(1)}$ on $R_g \setminus
\{ P_{\infty} , P_{\infty}^* \} $ .\newline
$\Omega ^{(2)}(P) :=  \int_{E_{2g+1}}^P \omega ^{(2)}$ on $R_g \setminus
\{ P_{\infty} , P_{\infty}^* \} $.

We recall the definition of the Riemann theta 
function (see e.g. \cite[VI.1]{FK}).
For a symmetric matrix $\chi$, with positive definite imaginary part,
\bd
\vartheta (v| \chi) := \sum_{m \in \ganzz^g}
e^{2 \pi i <m,v> + \pi i <m,\chi m>},
\ed
where $<u,v> := \sum_{j = 1}^g u_j v_j.$ 
The ``periodicity'' properties of the Riemann 
theta function are as follows. Let $1 \leq m \leq g$
and denote by $e_m$ the $m$ -th column of the $g \times g$ identity matrix
and by $\chi _m$ the $m$-th column of $\chi$. Then
(see e.g. \cite[VI.1.2]{FK})
\lfd
\begin{eqnarray}
\vartheta (v + e_m | \chi) &=&  \vartheta (v | \chi). \label{g2n.10} \\
\vartheta (v + \chi _m | \chi) &=&
\exp (- 2 \pi i v_m-\pi i \chi _{m,m}) \vartheta (v | \chi).
\lfd \label{g2n.15}
\end{eqnarray}

The main tool for the construction of the $g$-gap solutions
is the following existence and uniqueness theorem for the {\it Baker -
Akhiezer} function. Using the above notation, we have:

\lfd
\begin{theorem} 
\label{tg2.1}
(\cite{Kri1}, see also \cite{Akh},\cite{Dub}, \cite{IM})
For all $n \in \ganz, t \in \relz$, there is a unique (up to
multiplication by a constant ) function  $\psi _n(t,\cdot)$,
which is not identically equal to $0$ and satisfies
\begin{itemize}
\item[(i)]
$\psi _n(t,\cdot)$ is meromorphic on
$R_g \setminus \{ P_{\infty} , P_{\infty}^* \}. $
\item[(ii)]
$\psi _n(t,\cdot)$ has at most simple poles at $P_1, \ldots ,P_g  \in
R_g \setminus \{ P_{\infty} , P_{\infty}^* \} $ and is holomorphic 
elsewhere on $R_g \setminus \{ P_{\infty} , P_{\infty}^* \}$.
\item[(iii)]
$\psi _n(t,\cdot) E^{\pm n} exp(\pm \frac{t}{2} E)$ has a holomorphic 
continuation at $P_{\infty}$ /resp. $P_{\infty}^*.$
\end{itemize}
\end{theorem}

\begin{proof}
This result is well known and we only sketch the proof. First we prove
existence. Define
\be
\label{g2.0}
\Psi _n (t,P) :=
exp(n \Omega ^{(1)}(P) + t \Omega ^{(2)}(P))
\frac{\vartheta (A(P) + Un + Vt - Z | \tau)}
     {\vartheta (A(P) - Z | \tau)},
\ee
where $\Omega ^{(1)}, \Omega ^{(2)}, A, U, V$ and $\tau$ were defined above
and $Z$ will be introduced below.
Before we can check that this function satisfies all the conditions of 
the theorem, some more remarks are needed.
\begin{itemize}
\item
As already mentioned above the functions $A, \Omega ^{(1)}, \Omega ^{(2)}$
are multivalued. We will show below that nevertheless the function
$\Psi _n (t,P)$ is well defined if we insist that the path of integration 
is the same for all three functions.
\item
$Z := K + \sum_{1 \leq m \leq g} A(P_m)$. As we want $\Psi _n (t,P)$
to be well defined we now specify the path of integration from $E_{2g+1}$ to
$P_m$. We first integrate on the upper sheet along $\relz$ on the $\comz 
^+ := \{ E \in \comz: \mbox{ Im}(E) > 0 \}$ side 
from $E_{2g + 1}$ to the branchpoint $E_{2m-1}$ and then from 
$E_{2m-1}$ to $P_m$ along the real axis, where the path always first stays
on the upper sheet, and in case that $P_m$ lies on the lower sheet,
we switch the sheet at $E_{2m}$. Of course we could have chosen
any fixed path from $E_{2g+1}$ to $P_m$, but the preceeding choice leads
to an especially simple formula for $Z$ in (\ref{g3.82}) below.
\item
In the case that $g=0$ the theta functions are simply replaced by the factor 1.
One checks that all the steps of the proof given below are trivially
satisfied in this case.
\end{itemize}

In order to see that $\Psi _n (t,P)$ is well defined we have to examine
what happens if we add to the path of integration one of the cycles
$\alpha _k, \beta _k, 1 \leq k \leq g$ or $\gamma _j , j=1,2$, where the 
$\gamma$ 's are the cycles around $P_{\infty}$ /resp. $P_{\infty}^*$.
The $\gamma$ cycles only effect $\Omega ^{(1)}$ by adding a multiple
of $2 \pi i$, which does not change the value of 
$\Psi _n (t,P)$ because of the exponentiation.
The $\alpha$ cycles only change the entries of $A(P)$ by adding integers
which has no effect because
of the periodicity of the theta function in the real direction (see
(\ref{g2n.10})).
The $\beta$ cycles change all three integrals and it is straightforward to
show that all the factors cancel out by the monodromy property of the
theta function (\ref{g2n.15}) and by the choice of the vectors $U$ and $V$
(see (\ref{g2n.1}) and (\ref{g2n.2})).

As $A, \Omega ^{(1)}, \Omega ^{(2)}$ are holomorphic in
$R_g \setminus \{ P_{\infty} , P_{\infty}^* \}$, properties (i) and (ii)
are equivalent to showing
that the zeros of $\vartheta (A(P) - Z | \tau)$ are all simple and that they
are given by $P_1, P_2, \ldots , P_g$. By definition of $Z$ and by
\cite[Theorem b,VI.3.3]{FK},
(see also \cite[Thm 17.9]{Ges}) this is equivalent to proving that the 
divisor $P_1 P_2 \ldots P_g$ is nonspecial, i.e that there exists no
holomorphic differential on $R_g$ vanishing at $P_1, P_2, \ldots , P_g$
and which is not identically equal to zero.
But this is easily deduced from the characterisation of holomorphic
differentials as given in (\ref{aa.3}) and from the fact that
$\pi(P_i) \neq \pi(P_j)$ for $i \neq j$.

Property (iii) can be verified from the definitions of
$\omega ^{(1)}$ and $\omega ^{(2)}.$

This settles existence and we can now turn to the question of uniqueness.
Suppose $\tilde{\Psi} _n (t,\cdot)$ satisfies (i),(ii) and (iii). 
We consider the function 
$\frac{\tilde{\Psi} _n (t,\cdot)}{\Psi _n (t,\cdot)}$. It is meromorphic
and its poles are the zeros of
$\vartheta (A(P) + Un + Vt - Z | \tau)$. It will be shown below 
(see Remark \ref{rg3.5}) that for all
$n$ and $t$, $\vartheta (A(P) + Un + Vt - Z | \tau)$ has exactly
g zeros, one in each gap, which form a nonspecial divisor 
and by the Riemann-Roch theorem it follows that
$\frac{\tilde{\Psi} _n (t,\cdot)}{\Psi _n (t,\cdot)}$ must be 
a constant (see e.g. \cite[III.4.8]{FK}, \cite[Thm 16.11]{Ges}).
\end{proof}

We now introduce a normalization of the {\em Baker -Akhiezer} function.
Denote by $\psi _n (t,\cdot)$ the uniquely defined 
BA-function which has the following expansion at $P _{\infty}$.
\be
\label{g2.1}
\psi _n (t,P) =
E^{-n} exp(- \frac{t}{2} E) 
(1 + \sum_{s=1}^{\infty} \xi _s^+ (n,t) E^{-s}).
\ee
This is possible as we see from (\ref{g2.0}) and the position of the 
zeros of $\vartheta (A(\cdot) + Un + Vt - Z | \tau)$ described above,
that $\Psi _n (t,P_{\infty}) \neq 0$.
The expansion at $P _{\infty}^*$ is then written as
\be
\label{g2.2}
\psi _n (t,P) =
E^{n} exp(\frac{t}{2} E) 
(\xi _0^- (n,t) + \sum_{s=1}^{\infty} \xi _s^- (n,t) E^{-s}).
\ee
We can express $\xi _0^- (n,t)$ in the following way.
Let us expand $\Psi _n (t,P)$ as given in equation (\ref{g2.0})
around $P_{\infty}$ and $P_{\infty}^*$.
From the definition of $\omega ^{(1)}$ and $\omega ^{(2)}$ it is obvious
that we can expand their integrals in the following way around  
the infinities.
\be
\label{g2.9}
\Omega ^{(1)} (P) = 
\mp ( \ln E + \sum_{l=0}^{\infty} I_l^{\pm} E^{-l}) 
\quad \mbox{ around } P_{\infty} \mbox{ /resp. } P_{\infty}^*.
\ee
\be
\label{g2.10}
\Omega ^{(2)} (P) = 
\mp \frac{1}{2}( E + \sum_{l=0}^{\infty} R_l^{\pm} E^{-l}) 
\quad \mbox{ around } P_{\infty} \mbox{ /resp. } P_{\infty}^*.
\ee
In order to have the zero order terms $I_0^{\pm}$ and $R_0^{\pm}$
well defined we shall now fix the path of integration from $E_{2g+1}$ 
to $P_{\infty}$ and $P_{\infty}^*$ as the path along the real axis
from $E_{2g+1}$ to $ + \infty$ on the upper sheet /resp. on the lower
sheet.
The holomorphic continuation of
$ \Psi _n (t,P) E^n exp(\frac{tE}{2})$ takes at $P_{\infty}$ the value
\bd
exp(-n I_0^+ - \frac{t}{2} R_0^+)
\frac{\vartheta (A(P_{\infty}) + Un + Vt - Z | \tau)}
     {\vartheta (A(P_{\infty}) - Z | \tau)}.
\ed
The path in the integration of $\omega$ is chosen as above, which determines 
the value of $A(P_{\infty})$ uniquely. It follows
by the Riemann bilinear relations (see e.g. \cite[III(3.6.3)]{FK}),
that
\be
\frac{1}{2\pi i} \int_{\beta} \omega ^{(1)}
= \int _{P_{\infty}^*}^{P_{\infty}} \omega = A(P_{\infty}) - 
A(P_{\infty}^*).
\label{aa.30}
\ee
As $\omega$ differs on the different sheets only by a sign we obtain
with (\ref{g2n.1}), that
\be
\label{g2.53}
A(P_{\infty}) = \frac{1}{2} U = -A(P_{\infty}^*).
\ee
At $P_{\infty}^*$ the holomorphic continuation of
$ \Psi _n (t,P) E^{-n} exp(-\frac{tE}{2})$ takes the value
\bd
exp(n I_0^- + \frac{t}{2} R_0^-)
\frac{\vartheta (A(P_{\infty}^*) + Un + Vt - Z | \tau)}
     {\vartheta (A(P_{\infty}^*) - Z | \tau)}.
\ed
Hence
\lfd
\begin{eqnarray}
\xi _0^- (n,t) &=& exp \left(n (I_0^+ + I_0^-) 
+ \frac{t}{2} (R_0^+ + R_0^-) \right) 
\nonumber \\
&&
\frac{\vartheta ( U(n-\frac{1}{2}) + Vt - Z | \tau)
      \vartheta ( \frac{1}{2} U - Z | \tau)}
     {\vartheta ( U(n+\frac{1}{2}) + Vt - Z | \tau)
      \vartheta (- \frac{1}{2} U - Z | \tau)}. \label{g2.11} 
\end{eqnarray}
In Appendix \ref{g3} we will determine formulae for 
$I_0^+, I_0^-, R_0^+, R_0^-,
U, V, Z$, from which it is easy to see that all the parameters are real 
and hence $\xi _0^- (n,t) \in \relz \setminus \{0\}$. 
But $\xi _0^- (0,0) = 1$ and
by the hence by continuity (regard $n$ as an arbitrary real variable) 
$\xi _0^- (n,t) > 0$. We define
\be
\label{g2.3}
x_n (t) := \ln \xi _0^- (n,t)
\ee
as a real number. The next theorem shows that $x_n (t)$ is a solution
of the Toda lattice.

\lfd
\begin{theorem} 
\label{tg2.2}
(\cite{Kri1}, \cite{Ges2})
Let $x_n (t)$ be defined as in (\ref{g2.3}). Then $x_n$ is
twice differentiable (in fact infinite differentiable) 
for all $n \in \ganz$ and furthermore 
\bea
\forall n \in \ganz, t \in \relz &:&
\ddot{x} _n (t) = exp(x_{n-1}(t) - x_n(t)) - exp(x_n(t) - x_{n+1}(t)).\\
\eea
\end{theorem}

\begin{proof}
Recall from the Introduction that we can write the Toda 
equation in Lax pair form by using Flaschka variables 
(see (\ref{e.60}) -- (\ref{e.75})).
We will show below that for 
\be
\label{g2.50}
\tilde{\psi} _n(t,P) :=  
exp(- \frac{x_n}{2}) \psi _n(t,P), 
\ee
the following equations are 
satisfied for all $t \in \relz$ and $P \in  
R_g \setminus \{ P_{\infty} , P_{\infty}^* \}. $
\be
\tilde{L} \tilde{\psi} = \frac{1}{2} E \tilde{\psi}.
\label{aa.35}
\ee
\be
\frac{\partial}{\partial t} \tilde{\psi} = \tilde{B} \tilde{\psi}.
\label{aa.40}
\ee
Equation (\ref{e.75}) follows as the compatibility condition of
equations (\ref{aa.35}) and (\ref{aa.40}).\newline 
{\bf Proof of equations (\ref{aa.35}) and (\ref{aa.40}).}
\newline
We define the auxiliary function $\triangle _1:= ((\tilde{L} - \frac{E}{2})
\tilde{\psi})_n$.
Using (\ref{g2.1}) we can calculate the behavior of 
$\triangle _1 E^{n} exp( \frac{t}{2}E)$ around $P_{\infty}$ and we arrive at
\bd
\frac{1}{2} exp(-\frac{x_n}{2}) ( \xi _1^+ (n-1,t) - 
\xi _1^+ (n,t) - \dot{x} _n)  + O(E^{-1}).
\ed
In the same way we determine from equation (\ref{g2.2}) the behavior of
$\triangle _1 E^{-n} exp(- \frac{t}{2}E)$ around $P_{\infty}^*.$ It is
given by
\bd
\frac{1}{2} exp(\frac{x_n}{2}) ( e^{-x_{n+1}} \xi _1^- (n+1,t) - 
e^{-x_{n}} \xi _1^- (n,t) - \dot{x} _n)  + O(E^{-1}).
\ed
The quantity $\triangle _1$ satisfies conditions (i) and (ii) 
of Theorem \ref{tg2.1} and
thus we can conclude that
\be
\label{g2.7}
\hspace{.45in}
\xi _1^+ (n-1,t) - \xi _1^+ (n,t) - \dot{x} _n  =
e^{-x_{n+1}} \xi _1^- (n+1,t) - e^{-x_{n}} \xi _1^- (n,t) - \dot{x} _n.
\ee

Proceeding in an analog way for the auxiliary function
$\triangle _2= ((\frac{\partial}{\partial t} -\tilde{B}) 
\tilde{\psi})_n$ we conclude
again from Theorem \ref{tg2.1} that
\be
\label{g2.8}
\hspace{.5in}
\xi _1^+ (n-1,t) - \xi _1^+ (n,t) - \dot{x} _n  =
- e^{-x_{n+1}} \xi _1^- (n+1,t) + e^{-x_{n}} \xi _1^- (n,t) + \dot{x} _n.
\ee

Adding equations (\ref{g2.7}) and (\ref{g2.8}) we see that
$\xi _1^+ (n-1,t) - \xi _1^+ (n,t) - \dot{x} _n  = 0$ and again by 
Theorem \ref{tg2.1} it follows that $\triangle _1$ and $\triangle _2$ are
both equal to zero and this proves equations (\ref{aa.35}) and (\ref{aa.40}).
\end{proof}

\chapter{Evaluation of the g-gap solution for small gaps}
\label{g3}

The goal of this appendix is to prove Theorem \ref{tg1.5} , that is to produce 
a 2g-parameter  family of time periodic solutions
of the doubly infinite Toda lattice with period $\frac{2 \pi}{\gamma}$ 
satisfying certain regularity conditions. 
We use the g -gap solutions produced in the last section in the case 
that the gaps are small. 
First we introduce some notation.
We choose the variables $\lambda _j, p_j, 1 \leq j \leq g$ in such a way,
that $E_{2j-1} = \lambda _j - p_j$ and $E_{2j} = \lambda _j + p_j , 
1 \leq j \leq g$. In addition we call $a:= E_0$ and $b:=E_{2g + 1}$.
Hence the g-gap solutions are completely parameterized by the $3g + 2$
real quantities $a,b, \lambda _j, p_j, P_j , 1 \leq j \leq g$.
The plan of the present Appendix is as follows. We will determine 
the dependence of the quantities 
$\tau, U,
V, Z, I_0^{\pm}, R_0^{\pm}$ in equation (\ref{g2.11})
on the parameters  $a,b, \lambda _j, p_j, P_j , 1 \leq j \leq g$.
Then we will proceed to show that for small gaps there is indeed
a choice of parameters, such that the corresponding solution $x_n(t)$
as given by (\ref{g2.3}) has all the properties of Theorem
\ref{tg1.5}. Special emphasis will be given to the limits $p_j 
\rightarrow 0$, i.e. when gaps close. Even though the analytic 
expressions will achieve some limit value it is not clear a priori, that
this limit coincides with the formula for the corresponding 
lower-gap solution. Therefore this has to be checked separately.
\lfd
\begin{remark}
\label{rg3.1}
\begin{em}
During the Appendix we will always assume that the hypothesis
of Theorem \ref{tg1.5} holds, i.e.
\bd
g \gamma < 2e^{-\frac{c}{2}} < (g+1) \gamma.
\ed
We specify the range for the different parameters.
\lfd
\begin{eqnarray}
a &\in& (-2e^{-\frac{c}{2}} - \epsilon _a,  
-2e^{-\frac{c}{2}}+ \epsilon _a). \label{g3.3a} \\
b &\in& (2e^{-\frac{c}{2}} - \epsilon _b,  
2e^{-\frac{c}{2}}+ \epsilon _b). \lfd \label{g3.3b}  \\
p_j &\in& (0, \epsilon _{p_j}) . \lfd \label{g3.3c} \\
\lambda _j &\in& (\lambda _j^{(0)} - \epsilon _{\lambda _j} ,
                  \lambda _j^{(0)} + \epsilon _{\lambda _j} ), \mbox{ where }
\lfd \label{g3.3d}  \\
\lambda _j^{(0)} &:=& \frac{a+b}{2} - 
\sqrt{ \left( \frac{b-a}{2} \right)^2 - j^2 \gamma ^2 }, 1 \leq j \leq g.  
\lfd \label{g3.3} \\
\pi ( P_j ) &\in& [\lambda _j - p_j, \lambda _j + p_j ]. \lfd \label{g3.3f}
\end{eqnarray}
The main requirements which have to be satisfied for the choice 
of the various $\epsilon$ 's is that we have to ensure that the bands
do not vanish, i.e. 
$ a < \lambda _1 - p_1 < \cdots < \lambda _g + p_g < b $. Furthermore
we want $\lambda _j^{(0)}$ to be real, which can be achieved by
the assumption (\ref{g1.30}) above together with the choice of $a$ and $b$.
 During the
calculations some other conditions on the smallness of the $\epsilon$ 's
will occur (e.g. induced by the use of the implicit function theorem)
and of course we want them to be satisfied as well.
\end{em}
\end{remark}

\section{The holomorphic differentials and the $\tau$ matrix} 
\label{g301}

The holomorphic differentials on $R_g$ can be written as
$\frac{p(E)}{\sqrt{R(E)}} dE$ (compare with \ref{aa.3}), 
where $p$ is a polynomial 
of degree $\leq g-1$ ,
\bd
R(E) =  \prod_{j=0}^{2g+1} (E-E_j) = (E-a)(E-b)
\prod_{j=1}^g ((E-\lambda _j)^2 - p_j^2), 
\ed
and $\sqrt{R(E)}$ is defined on $R_g$ with the usual convention, that
on the upper sheet $\sqrt{R(E)} \rightarrow +\infty ,$ as  
$ E \rightarrow + \infty$.
Our first goal is to determine the canonical basis of the holomorphic 
differentials $\omega _l = \frac{r_l (E)}{\sqrt{R(E)}} dE$, which has 
to be chosen such that $\int_{\alpha _k} \omega _l = \delta _{k,l}.$

For the calculations it turns out that the following basis of 
polynomials of degree $\leq g-1$ is useful.
We denote
\be
\label{g3.5}
e_j (E) := \prod_{\scriptsize \begin{array}{c} 
        m=1 \\ m \neq j \end{array} \normalsize }^g
\frac{E - \lambda _m}{\lambda _j - \lambda _m}.
\ee

We begin with the calculation of some elementary integrals.
The first candidate is
\be
\label{g3.6}
\tilde{b}_{k,j} := 
\int_{E_{2k-2}}^{E_{2k-1}}
\frac{e_j(E)}{\sqrt{|R(E)|}} dE, 1 \leq k \leq g,
\ee
and  
\be
\label{g3.7}
\tilde{B}:=(\tilde{b}_{k,j})_{k,j=1}^g.
\ee
Here the square root in the denominator always denotes the 
positive root. In order to demonstrate how this integral will be 
analyzed (especially in the limit as $p_k \rightarrow 0$) we split
the integral in two parts. Let
$d$ be any fixed number in the interval $[E_{2k-2}, E_{2k-1}]$, 
which is away from the edges of the bands independently of the 
choice of the parameters as they vary over the allowed regions.
Let
\be
\label{g3.10}
f_{j,k} (E) :=
e_j (E) \frac{1}{\sqrt{(E-a)(b-E)}}
\prod_{\scriptsize \begin{array}{c} 
        m=1 \\ m \neq k \end{array} \normalsize }^g
\frac{1}{\sqrt{(E-\lambda _m)^2 - p_m^2}},
\ee
then there exists a smooth function $\tilde{f} _{j,k}$ such 
that we can write 
\be
\label{g3.11}
f_{j,k} (E) = f_{j,k} (\lambda _k) + (E- \lambda _k) \tilde{f} _{j,k} (E)
\ee
for $E \in [d, \lambda _k]$. By simple changes of variables we have
\bea
\int_{d}^{\lambda _k - p_k}
\frac{1}{\sqrt{(E-\lambda _k)^2 - p_k^2}} dE   &=&
\int_{1}^{\frac{\lambda_ k -d}{p_k}} 
\frac{ds}{\sqrt{s^2 - 1}} \\
&=&
\mbox{ arccosh } \left(  \frac{\lambda_ k - d}{p_k}  \right) \\
&=&
\ln 2 \frac{\lambda_ k - d}{p_k}  +
\ln \left(
\frac{1}{2} + \sqrt{ \frac{1}{4} - \frac{p_k^2}{4(\lambda _k - d)^2} } 
\right).
\eea
\bea
\int_{d}^{\lambda _k - p_k}
\frac{(E- \lambda _k) \tilde{f} _{j,k} (E)}
     {\sqrt{(E-\lambda _k)^2 - p_k^2}} dE
&=&
- \int_{0}^{\sqrt{(d-\lambda _k)^2 - p_k^2}}
\tilde{f} _{j,k} \left( \lambda _k + \sqrt{s^2 + p_k^2} \right) ds.
\eea
Note that the last quantity is a $C^{1}$ function of all the parameters,
namely $(a, b, \lambda _j, p_j)$.
Furthermore we see that the expression just depends on $p_j^2$ and 
therefore we can define it for $p_j \leq 0$ such that the smooth
dependence on the parameters is preserved and we obtain a function,
which is even in all (small) $p_j$ 's. This continuation is purely formal.
Using the abbreviation
\be
\label{g3.15}
h_k := \frac{1}{\sqrt{(\lambda _k -a)(b- \lambda _k)}}  
\prod_{\scriptsize \begin{array}{c} 
        m=1 \\ m \neq k \end{array} \normalsize }^g
\frac{1}{\sqrt{(\lambda _k - \lambda _m)^2 - p_m^2}},
\ee
we state

\lfd
\begin{lemma}
\label{lg3.1}
Given the definitions as above, then
\be
\label{g3.20}
\tilde{b} _{k,j} =
\delta _{k,j} h_{k} \ln \frac{1}{p_k} +
\delta _{k-1,j} h_{k-1} \ln \frac{1}{p_{k-1}} +
\tilde{b} _{k,j}^{(reg)} ,
\ee
where $\tilde{b} _{k,j}^{(reg)}$ is a $C^{1}$ function of the parameters
$(a, b, \lambda _j, p_j)$, which can be extended for nonpositive 
values of the $p_j$ 's in an even way, preserving the smoothness.
\end{lemma}

\begin{proof}
We split the integrals as described above. The integrals from
$E_{2k - 2}$ to $d$ can be dealt with in a similar way. 
\end{proof}

There is a second elementary integral which we wish to discuss.
Define 
\be
\label{g3.25}
\tilde{a} _{k,j} :=
\int_{\lambda _k - p_k}^{\lambda _{k} + p_{k}}
\frac{e_j(E)}{\sqrt{|R(E)|}} dE .
\ee
Utilizing the notation introduced in equation (\ref{g3.10}) 
we obtain without effort
\be
\label{g3.26}
\tilde{a} _{k,j} =
\int_{-1}^{1}
\frac{f_{j,k}(\lambda _k + p_k s)}{\sqrt{1 - s^2}} ds ,
\ee
from which it is obvious that $\tilde{a} _{k,j}$ is a smooth function 
of the parameters and can be extended to (small) nonpositive values of $p_j$
as an even and smooth function. From the formula (\ref{g3.26}) we can 
also read  off that
\be
\label{g3.27}
\tilde{a} _{k,j} =
\delta _{k,j} \pi h_k + \tilde{a} _{k,j}^{(ho)},
\ee
with $\tilde{a} _{k,j}^{(ho)} = O(p_k^2)$.( (ho) stands for higher order).
Let us introduce the matrices $\tilde{A} := (\tilde{a} _{k,j})_{k,j=1}^g$,
$\tilde{A}^{(ho)} := (\tilde{a}^{(ho)} _{k,j})_{k,j=1}^g$
and $\tilde{C} :=  \tilde{A}^{-1}$. We can again state a lemma.

\lfd
\begin{lemma}
\label{lg3.2}
Given the definitions as above. Then 
\be
\label{g3.30}
\tilde{A} =  \mbox{ diag } (\pi h_k) + \tilde{A}^{(ho)}.
\ee
\be
\label{g3.31}
\tilde{C} =  \mbox{ diag } (\frac{1}{\pi h_k} ) + \tilde{C}^{(ho)}.
\ee
All entries of the matrices are smooth in the parameters $(a, b, \lambda _j
,p_j)$ and have an even and smooth extension for 
(small) nonpositive $p_j$. Finally
we have the estimates
\bd
\forall 1\leq j, k \leq g : \tilde{a} _{k,j}^{(ho)} = O(p_k^2) \mbox{ and }
\tilde{c} _{k,j}^{(ho)} = O(p_k^2).
\ed
\end{lemma}

\begin{proof}
The claim has been shown for $\tilde{A}$. Expanding
$\tilde{C}$ in a Neumann series completes the proof.
\end{proof}

Using these two lemmas we can evaluate the corresponding integrals
over the cycles of the canonical homology basis.
Let
\be
\label{g3.35}
a_{k,j} := 
\int_{\alpha _k} 
\frac{e_j(E)}{\sqrt{R(E)}} dE  \mbox{ and } A:=(a_{k,j})_{k,j=1}^g.
\ee
\be
\label{g3.36}
b_{k,j} := 
\int_{\beta _k} 
\frac{e_j(E)}{\sqrt{R(E)}} dE  \mbox{ and } B:=(b_{k,j})_{k,j=1}^g.
\ee
Looking at Figure \ref{fa1} in the beginning of Appendix \ref{g2}
and keeping track of the signs we obtain 
the following formulae
\be
\label{g3.45}
A = 
\mbox{ diag } ( 2 (-1)^{g-k} ) \tilde{A}.
\ee
\lfd
\beqa
B &=& 
\scriptsize
\left(
\begin{array}{ccc}
1&&0 \\
\vdots & \ddots & \\
1 & \cdots & 1
\end{array}
\right)
\normalsize
\mbox{ diag } ( 2 i (-1)^{g-k} ) \tilde{B}. \nonumber \\
&=&
B^{(sing)} + B^{(reg)}  , \label{g3.40}
\eeqa
with 
\lfd
\beqa
B^{(sing)} &=& 
\mbox{ diag } (-2i (-1)^{g-k} h_k \ln p_k)  \label{g3.41}, \\
B^{(reg)}  &=&
\scriptsize
\left(
\begin{array}{ccc}
1&&0 \\
\vdots & \ddots & \\
1 & \cdots & 1
\end{array}
\right)
\normalsize
\mbox{ diag } (2i (-1)^{g-k} ) \tilde{B}^{(reg)}. \lfd \label{g3.42} 
\eeqa
We return to determining the basis of the holomorphic differentials
\be
\label{g3.50}
\omega _l = \frac{r_l (E)}{\sqrt{R(E)}} dE.
\ee
Obviously we can write $r_l (E) = \sum_{j=1}^g r_l (\lambda _j) e_j(E)$, 
i.e. using vector notation 
\be
\label{g3.51}
r = R^T e,\quad \mbox{ with } R=(r_{j,l})_{j,l=1}^g , r_{j,l} := r_l(\lambda _j).
\ee
The normalization condition for the $\omega _l$ reduces to
$\sum_{j=1}^g a_{k,j} r_{j,l} = \delta _{k,l}$ or equivalently
\be
\label{g3.52}
R = A^{-1}.
\ee
Note that the last three equations 
together with Lemma \ref{lg3.2} determine the holomorphic differential
completely. 

Finally the matrix of $\beta$ periods of the $\omega _l$ 's 
can be expressed as
\bd
\tau _{k,l} = \int_{\beta _k} \frac{r_l (E)}{\sqrt{R(E)}}dE =
\sum_{j=1}^g b_{k,j} r_{j,l}.
\ed
Therefore
\be
\label{g3.55}
\tau = B R = B A^{-1}.
\ee
We combine equations 
(\ref{g3.45}),(\ref{g3.40}),(\ref{g3.41}),(\ref{g3.42}), (\ref{g3.55})
and Lemma \ref{lg3.2} to compute $\pi i \tau$.
\bea
\pi i \tau &=&
\left[
\mbox{ diag } (2 \pi (-1)^{g-k} h_k \ln p_k) +
\scriptsize
\left(
\begin{array}{ccc}
1&&0 \\
\vdots&\ddots&  \\
1&\cdots &1 
\end{array}
\right)
\normalsize
\mbox{ diag } (-2 \pi (-1)^{g-k} ) \tilde{B}^{(reg)}
\right]  \\
&&
\times \tilde{C} \mbox{ diag } ( \frac{1}{2} (-1)^{g-k} ).
\eea
Introducing the matrices
\be
\label{g3.60}
\tau_1 :=
\mbox{ diag } ( \pi (-1)^{g-k} h_k \ln p_k) 
\tilde{C} ^{(ho)}  \mbox{ diag } ( (-1)^{g-k} ),
\ee
\be
\label{g3.61}
\tau_2 :=
\scriptsize
\left(
\begin{array}{ccc}
1&&0 \\
\vdots&\ddots&  \\
1&\cdots &1 
\end{array}
\right)
\normalsize
\mbox{ diag } (- \pi (-1)^{g-k} ) \tilde{B}^{(reg)}
\tilde{C} \mbox{ diag } ( (-1)^{g-k} ),  
\ee
\be
\label{g3.63}
\tau^{(reg)} := \tau_1 + \tau_2,
\ee
we can summarize the information about $\pi i \tau$ as follows.

\lfd
\begin{lemma}
\label{lg3.3}
Given the definitions made above, we can express
\bd
\pi i \tau =
\mbox{ diag } (\ln p_k) + \tau^{(reg)} ,
\ed
where the entries of the matrix $\tau^{(reg)}$ 
are $C^{1}$ functions of the parameters $(a, b, \lambda _j, p_j)$.
Furthermore they have an even extension for nonpositive values of $p_j$,
which preserves the $C^{1}$ regularity.
\end{lemma}

\begin{proof}
It is only the matrix $\tau_1$, which needs some consideration,
namely
the questions of a smooth and even extension for the terms 
$(\ln p_k) \tilde{c}^{(ho)}_{k,j}$. But this is an immediate consequence
of Lemma \ref{lg3.2}.
\end{proof}

\section{Frequencies and phase - $U, V$ and $Z$}
\label{g302}

As already remarked in Appendix \ref{g2} (\ref{g2.53}), the vector
$U$ can be expressed as an integral of the $\omega _l$ 's.
More precisely
\bea
U &=& \frac{1}{2 \pi i} \int_{\beta} \omega ^{(1)}  \\
&=&
\int_{P_{\infty}^*}^{P_{\infty}} \omega  \\
&=&
R^T \left( 2 
\int_b^{\infty} \frac{e(E)}{\sqrt{R(E)}} dE \right).
\eea
The last integration is performed on the upper sheet along the
real axis. Using equations  
(\ref{g3.45}),  and (\ref{g3.52}) we conclude the following.
\lfd
\begin{lemma}
\label{lg3.4}
Given the notation as above. Then
\be
\label{g3.65}
U =
\mbox{ diag } ( \frac{1}{2} (-1)^{g-k} ) \tilde{C}^T  \left(2
\int_b^{\infty}  \frac{e(E)}{\sqrt{R(E)}} dE \right).
\ee
$U$ is a smooth function of all the parameters 
$(a, b, \lambda _j, p_j)$.
Again there is an extension for nonpositive values of the $p_j$ 's, which
is even in each $p_j$ and preserves the smoothness.
\end{lemma}

\begin{proof}
The methods are the same as in all the other lemmas. Looking at the formula
one might expect some difficulties for the differentiation with respect to
$b$, but the transformation $E \rightarrow E+b $ resolves the situation. 
Secondly one might wonder whether the infinite domain of integration causes
problems, but it is easy to check that all derivatives have uniformly 
at least as good a decay as $\frac{1}{E^2}$, which is integrable.
\end{proof}

In order to evaluate $V$ (see (\ref{g2n.2})), we have to determine the zero
order coefficient of $\omega _l$ at $P_{\infty}$ and $P_{\infty}^*$
(see \cite[III.3.8 (3.8.2)]{FK}).
 Denote
$r_l(E) = \sum_{j=0}^{g-1} d_{l,j} E^j$ and $\xi = \frac{1}{E}$ 
as the local coordinates at the infinities. Then
\bd
\omega _l = -
\frac{\sum_{j=0}^{g-1} d_{l,j} \xi^{g-j-1}}
     {\sqrt{\prod_{j=0}^{2g+1} (1 - \xi E_j)}} d\xi.
\ed
This implies
\bea
\omega _l =  (- d_{l,g-1} +  O(\xi)) d \xi \mbox{ at } P_{\infty},&& \\
\omega _l =  ( d_{l,g-1} +  O(\xi)) d \xi \mbox{ at } P_{\infty}^*.
\eea
We obtain
\bd
V_l = \frac{1}{2 \pi i} \int_{\beta} \omega ^{(2)} = - d_{l,g-1}.
\ed
Introducing the abbreviation
\be
\label{g3.70}
g_k := 
\prod_{\scriptsize \begin{array}{c} 
        m=1 \\ m \neq k \end{array} \normalsize }^g
\frac{1}{\lambda _k - \lambda _m},
\ee
it follows from (\ref{g3.5}) and (\ref{g3.50}) -- (\ref{g3.52}),  that
\bd
V_l  =  - \sum_{j=1}^g r_{j,l} g_j.
\ed
We can rewrite this equation as
\be
\label{g3.71}
V = -R^T \mbox{ diag } (g_k) 
\scriptsize
\left(
\begin{array}{c}
1 \\
\vdots \\
1 \\
\end{array}
\right).
\normalsize
\ee
\be
\label{g3.72}
V = - \mbox{ diag } ( \frac{1}{2} (-1)^{g-k} )
\tilde{C}^T \mbox{ diag } (g_k) 
\scriptsize
\left(
\begin{array}{c}
1 \\
\vdots \\
1 \\
\end{array}
\right).
\normalsize
\ee
The next lemma will finally tell us that we can choose the 
$\lambda _j$ 's as functions of the other parameters such that
the solution of the Toda lattice which we have constructed in
Appendix \ref{g2} is time periodic with frequency $\gamma$.

\lfd
\begin{lemma}
\label{lg3.5}
There are smooth functions $\lambda _j (a, b, p_1, \ldots, p_g), 1 \leq j
\leq g$, which are even in each $p_k$ such that
\be
\label{g3.73}
V = -\frac{1}{2 \pi}
\scriptsize
\left(
\begin{array}{c}
\gamma \\
\vdots \\
g \gamma \\
\end{array}
\right).
\normalsize
\ee
\end{lemma}

\begin{proof}
This lemma is a consequence of the implicit function theorem. To see
this more explicitly, we use Lemma \ref{lg3.2} and equations 
(\ref{g3.15}), (\ref{g3.70}) to evaluate the formula (\ref{g3.72}).
\be
\label{g3.74}
V
=    \left[
\mbox{ diag } \left(
- \frac{1}{2 \pi} \sqrt{(\lambda _k -a)(b - \lambda _k)} 
\prod_{\scriptsize \begin{array}{c} 
        m=1 \\ m \neq k \end{array} \normalsize }^g
\sqrt{1 - \scriptsize \frac{p_m^2}{(\lambda _k - \lambda _m)^2} \normalsize}
\right) + O(p^2)  \right]
\left(
\begin{array}{c}
1 \\
\vdots \\
1 \\
\end{array}
\right).
\ee
If all the $p_k$ 's are equal to zero, we choose
\bd
\lambda _j = \lambda _j^{(0)} = \frac{a+b}{2} - 
\sqrt{ \left( \frac{b-a}{2} \right)^2 - j^2 \gamma ^2 }
\ed
($\lambda _j^{(0)}$ is real, see (\ref{g3.3}) and remark below),
in order to solve equation (\ref{g3.73}). 
Observe that the $\lambda _j^{(0)}$ 's are distinct and lie in 
$ (a, \frac{a+b}{2}) \subset (a, b)$.
We compute the derivative
\bd
\partial _{\lambda} V (p=0,a,b,\lambda _j^{(0)}) =
\mbox{ diag } \left(- 
\frac{ -2 \lambda _k^{(0)} +a +b}
{4 \pi \sqrt{(\lambda _k^{(0)} -a)(b - \lambda _k^{(0)})}}
\right),
\ed
which is invertible as $k \gamma < 2 e^{-\frac{c}{2}}$ for all
$1 \leq k \leq g$ and one can therefore choose $\epsilon _a$ and
$\epsilon _b$ in (\ref{g3.3a}), (\ref{g3.3b}) small enough such that 
$k \gamma < \frac{b-a}{2}$ for all $1 \leq k \leq g$. 
The evenness of the $\lambda _j$ 's in the $p_k$ 's
is a consequence of the evenness of all the terms in the equation and
the uniqueness of the solution.
\end{proof}
\lfd
\begin{remark}
\label{rg3.12}
We have chosen for each $\lambda _j^{(0)}$ only one of the two 
possible solutions
of the equation (\ref{g3.73}) 
for $p = 0$. We have seen in 
Remark \ref{rg3.3}, that the physical
reason for this lies in the direction in which energy is transported in
the corresponding $g$-gap solution. The reader may recall that
exactly the same situation occurred in Chapter \ref{l} with the choice
of the spatial frequency $\beta$.
\end{remark}

Next we compute the phase $Z$. Therefore it is necessary to determine
the vector of Riemann constants $K$. Proceeding as in \cite[VII.1.2]{FK}
we obtain
\bd
K = - \sum _{m=1}^g A(E_{2m -1}),
\ed
as it is easy to check that the zeros of 
$P \rightarrow \vartheta (A(P))$ are given by 
$E_{2m -1}, 1 \leq m \leq g$.
This implies 
\lfd
\beqa
Z &=& \sum _{m=1}^g A(P_m) + K \label{g3.82a} \\
&=& \sum_{m=1}^g \int_{E_{2m - 1}}^{P_m} \omega =
R^T \sum_{m=1}^g \int_{E_{2m - 1}}^{P_m} \frac{e(E)}{\sqrt{R(E)}} dE
\nonumber \\
&=&
- \mbox{ diag } ( \frac{1}{2} (-1)^{g-k} )
\tilde{C}^T  \left\{
\int_{E_{2m-1}}^{P_m} \frac{e_j(E)}{\sqrt{R(E)}} dE
\right\}_{j,m=1}^g
\scriptsize
\left(
\begin{array}{c}
1 \\
\vdots \\
1 \\
\end{array}
\right).
\normalsize
\lfd \label{g3.82}
\eeqa
The path of integration from $E_{2m-1}$ to $P_m$ is chosen along the real axis,
beginning on the upper sheet and in case that $P_m$ lies on the lower sheet
we switch sheets at $E_{2m}$.
This choice is consistent with the description at
the beginning of the proof of Theorem \ref{tg2.1}.

The following lemma states that we can
choose the $P_j$ 's such that all phases $Z$ are obtained. Recall that 
$P_j$ is a point on $R_g$ such that $\pi (P_j) \in [E_{2j-1},E_{2j}]$,
(see (\ref{g3.3f})) i.e. each $P_j$ lies on a cycle diffeomorhic to $S^1$.

\lfd
\begin{lemma}
\label{lg3.6}
The map $(P_1,\ldots,P_g) \longmapsto Z$ is a surjective map
from $(S^1)^g$ to $\relz ^g / \ganz ^g$
for all choices of parameters $(a, b, \lambda _j, p_j)$. 
\end{lemma}
\begin{proof}
By inspection of equation (\ref{g3.82}) it is clear that $Z
\in \relz$ (see Lemma \ref{lg3.2}).
It suffices to show that the image of the map is open and closed.
First of all we have to convince ourselves that the map is 
differentiable. The only problem may occur at points 
$P_m = E_{2m-1}$ as the path of integration changes discontinously 
at this point. More precisely, we have to investigate what 
happens if we add to
the path of integration a cycle $\alpha _m$.
By equation (\ref{g3.82a}) and the definition of $\omega$ 
we see that this adds a vector to
$Z$, which is 
zero in $\relz ^g / \ganz ^g$. This settles the question of
differentiability.
Furthermore we know (see \cite[III.11.11 (Remark 1)]{FK} 
or \cite[Thm 17.20]{Ges}) that the differential of the
map has maximal rank $g$ at any point. This
shows that the image is open. That the image is closed is a consequence 
of the well known fact that the image of a compact set under a continuous
map is again compact and hence closed.
\end{proof}
\lfd
\begin{remark}
\label{rg3.5}
The zeros of $\vartheta (A(\cdot) + Un + Vt - Z | \tau )$.
\end{remark}
We now prove the asssertion, that the function
$P \mapsto \vartheta (A(P) + Un + Vt - Z | \tau )$ has exactly 
$g$ one zero in each gap  for all $n \in \relz, t \in \relz$ 
,which was used in the proof of Theorem \ref{tg2.1}.
By \cite[VI.3.3, Theorem b]{FK} it suffices  to show that there exists a 
choice of $P_j'$ with $\pi (P_j') \in [E_{2j-1}, E_{2j}]$,
such that 
\bd
Un + Vt - Z(P_1, \ldots , P_g) = - Z(P_1', \ldots , P_g') 
\mbox{ in } \relz ^g / \ganz ^g.
\ed 
As $U, V, Z(P_1, \ldots , P_g)$ are easily seen to be real valued vectors,
the existence of such points follows from the Lemma \ref{lg3.6} above. 

\section{$I_0^{\pm}$ and $R_0^{\pm}$}
\label{g303}

{\bf $I_0^{\pm}$:} \newline
Recall the definition of $I_0^{\pm}$ in (\ref{g2.9}). Let 
\be
\label{g3.90}
\omega_0^{(1)} :=
- \frac{\prod_{j=1}^g (E - \lambda _j)}{\sqrt{R(E)}} dE.
\ee
This differential has the desired behavior at the poles and hence we 
only have to normalize in order to obtain $\omega ^{(1)}$.
\be
\label{g3.91}
\omega^{(1)} = \omega_0^{(1)} -
\sum_{j=1}^g \left( \int_{\alpha _j} \omega_0^{(1)} \right) \omega _j.
\ee
Note that $\omega ^{(1)}$ just changes sign if we switch from the 
upper to the lower sheet, and therefore
$\forall s \geq 0: I_s^+ = I_s^-$ (see equation (\ref{g2.9})).
We can express
\bea
I_0^+ &=& \lim_{E \rightarrow \infty}
- \left( \int_b^E \omega ^{(1)} + \ln E \right)  \\
&=&
\lim_{E \rightarrow \infty}
- \left( \int_b^E \omega_0 ^{(1)} + \ln E \right) +   
< \int_{\alpha} \omega_0^{(1)}, \frac{1}{2} U >.
\eea
The above integration takes place on the upper sheet.
We recall, that the $g$-gap solution should be brought into the
form $x_n (t) = c n +$ small. Comparing with (\ref{g2.11}) and (\ref{g2.3}),
we see that this implies a condition on the following quantity.
\be
\label{g3.95}
I := I_0^+ + I_0^- .
\ee

\lfd
\begin{lemma}
\label{lg3.7}
\be
\label{g3.96}
I =
\lim_{E \rightarrow \infty}
-2 \left( \int_b^E \omega_0 ^{(1)} + \ln E \right) +
< \int_{\alpha} \omega_0^{(1)}, U >
\ee
is a smooth function of all the parameters $(a, b, \lambda _j, p_j)$
and has an even and smooth extension for nonpositive values of 
the $p_j$ 's. Furthermore there exists a unique smooth function
$b(a,p_1,\ldots,p_g)$, which is even in each $p_k$, such that
\be
\label{g3.97}
I(a,b, \lambda _j, p_j) = c.
\ee
It is understood that $\lambda _j$ is a function of $a,b, p_k$
as determined in Lemma \ref{lg3.5}.
\end{lemma}

\begin{proof}
We break $I$ into several pieces which we examine seperately.
During this proof it is understood that all the square roots  that 
appear take on positive values. Let us first examine the dependence
of $I$ on the parameters.
\begin{itemize}
\item
\bea
&& - \int_b^E \frac{d \lambda}{\sqrt{(\lambda - a)(\lambda - b)}}
+ \ln E \\
&=&
- \mbox{ arccosh } \left(
\frac{2}{b-a} \left( E- \frac{a+b}{2} \right) \right) + \ln E \\
&=&
- \left(
\ln E + \ln \left( \frac{4}{b-a} \left(1 - \frac{a+b}{2E} \right) \right) +
\ln \left( 
\frac{1}{2} + \sqrt{ \frac{1}{4} - 
\frac{(b-a)^2}{4(a+b-2E)^2}}   \right)
\right)  \\ 
&& + \ln E,
\eea
which tends to $\ln \frac{b-a}{4}$ as $E \rightarrow \infty$.
\item
The remainder of the first term of $I$ 
in equation (\ref{g3.96}) is given by
\bd
\int_b^{\infty} \frac{1}{\sqrt{(\lambda - a)(\lambda - b)}} \left(
\prod_{j=1}^g
\frac{1}{\sqrt{1 - \frac{p_j^2}{(\lambda - \lambda _j)^2}}} 
 -1 \right) d \lambda.
\ed
In order to faciliate the differentiaton with respect to $b$, we 
shift the variable of integration by $b$. Furthermore it is not
too complicated to see that all the derivatives of the integrands with
respect to the parameters are uniformly bounded by $O(\frac{1}
{\lambda ^3})$ on the domain of integration and hence integrable.
\item
For the last term $< \int_{\alpha} \omega_0^{(1)}, U >$,
it is enough to cite Lemma \ref{lg3.4} and to refer to the 
techniques introduced in the proof of 
Lemma \ref{lg3.2}.
\end{itemize}
Our second goal is to determine $b$ from the implicit function theorem.
Let $\lambda _j$ depend on $a, b, p_k$ as given in Lemma \ref{lg3.5}.
Then we can write $I = I(a, b, p_j)$, a smooth function satisfying
\bd
I(a, b, 0) = - 2 \ln \frac{b-a}{4}.
\ed
Let
\be
\label{g3.98} 
b^{(0)} := a + 4 e^{- \frac{c}{2} }, 
\ee
then
\bea
I(a, b^{(0)}, 0) &=& c. \\
\frac{\partial}{\partial b} I(a, b^{(0)} , 0) &=& \frac{e^{\frac{c}{2} }}{2}
\neq 0.
\eea
The implicit function theorem yields the remaining claims of the Lemma.
\end{proof}

{\bf $R_0^{\pm}:$}  \newline
Recall the definition of $R_0^{\pm}$ in (\ref{g2.10}). We proceed as above. Define
\be
\label{g3.100}
\omega_0^{(2)} :=
- \frac{\left(E - \frac{a+b}{2} \right) \prod_{j=1}^g (E - \lambda _j)}
{2 \sqrt{R(E)}} dE.
\ee
Then
\be
\label{g3.101}
\omega^{(2)} = \omega_0^{(2)} -
\sum_{j=1}^g \left( \int_{\alpha _j} \omega_0^{(2)} \right) \omega _j.
\ee
Equation (\ref{g2.10}) yields
\bea
R_0^+ = R_0^- &=& \lim_{E \rightarrow \infty}
 \left( -2 \int_b^E \omega ^{(2)} - E \right)  \\
&=&
\lim_{E \rightarrow \infty}
\left( -2 \int_b^E \omega_0 ^{(2)} - E \right) +   
< \int_{\alpha} \omega_0^{(2)}, U >.
\eea
The following lemma states that we can choose the parameter $a$
in such a way, that the average speed of the solution will be $0$.
\lfd
\begin{lemma}
\label{lg3.8}
Denote
\be
\label{g3.105}
R := R_0^+ + R_0^- .
\ee
Then 
\be
\label{g3.106}
R =
\lim_{E \rightarrow \infty}
\left( - 4 \int_b^E \omega_0 ^{(2)} - 2E \right) +
2 < \int_{\alpha} \omega_0^{(2)}, U >.
\ee
There exists a unique smooth function $a(p_1,\ldots,p_g)$,
which is even in each $p_k$, such that
\be
\label{g3.107}
R(a, b, \lambda _j, p_j) = 0,
\ee
where it is understood, that the $\lambda _j$ 's depend on 
$a, b$ and the $p_k$ 's as it is determined in Lemma \ref{lg3.5}
and $b$ depends on $a, p_k$ as described in the previous Lemma.
\end{lemma}

\begin{proof}
We proceed as in the proof of Lemma \ref{lg3.7}.
The only difference is that the first term in equation (\ref{g3.106})
is split up in
\bd
2 \int_b^E \frac{\lambda - \frac{a+b}{2}}
{\sqrt{(\lambda - a)(\lambda - b)}} d \lambda  -2E
\ed
and a remainder
\bd
2 \int_b^{\infty} \frac{\lambda - \frac{a+b}{2}}
{\sqrt{(\lambda - a)(\lambda - b)}} \left(
\prod_{j=1}^g
\frac{1}{\sqrt{1 - \frac{p_j^2}{(\lambda - \lambda _j)^2}}} 
 -1 \right) d \lambda.
\ed
The first part tends to $-(a+b) $, as $E \rightarrow \infty$,
whereas the second part can be dealt with as in the previous lemma.
Furthermore this part vanishes identically for $p_1 = \cdots = p_g = 0$.
With the techniques of the proof of Lemma \ref{lg3.2} we observe
that $ \int_{\alpha _k} \omega_0^{(2)} = O(p_k^2) $.
Finally we examine $R$ as a function of $a$ and $p_k$.
From what was just said it is immediate, that
\bd
R(a, p_j = 0) = - (a + b(a,0)) 
= - (2a + 4 e^{- \frac{c}{2}}).
\ed
defining $a^{(0)} := -2 e^{- \frac{c}{2}}$ we obtain
\bea
R(a^{(0)}, 0) &=& 0, \\
\frac{\partial}{\partial a} R(a^{(0)} , 0) &=& -2,
\eea
which allows the use of the implicit function theorem.
\end{proof}

Note that (\ref{g2.11}), (\ref{g2.3}) together with Lemmas \ref{lg3.3}, 
\ref{lg3.4}, \ref{lg3.5}, \ref{lg3.6}, \ref{lg3.7}, \ref{lg3.8}  
already prove Theorem \ref{tg1.5}
with the exception of the last statement, concerning the limits
$p_j \rightarrow 0$.

\section{The closing of gaps}
\label{g304}

All the formulae of the last sections were derived under the 
assumption that all the
the gaps are open, i.e. for all $1 \leq j \leq g: p_j > 0.$
Nevertheless we have made it a point that almost all of the 
analytical expressions which we have calculated have  a smooth
(at least $C^1$)
continuation for nonpositive values of the $p_j$ 's. It is not 
a priori clear that the limits we obtain if we let some of the 
$p_j$ 's tend to zero, coincide with the formulae for the corresponding
lower gap solution. It is the goal of this section to convince
ourselves that the expression for $\xi _0^- (n,t)$ which is given
in equation (\ref{g2.11}) has a continuous limit if some or all gaps
close as described above. Clearly it suffices to
examine the case that only one of the gaps closes at a time.

To be more specific, fix
$p_j > 0$ for $1 \leq j \leq g, j \neq \nu$ and $p_{\nu} = 0$, 
as well as the other parameters $a, b, \lambda _j,P_j$. We further
assume that $g > 1$. The case that the last gap closes will be dealt
with at the end of this subsection.

We use the following notation. All quantities which have to be 
evaluated will in two versions, namely with or without $'$
(e.g. $\tilde{A}, \tilde{A}'$). Without $'$ denotes the quantity  
for the g-gap case, in the limit $p_{\nu} \rightarrow 0$, whereas 
the quantity with $'$ stands for the corresponding $g-1$ gap expression
with the same choice of the remaining parameters. Note that the 
parameters $\lambda _{\nu}, p_{\nu}$ and $P_{\nu}$ do not appear in
the $g-1$ gap case.

For a $g \times g$ matrix $M$, we define by $M_{k,j \neq \nu}$ 
the $(g-1) \times (g-1)$ matrix which is obtained from $M$ by cancelling the 
$\nu$ -th row and column. Similarily for a vector
$v$, we denote by $v_{j \neq \nu}$ the vector where the 
$\nu$ -th entry is cancelled.
Finally the $(k,j)$ entry of a matrix $M$ will sometimes be denoted by
$M(k,j)$.
Let us now collect all the technical information, we 
will need.

\lfd
\begin{proposition}
\label{pg3.1}

With the notation introduced above, the following relations
hold for $g > 1$.
\begin{itemize}
\item[(i)]
$
\tilde{A} _{k,j \neq \nu} =
\mbox{ diag } (\mbox{ sgn }(\lambda _k - \lambda _{\nu}))
\tilde{A}'
\mbox{ diag } ( \frac{1}{ \lambda _k - \lambda _{\nu} } ).
\newline
\tilde{A} (\nu,j) = 0, \mbox{ for } j \neq \nu.
$
\item[(ii)]
$
\tilde{C} _{k,j \neq \nu} =
\mbox{ diag } ( \lambda _k - \lambda _{\nu} )
\tilde{C}'
\mbox{ diag } (\mbox{ sgn }(\lambda _k - \lambda _{\nu})).
\newline
\tilde{C} (\nu,j) = 0, \mbox{ for } j \neq \nu.
$
\newline
$
\tilde{C}^{(ho)} _{k,j \neq \nu} =
\mbox{ diag } ( \lambda _k - \lambda _{\nu} )
\tilde{C}'^{(ho)}
\mbox{ diag } (\mbox{ sgn }(\lambda _k - \lambda _{\nu})).
\newline
\tilde{C}^{(ho)} (\nu,j) = 0, \mbox{ for } j \neq \nu.
$
\item[(iii)]
$
U _{j \neq \nu} = U'.
$
\item[(iv)]
$
V _{j \neq \nu} = V'.
$
\item[(v)]
$
Z _{j \neq \nu} = Z '.
$
\item[(vi)]
$
\tau _{k,j \neq \nu}^{(reg)} = \left( \tau^{(reg)} \right) '.
$
\item[(vii)]
$
I = I'.
$
\item[(viii)]
$
R = R'.
$
\end{itemize}
\end{proposition}

\begin{proof}
\begin{itemize}
\item[(i)]
\bea
j,k \neq \nu : \tilde{a} _{k,j}
&=&
\int_{\lambda _k - p_k}^{\lambda _k + p_k}
\frac{e_j(E)}{\sqrt{|R(E)|}} dE \\
&=&
\frac{1}{ \lambda _j - \lambda _{\nu} }
\int_{\lambda _k - p_k}^{\lambda _k + p_k}
\frac{e_j'(E)}{\sqrt{|R'(E)|}} 
\mbox{ sgn }(E - \lambda _{\nu}) dE. \\
k = \nu : \tilde{a} _{k,j}
&=&
\delta _{\nu ,j} \pi h_{\nu} \mbox{ (see equation (\ref{g3.27}) ) }.
\eea
\item[(ii)]
We obtain the information about $\tilde{C}$ from evaluating the 
relation $\tilde{A} \tilde{C} = I$ in the following order.
From the $(\nu, \nu)$ entry conclude that $\tilde{C} (\nu, \nu)
= \frac{1}{\pi h_{\nu}}$. Looking at the $\nu$ -th row we see that
$\forall j \neq \nu: \tilde{C}(\nu, j) = 0$. This shows that
$\tilde{C} _{k,j \neq \nu}$ is the inverse of $\tilde{A} _{k,j \neq \nu}$,
and this yields the claim for $\tilde{C}$. To prove the claim for
$\tilde{C}^{(ho)}$ it is enough (see (\ref{g3.31})) to observe from
(\ref{g3.15}) that for $k \neq \nu$
\bd
\frac{1}{\pi h_k} = (\lambda _k - \lambda _{\nu})
\frac{1}{\pi h_k'} \mbox{ sgn } (\lambda _k - \lambda _{\nu}).
\ed
\item[(iii)]
For $j \neq \nu$, we check that 
$\int_b^{\infty} \frac{e_j(E)}{\sqrt{R(E)}} dE  =
\frac{1}{ \lambda _j - \lambda _{\nu} }
\int_b^{\infty} \frac{e_j'(E)}{\sqrt{R'(E)}} dE$. The 
claim then follows from equation (\ref{g3.65}) and from
property (ii) (esp. $\tilde{C}^T(k,\nu) = 0$ for $k \neq \nu)$.
\item[(iv)]
The proof is similar to (iii). Here it suffices 
by equation (\ref{g3.72}) to show that
\bd
\left[
\mbox{ diag }(g_k)
\scriptsize
\left(
\begin{array}{c}
1 \\
\vdots \\
1 
\end{array}
\right)
\normalsize
\right]_{j \neq \nu} =
\mbox{ diag }(\frac{1}{\lambda _k - \lambda _{\nu}})
\mbox{ diag }(g_k')
\scriptsize
\left(
\begin{array}{c}
1 \\
\vdots \\
1 
\end{array}
\right).
\normalsize
\ed
But this follows directly from (\ref{g3.70}).
\item[(v)]
Looking at formula (\ref{g3.82}) it is sufficient to verify that
\bea
&&
\left\{
\left[
\left(
\int_{E_{2m-1}}^{P_m}
\frac{e_j(E)}{\sqrt{R(E)}} dE
\right)_{j,m}
\right]
\scriptsize
\left(
\begin{array}{c}
1 \\
\vdots \\
1 
\end{array}
\right)
\normalsize
\right\}_{j \neq \nu} \\
&=&
\mbox{ diag }(\frac{1}{\lambda _k - \lambda _{\nu}})
\left[
\left(
\int_{E_{2m-1}}^{P_m}
\frac{e_j'(E)}{\sqrt{R'(E)}} dE
\right)_{j,m \neq \nu}
\right]
\scriptsize
\left(
\begin{array}{c}
1 \\
\vdots \\
1 
\end{array}
\right).
\normalsize
\eea
It is not difficult to check this relation, if one observes that
for $j \neq \nu:
\int_{E_{2\nu}}^{P_{\nu}}
\frac{e_j(E)}{\sqrt{R(E)}} dE = 0$ 
in the limit of $p_{\nu} \rightarrow 0$.

\item[(vi)]
We have to examine the behavior of $\tau _1$ and of
$\tau _2$ (see equations (\ref{g3.60}),(\ref{g3.61})). 
The fact that
$\left( \tau _1 \right) _{k,j \neq \nu} = \left( \tau _1 \right) '$
follows directly from (ii) above and (\ref{g3.60}), using the fact that
for $k \neq \nu: h_k = \frac{1}{| \lambda _k - \lambda _{\nu} |} h_k'$.
The proof that 
$\left( \tau _2 \right) _{k,j \neq \nu} = \left( \tau _2 \right) '$
needs more consideration. First we remark, that
because of $\tilde{C}(\nu, j) = 0$ for $j \neq \nu$  (see (ii)), it 
suffices to prove the following equality. 
\bea
&&\left[
\scriptsize
\left(
\begin{array}{ccc}
1 & & 0 \\
\vdots& \ddots &  \\
1&\cdots &1 
\end{array}
\right)
\normalsize
\mbox{ diag }(-2 \pi (-1)^{g-k}) \tilde{B}^{(reg)}
\right]_{k,j \neq \nu} \\
&=&
\scriptsize
\left(
\begin{array}{ccc}
1 & & 0 \\
\vdots& \ddots &  \\
1&\cdots &1
\end{array}
\right)
\normalsize
\mbox{ diag }(-2 \pi (-1)^{g-1-k}) \left( \tilde{B}^{(reg)} \right)'
\mbox{ diag }(\frac{1}{\lambda _k - \lambda _{\nu}}).
\eea
Write $\tilde{B}^{(reg)} =
\tilde{B} - \tilde{B}^{(sing)}$. Observe that by equation (\ref{g3.40})
\bea
&&\left[
\scriptsize
\left(
\begin{array}{ccc}
1 & & 0 \\
\vdots& \ddots &  \\
1&\cdots &1 
\end{array}
\right)
\normalsize
\mbox{ diag }(-2 \pi (-1)^{g-k}) \tilde{B}^{(sing)}
\right]_{k,j \neq \nu} \\
&=&
\left[
\mbox{ diag }(2 \pi (-1)^{g-k} h_k \ln p_k)
\right]_{k,j \neq \nu},
\eea
and one checks easily that this equals
\bd
\scriptsize
\left(
\begin{array}{ccc}
1 & & 0 \\
\vdots& \ddots &  \\
1&\cdots &1 
\end{array}
\right)
\normalsize
\mbox{ diag }(-2 \pi (-1)^{g-1-k}) \left( \tilde{B}^{(sing)} \right)'
\mbox{ diag }(\frac{1}{\lambda _k - \lambda _{\nu}}).
\ed
To complete the proof we must check the relation for $\tilde{B}$ itself. It 
suffices to compute $\tilde{b} _{k,j}$ for $j \neq \nu$. It 
is enough to verify the following two relations.
\begin{itemize}
\item[$\bullet$]
For $k \neq \nu$ and $k \neq \nu + 1$:
$\tilde{b}_{k,j} =
\frac{ \mbox{ sgn }(\lambda _k - \lambda _{\nu})}
     { \lambda _j - \lambda _{\nu}}
\tilde{b}'_{k',j'}$, \newline
where $k' = k$, for $k < \nu$ and  $k' = k-1$, for $k > \nu$, \newline
and $j' = j$, for $j < \nu$ and  $j' = j-1$, for $j > \nu$.
\item[$\bullet$]
$- \tilde{b}_{\nu,j} + \tilde{b}_{\nu + 1,j} =
\frac{1}{\lambda _j - \lambda _{\nu}} \tilde{b}'_{\nu,j'} $.
\end{itemize}
The first relation follows directly from the definitions. 
The crucial point for proving the second 
relation is the following calculation. In order to investigate the
limit $p_{\nu} \rightarrow 0$ explicitly, let us now assume for a moment that 
$p_{\nu} > 0$. Let $\epsilon > 0$ be an arbitrary but small number
and let $p_{\nu} < \epsilon$. Then
\bea
&&- \int_{\lambda _{\nu} - \epsilon}^{\lambda _{\nu} - p_{\nu}}
\frac{e_j(E)}{\sqrt{|R(E)|}} dE  
+ \int_{\lambda _{\nu} + p_{\nu}}^{\lambda _{\nu} + \epsilon}
\frac{e_j(E)}{\sqrt{|R(E)|}} dE  \\
&=&
\frac{1}{\lambda _j - \lambda _{\nu}}
\left(
- \int_{\lambda _{\nu} - \epsilon}^{\lambda _{\nu} - p_{\nu}}
\frac{e_j'(E)(E-\lambda _{\nu})}
     {\sqrt{|R'(E)|}  
     \sqrt{(E - \lambda _{\nu})^2-p_{\nu}^2} } dE  
\right. \\
&&  \left.
+ \int_{\lambda _{\nu} + p_{\nu}}^{\lambda _{\nu} + \epsilon}
\frac{e_j'(E)(E-\lambda _{\nu})}
     {\sqrt{|R'(E)|} 
     \sqrt{(E - \lambda _{\nu})^2-p_{\nu}^2} } dE  
\right)  \\
&=&
\frac{1}{\lambda _j - \lambda _{\nu}}
\left(
- \int_{\sqrt{\epsilon ^2 - p_{\nu}^2}}^0
\frac{e_j'(\lambda _{\nu} - \sqrt{s^2 + p_{\nu}^2})}
     {\sqrt{|R'(\lambda _{\nu} - \sqrt{s^2 + p_{\nu}^2})|}} ds 
\right. \\
&& \left.
+
\int_0^{\sqrt{\epsilon ^2 - p_{\nu}^2}}
\frac{e_j'(\lambda _{\nu} + \sqrt{s^2 + p_{\nu}^2})}
     {\sqrt{|R'(\lambda _{\nu} + \sqrt{s^2 + p_{\nu}^2})|}} ds 
\right) \\
& \rightarrow &
\frac{1}{\lambda _j - \lambda _{\nu}}
\left(
\int_{-\epsilon}^{\epsilon}
\frac{e_j'(\lambda _{\nu} + s)}
     {\sqrt{|R'(\lambda _{\nu} + s )|}} ds 
\right), \mbox{ as } p_{\nu} \rightarrow 0.
\eea
\item[(vii),] (viii) \newline
The proof for $I$ and $R$ follows trivially from
the formulae which were produced in the proofs of Lemma \ref{lg3.7} and
Lemma \ref{lg3.8},
properties (iii) and (iv) of this proposition and the observation,
that $\int_{\alpha_{\nu}} \omega _0^{(m)} \rightarrow 0$, as 
$p_{\nu} \rightarrow 0$ for $m=1,2$.
\end{itemize}
\end{proof}
{\bf Remarks:}

(1) In the above proposition we have kept $a,b, \lambda _j$ fixed as
$p_{\nu} \rightarrow 0$, but this is clearly not necessary. Indeed,
if $a,b, \lambda _j$ are given continuous functions of $p_{\nu}$
then the obvious analog of the proposition holds true.
For example in the formula (vii),
\bd
I(a,b,\lambda _j,p_1,\ldots,0,\ldots,p_g) =
I'(a,b,\lambda _j,p_{j \neq \nu}),
\ed
simply replace $a,b,\lambda _j$ by their limiting values at 
$p_{\nu} = 0$.
 
(2) During the proposition we assumed that $g > 1$, as most of the terms
do not carry any meaning for the $0$-gap case. Only the quantities 
$I$ and $R$ are also well defined for $g = 0$ and without changing
the proof of property (vii) and (viii), 
we see that the proposition also holds in the transition
from the $1$ -gap to the $0$-gap situation.

We finally turn to the question of basic interest, namely how the 
formula for the solutions $x_n(t)$ in Theorem \ref{tg1.5} (see 
equation (\ref{g3.110}))  behaves, if we let some
or all of the $p_{\nu}$ tend to zero. Recall that we have determined the 
parameters $a,b$ and $\lambda _j$ as functions of $p_1, \ldots, p_g$ such
that $x_n(t)$ is a periodic function in $t$, which is of the form
\bd
x_n (t) =
c n + \ln
\frac{\vartheta (\frac{1}{2} U - Z |\tau)
      \vartheta ((n-\frac{1}{2}) U + tV - Z |\tau)}
     {\vartheta (- \frac{1}{2} U - Z |\tau)
      \vartheta ((n+\frac{1}{2}) U + tV - Z |\tau)}.
\ed
Let us now again investigate what happens if one of the gaps closes,
i.e. $p_{\nu} \rightarrow 0$. 
The first question is, whether the choice of 
the parameters $a, b$ and $\lambda _j$ as functions of the $p_k$ 's have 
the proper limit as $p_{\nu} \rightarrow 0$, i.e. whether for example
\bd
\lim _{p_{\nu} \rightarrow 0} a(p_1, \ldots, p_{\nu}, \ldots, p_g) =
a'(p_1, \ldots, p_{\nu - 1}, p_{\nu + 1}, \ldots, p_g).
\ed
But this can be seen from properties (iv),(vii) and (viii) of  
Proposition \ref{pg3.1}, as the solution of the determining equations
(\ref{g3.73}), (\ref{g3.97}) and (\ref{g3.107}) for 
$a = a(p_1, \ldots, p_g), b = b(p_1, \ldots, p_g)$ and
$\lambda _j = \lambda _j(p_1, \ldots, p_g)$ are unique.
We are now in the position to apply the Proposition \ref{pg3.1} in order to
investigate equation (\ref{g3.110}). It suffices to look
at the behavior of the theta functions.  Using Lemma \ref{lg3.3}
and Lemma \ref{lg3.5}
we obtain
\lfd
\begin{eqnarray}
&&\vartheta ((n-\frac{1}{2}) U + tV - Z |\tau)
\nonumber \\
&=&
\sum_{l \in \ganz ^g}
p_1^{l_1^2} \cdot \ldots \cdot p_g^{l_g^2}
\exp \left(2 \pi i (<l,U>(n-\frac{1}{2})
  - <l, Z >)
  +<l, \tau^{(reg)} l> \right) \nonumber \\
&& \exp \left(
  -i <l,\tiny \left( \begin{array}{c} 1 \\ 2 \\ \vdots \\ g \end{array} 
  \right) \normalsize > \gamma t \right). \label{g3.115}
\end{eqnarray}
We see, that in the limit $p_{\nu} \rightarrow 0$ only those terms in 
the sum survive, for which $l_{\nu} = 0$. Using in addition 
the relevant results from Proposition \ref{pg3.1}
 we see that the limit of the $g$-gap theta function is the appropriate 
$g-1$-gap theta function.
The other three theta functions in equation
(\ref{g3.110}) can be dealt with in exactly the same way.  
Hence we have completed the proof of Theorem \ref{tg1.5}.

\chapter{Numerical experiments}
\label{c}
\section{Figures of lattice motion}
\label{c1}

Figures \ref{fc1}-\ref{fc6} below display the motion of the first
ten particles ($x_0$ -- $x_9$ with the zeroth particle on top
and the nineth particle on the bottom of each figure) of lattices,
which are described by the following system of equations:
\be
\ddot{x}_n = F(x_{n-1} - x_n) - F(x_n - x_{n+1}), \quad n \geq 1,
\label{cc.5}
\ee
with driver $x_0$ of the form 
\be
x_0(t) = t + \varepsilon (\sin \gamma t + 0.5 \cos 2 \gamma t)
\label{cc.10}
\ee
and initial values given at $t=0$
\be 
x_n(0) =\dot{x}_n (0) = 0, \quad  n\geq 1.
\label{cc.15}
\ee
We remark that we have also made experiments with driving particles
$x_0 (t) = 2at + h(\gamma t)$ and $h$ being periodic functions
different from type (\ref{cc.10}) (e.g. $h$ piecewise linear) and we
have always obtained results similar to those described below. 
The choice of parameters $\varepsilon, \gamma, F$ is made as follows.
On each page there are four figures. They correspond to different 
force functions $F$:
\lfd 
\bea
\mbox{ top left }&:& F(x) = e^x \mbox{ (Toda lattice) } \\
\mbox{ top right }&:& F(x) = 2.25 x \mbox{ (linear lattice) } \\
\mbox{ bottom left }&:& F(x) = 1.71 (x + 0.2 x^3)  \\
\mbox{ bottom right }&:& F(x) = \frac{2.53}{1 - 0.4 x}   
\eea
The parameters for these four types of force functions were chosen
such that in the case of $\epsilon = 0$, the system behaves subcritical
(i.e. $0.5 < a_{{\rm crit}}(F)$) and the lattice comes to rest ($x_n (t) 
\rightarrow c n$ as $t \rightarrow \infty$). For better comparison of the 
different force laws, we  ensured in addition that 
$F'(-c) \approx 2.25$ in all for cases. Therefore we expect from
the linear calculations (cf equation (\ref{e.130})) that the 
threshold values for the frequencies should be approximately the same in 
all four cases (at least for small $\epsilon$), namely 
$\gamma _k \approx 3/k$. 
Hence the different values which we selected for the driving frequency,
$\gamma = 3.1, 2.1, 1.2$ satisfy
$3.1 > \gamma _1 > 2.1 > \gamma _2 > 1.2 > \gamma _3$. 
Finally all experiments were made for two different values of the
driver' s amplitude, namely $\varepsilon = 0.2$ 
(see Figures \ref{fc1}, \ref{fc2}, \ref{fc3}) and $\varepsilon = 0.5$
(see Figures \ref{fc4}, \ref{fc5}, \ref{fc6}).

\section{Spectral densities}
\label{c2}

We consider equations (\ref{cc.5}), (\ref{cc.15}) in the 
case of the Toda lattice ($F(x) = e^x$). As described in the 
Introduction we observe numerically that the spectrum of the corresponding
Lax operator obtains a band-gap structure as $t \rightarrow \infty$
(cf Figures \ref{fe3} -- \ref{fe6}). In Chapter \ref{s} we derived under
various assumptions
an integral equation for the time-asymptotic spectral density
\be
J(\lambda)=\mathop{\lim}\limits_{t\rightarrow\infty}\,\,\,\frac{\sharp
\,\,\{{\rm eigenvalues}\,\,{\rm of}\,\,L(t)\,\,{\rm that}\,\,{\rm
are}
\,\,<\lambda\}}{t},
\label{c2.5}
\ee
(cf (\ref{e.165}), Ansatz \ref{as2.5} and Theorem \ref{ts2.5}), 
which we solved explicitly
in Theorem \ref{ts3.5}, provided a discrete number of data is known,
namely the number and endpoints of the bands. In order to test the 
results of Chapter \ref{s} we proceed as follows. 
We compute $J(\lambda)$ for $\lambda < \inf\sigma_{ess}(L(t))$
numerically by evaluating (\ref{c2.5}) for large times $t$.
This computation also yields a good approximation of
the position of the endpoints of the spectral bands. Finally
we determine the ``predicted spectral density'' by solving the 
integral equation given by (\ref{s3.40})-(\ref{s3.50}), using
the numerically obtained knowledge about the endpoints of the bands.
In Figures \ref{fc7} and \ref{fc8} below we display the results
of these experiments for two different drivers $x_0$ of type (\ref{cc.10}),
representing the cases
$\gamma _1 > \gamma > \gamma _2$ and  $\gamma _2 > \gamma > \gamma _3$,
i.e. the one-gap and the two-gap situation.
Both figures contain
\bea
\mbox{ top left }&:& \mbox{ the lattice motion } \\
\mbox{ top right }&:& 
\mbox{ the time evolution of the spectrum of the corresponding } \\
&& \mbox{ Lax
operator }\\
\mbox{ bottom left }&:& 
\mbox{ the numerically computed spectral density at } t = 200 \\
\mbox{ bottom right }&:& 
\mbox{ the ``predicted spectral density'' }
\eea
One observes very good agreement of numerical and predicted spectral
densities, which a posteriori justifies the assumptions introduced
in Chapter \ref{s}.

\clearpage
\lfd
\begin{figure}[h]
\leavevmode \epsfysize=18cm
\epsfbox{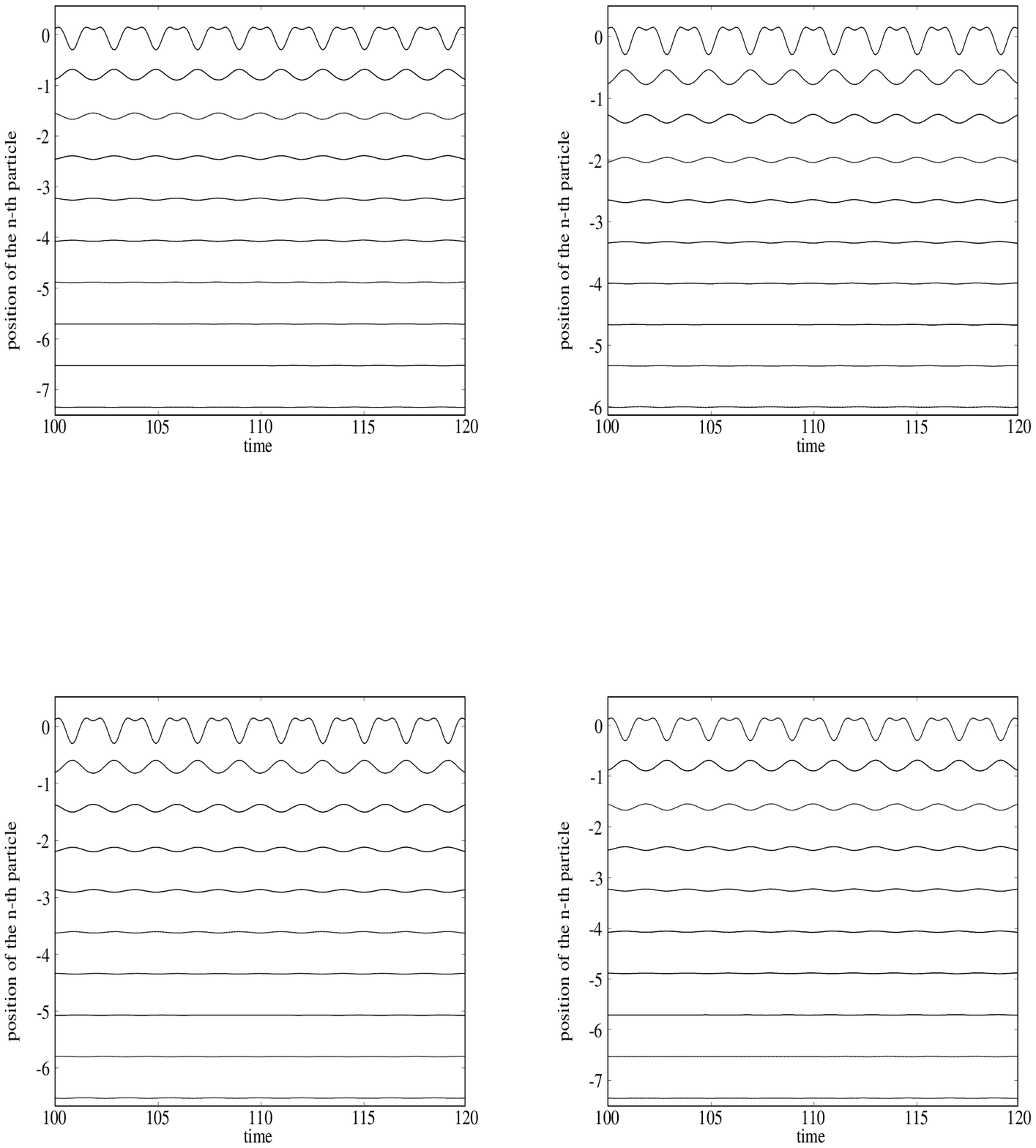}
\unitlength=1cm
\begin{picture}(13,0) 
\put(-11.5, 16.8){\makebox{$\textstyle {F(x) = e^x}$ (Toda lattice)}}
\put(-3.5, 7.8){\makebox{$\textstyle {F(x) = \frac{2.53}{1-0.4 x}}$}}
\put(-11.3, 7.8){\makebox{$\textstyle {F(x) = 1.71 (x + 0.2 x^3)}$ }}
\put(-4.75, 16.8){\makebox{$\textstyle {F(x) = 2.25 x}$ (linear lattice)}}
\end{picture}
\caption{{\em Motion of lattices (cf (C.1) -- (C.3)) with 
$\varepsilon = 0.2, \gamma = 3.1$}}
\label{fc1}
\end{figure}

\clearpage
\lfd
\begin{figure}[h]
\leavevmode \epsfysize=18cm
\epsfbox{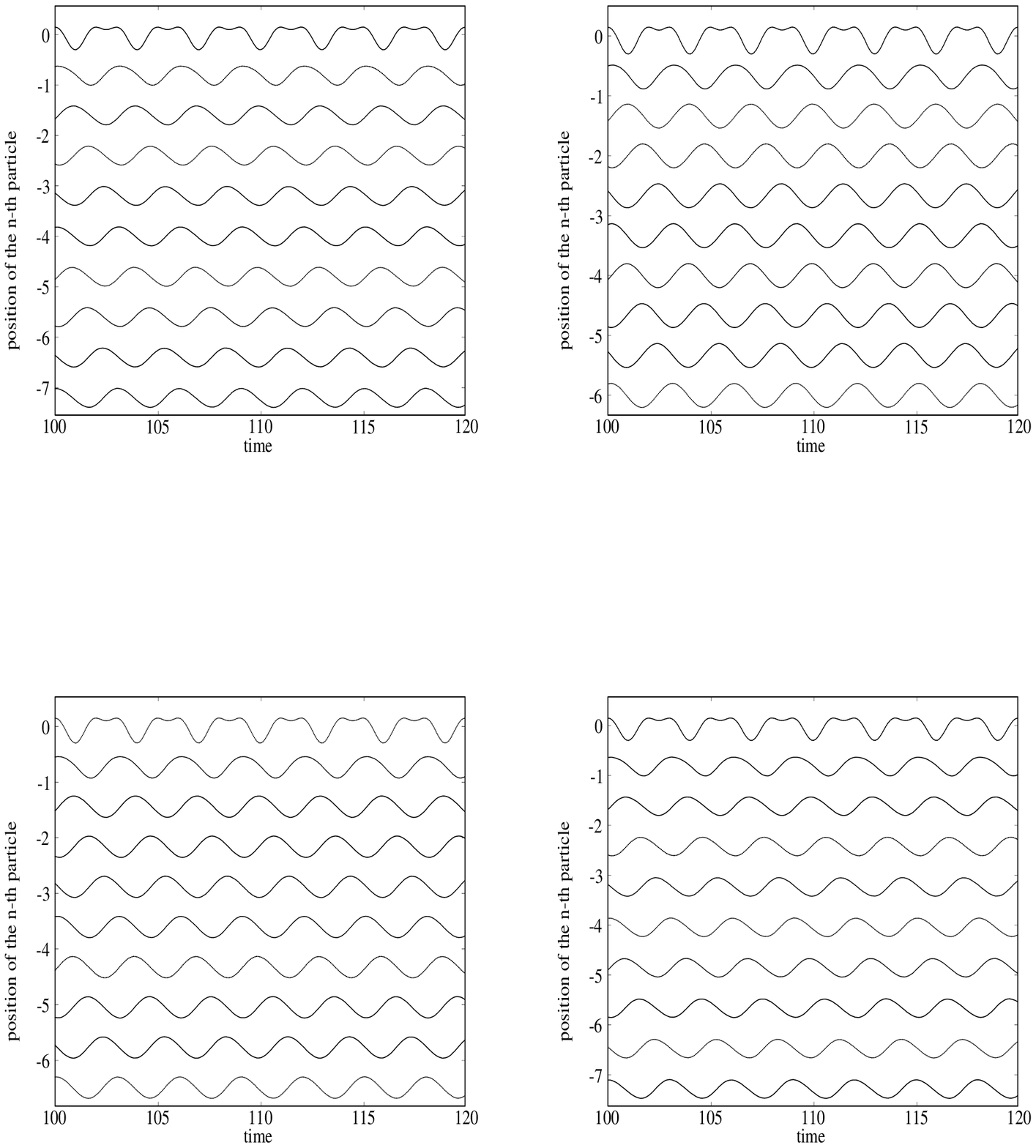}
\unitlength=1cm
\begin{picture}(13,0) 
\put(-11.5, 16.8){\makebox{$\textstyle {F(x) = e^x}$ (Toda lattice)}}
\put(-3.5, 7.8){\makebox{$\textstyle {F(x) = \frac{2.53}{1-0.4 x}}$}}
\put(-11.3, 7.8){\makebox{$\textstyle {F(x) = 1.71 (x + 0.2 x^3)}$ }}
\put(-4.75, 16.8){\makebox{$\textstyle {F(x) = 2.25 x}$ (linear lattice)}}
\end{picture}
\caption{{\em Motion of lattices (cf (C.1) -- (C.3)) with 
$\varepsilon = 0.2, \gamma = 2.1$}}
\label{fc2}
\end{figure}

\clearpage
\lfd
\begin{figure}[h]
\leavevmode \epsfysize=18cm
\epsfbox{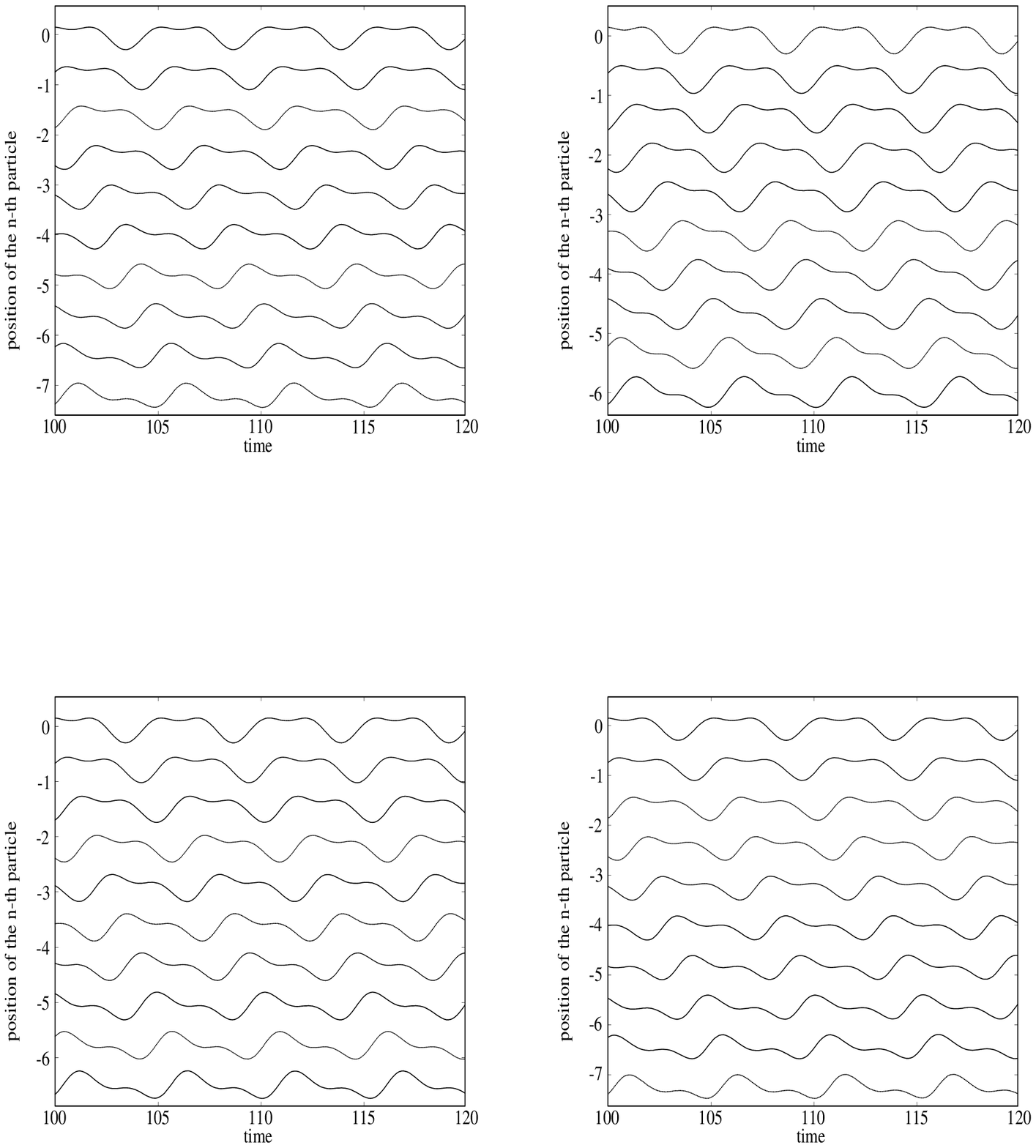}
\unitlength=1cm
\begin{picture}(13,0) 
\put(-11.5, 16.8){\makebox{$\textstyle {F(x) = e^x}$ (Toda lattice)}}
\put(-3.5, 7.8){\makebox{$\textstyle {F(x) = \frac{2.53}{1 - 0.4 x}}$}}
\put(-11.3, 7.8){\makebox{$\textstyle {F(x) = 1.71 (x + 0.2 x^3)}$ }}
\put(-4.75, 16.8){\makebox{$\textstyle {F(x) = 2.25 x}$ (linear lattice)}}
\end{picture}
\caption{{\em Motion of lattices (cf (C.1) -- (C.3)) with 
$\varepsilon = 0.2, \gamma = 1.2$}}
\label{fc3}
\end{figure}

\clearpage
\lfd
\begin{figure}[h]
\leavevmode \epsfysize=18cm
\epsfbox{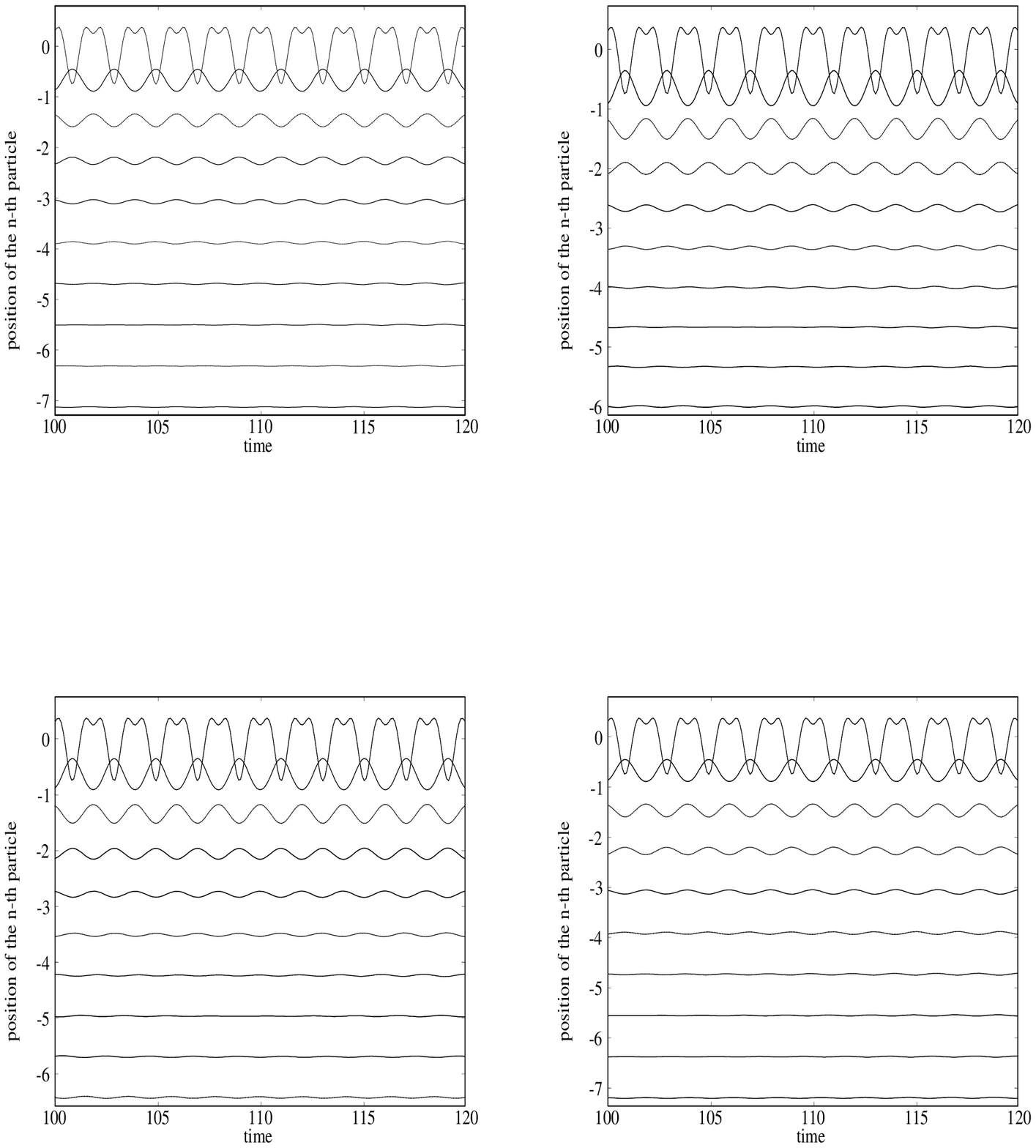}
\unitlength=1cm
\begin{picture}(13,0) 
\put(-11.5, 16.8){\makebox{$\textstyle {F(x) = e^x}$ (Toda lattice)}}
\put(-3.5, 7.8){\makebox{$\textstyle {F(x) = \frac{2.53}{1 - 0.4 x}}$}}
\put(-11.3, 7.8){\makebox{$\textstyle {F(x) = 1.71 (x + 0.2 x^3)}$ }}
\put(-4.75, 16.8){\makebox{$\textstyle {F(x) = 2.25 x}$ (linear lattice)}}
\end{picture}
\caption{{\em Motion of lattices (cf (C.1) -- (C.3)) with 
$\varepsilon = 0.5, \gamma = 3.1$}}
\label{fc4}
\end{figure}

\clearpage
\lfd
\begin{figure}[h]
\leavevmode \epsfysize=18cm
\epsfbox{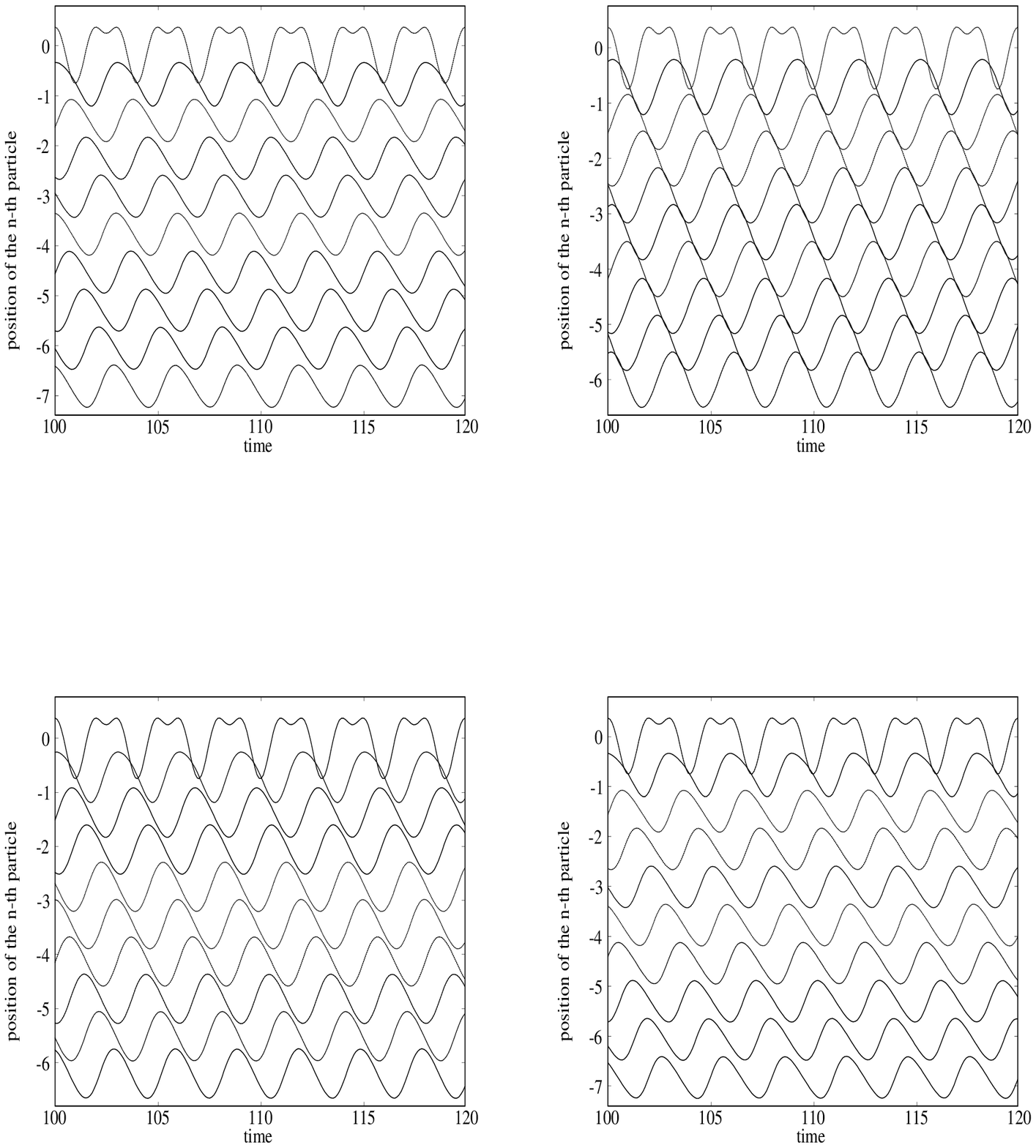}
\unitlength=1cm
\begin{picture}(13,0) 
\put(-11.5, 16.8){\makebox{$\textstyle {F(x) = e^x}$ (Toda lattice)}}
\put(-3.5, 7.8){\makebox{$\textstyle {F(x) = \frac{2.53}{1 - 0.4 x}}$}}
\put(-11.3, 7.8){\makebox{$\textstyle {F(x) = 1.71 (x + 0.2 x^3)}$ }}
\put(-4.75, 16.8){\makebox{$\textstyle {F(x) = 2.25 x}$ (linear lattice)}}
\end{picture}
\caption{{\em Motion of lattices (cf (C.1) -- (C.3)) with 
$\varepsilon = 0.5, \gamma = 2.1$}}
\label{fc5}
\end{figure}

\clearpage
\lfd
\begin{figure}[h]
\leavevmode \epsfysize=18cm
\epsfbox{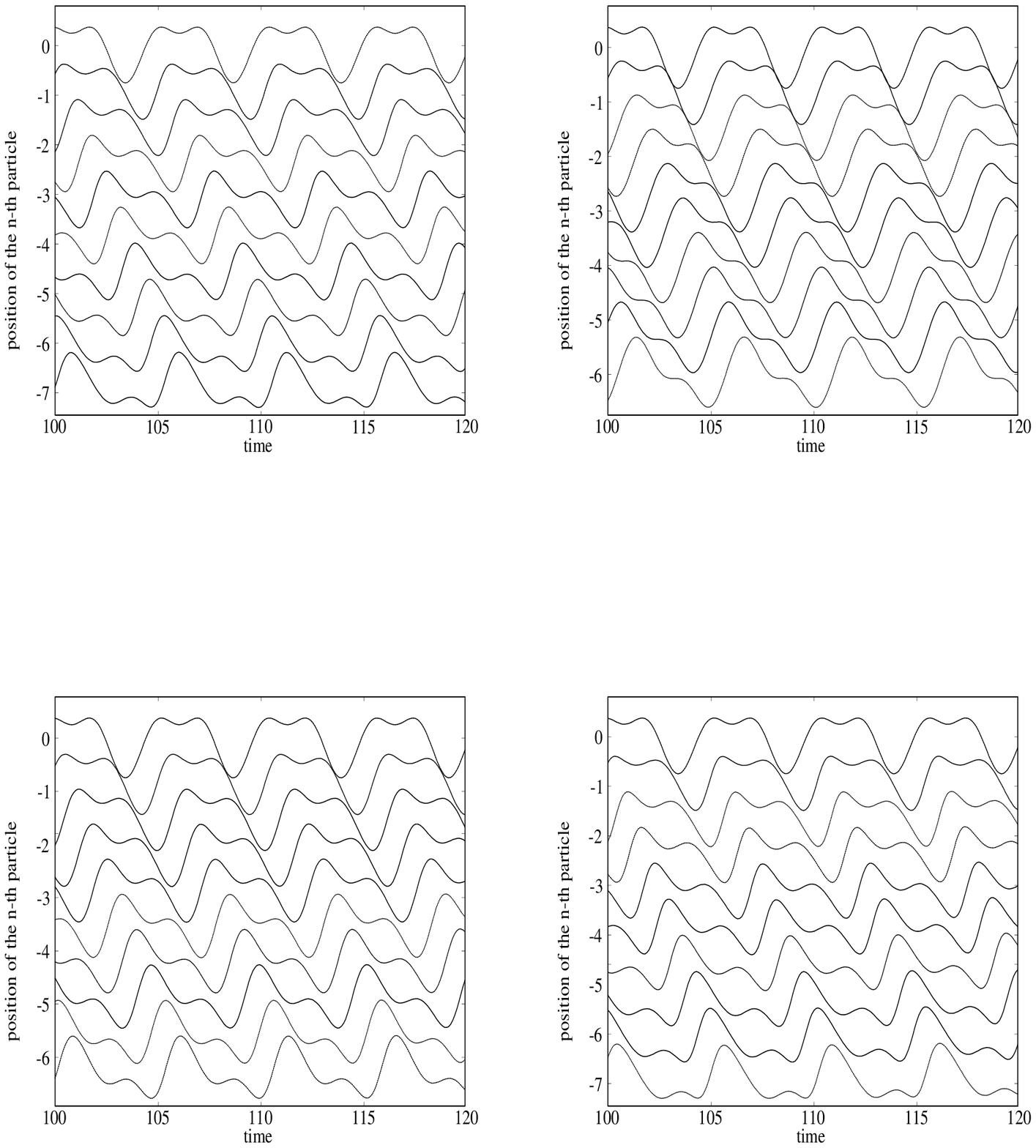}
\unitlength=1cm
\begin{picture}(13,0) 
\put(-11.5, 16.8){\makebox{$\textstyle {F(x) = e^x}$ (Toda lattice)}}
\put(-3.5, 7.8){\makebox{$\textstyle {F(x) = \frac{2.53}{1 - 0.4 x}}$}}
\put(-11.3, 7.8){\makebox{$\textstyle {F(x) = 1.71 (x + 0.2 x^3)}$ }}
\put(-4.75, 16.8){\makebox{$\textstyle {F(x) = 2.25 x}$ (linear lattice)}}
\end{picture}
\caption{{\em Motion of lattices (cf (C.1) -- (C.3)) with 
$\varepsilon = 0.5, \gamma = 1.2$}}
\label{fc6}
\end{figure}

\clearpage
\lfd
\begin{figure}[h]
\leavevmode \epsfysize=18cm
\epsfbox{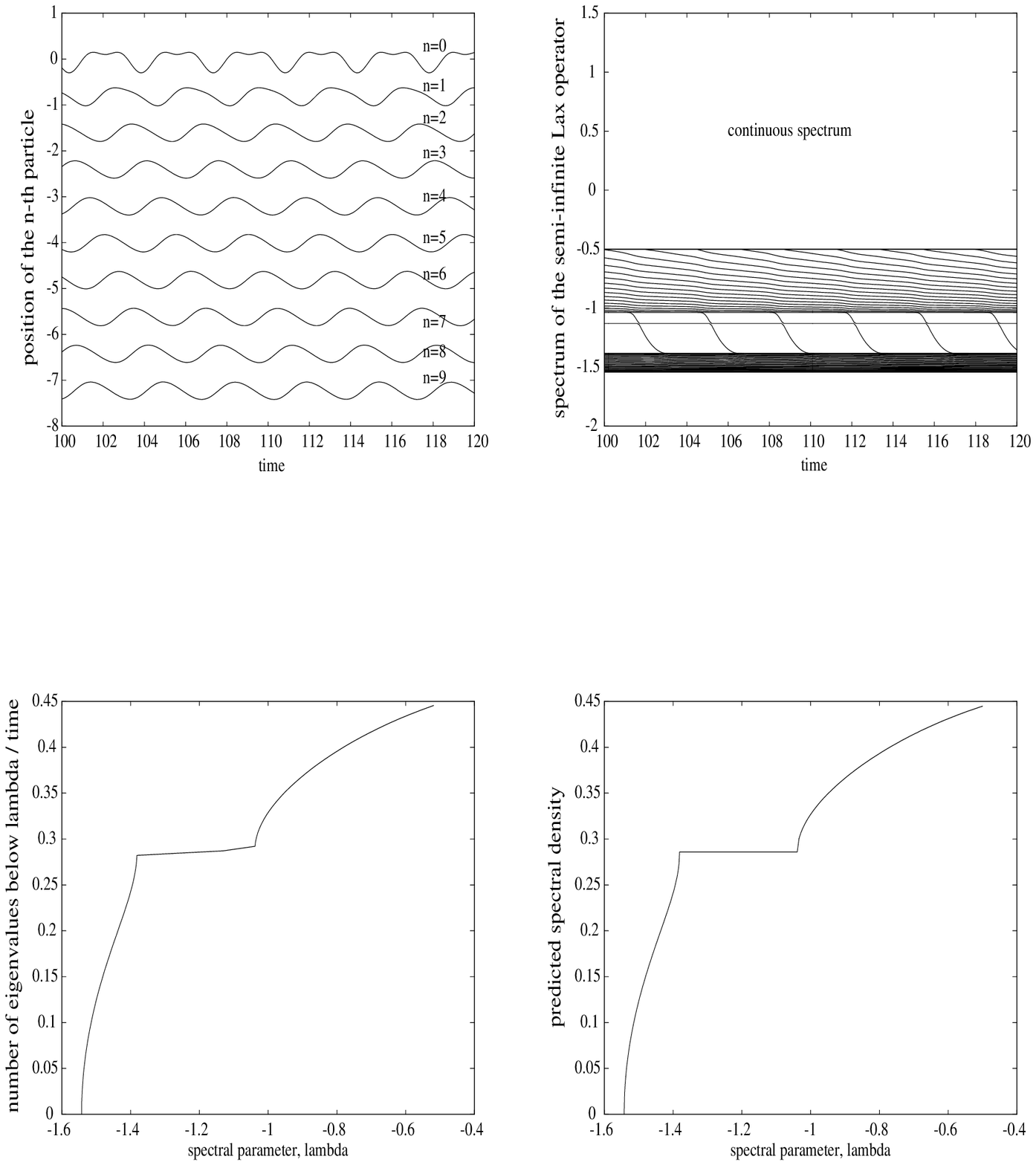}
\unitlength=1cm
\begin{picture}(13,0) 
\put(-10.5, 16.8){\makebox{Lattice motion}}
\put(-4.5, 8.0){\makebox{Predicted spectral density}}
\put(-12.75, 8.0){\makebox{Numerically observed spectral density}}
\put(-3.5, 16.8){\makebox{The spectrum}}
\end{picture}
\caption{{\em Toda lattice with driver 
$x_0(t) = t + 0.2 \sin \gamma t + 0.1 \cos \gamma t, \gamma = 1.8$, i.e.
$\gamma _1 > \gamma > \gamma _2$; see Section C.2 for detailed 
description.}}
\label{fc7}
\end{figure}

\clearpage
\lfd
\begin{figure}[h]
\leavevmode \epsfysize=18cm
\epsfbox{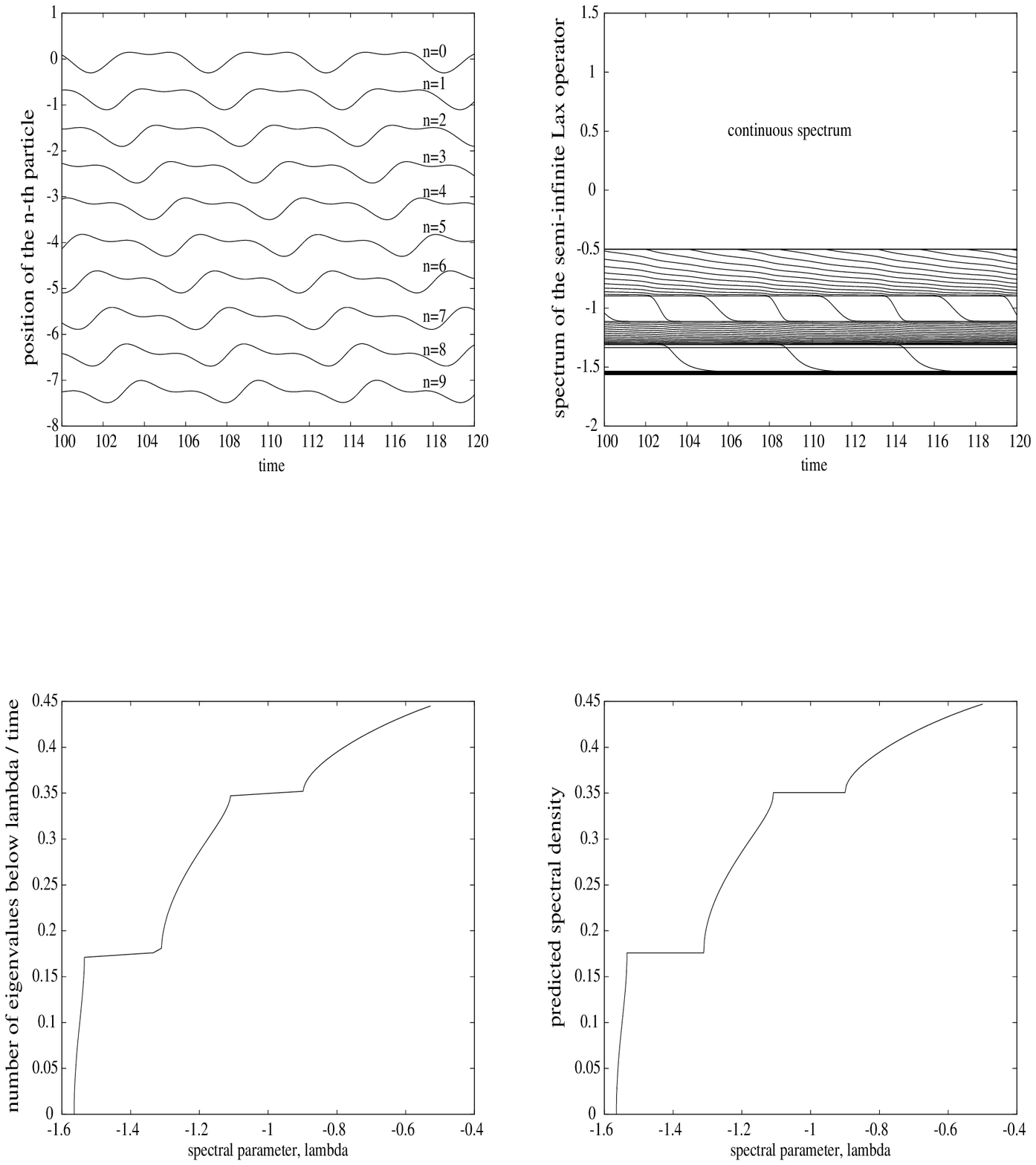}
\unitlength=1cm
\begin{picture}(13,0) 
\put(-10.5, 16.8){\makebox{Lattice motion}}
\put(-4.5, 8.0){\makebox{Predicted spectral density}}
\put(-12.75, 8.0){\makebox{Numerically observed spectral density}}
\put(-3.5, 16.8){\makebox{The spectrum}}
\end{picture}
\caption{{\em Toda lattice with driver 
$x_0(t) = t + 0.2 \sin \gamma t + 0.1 \cos \gamma t, \gamma = 1.1$, i.e.
$\gamma _2 > \gamma > \gamma _3$; see Section C.2 for detailed 
description.}}
\label{fc8}
\end{figure}

\end{document}


\newcommand{\D}{\displaystyle}
\newcommand{\SC}{\scriptstyle}
\newcommand{\T}{\textstyle}
\newcommand{\be}{\lfd \begin{equation}}
\newcommand{\ee}{\end{equation}}
\newcommand{\bea}{\begin{eqnarray*}}
\newcommand{\eea}{\end{eqnarray*}}
\newcommand{\beqa}{\begin{eqnarray}}
\newcommand{\eeqa}{\end{eqnarray}}
\newcommand{\mref}[1]{(\ref{#1})}
\newcommand{\bd}{\begin{displaymath}}
\newcommand{\ed}{\end{displaymath}}
\newcommand{\bef}{\lfd \begin{figure}[htb]}
\newcommand{\eef}{\end{figure}}

\newcommand{\summa}{\sum_{k=1}^{\infty}}
\newcommand{\summb}{\sum_{k=2}^{\infty}}
\newcommand{\rhof}{\rho_{F,c}}

\newcommand{\raua}{{\cal L} _{0,w}}
\newcommand{\raub}{{\cal L} _{\sigma ,w}}
\newcommand{\raubo}{{\cal L} _{\sigma _0,w}}
\newcommand{\rauc}{\ell _{1,w}}
\newcommand{\nora}{\| _{0,w}}
\newcommand{\norb}{\| _{\sigma ,w}}
\newcommand{\norbo}{\| _{\sigma _0,w}}
\newcommand{\norc}{\| _{\ell _{1,w}}}

\newcommand{\tue}{T _{u,\epsilon}}
\newcommand{\bue}{B _{u,\epsilon}}
\newcommand{\vue}{\tilde{v}}

\newcommand{\ganz}{{\db Z}}
\newcommand{\ganzz}{{\scriptsize {\db Z}}}
\newcommand{\natz}{{\db N}}
\newcommand{\relz}{{\db R}}
\newcommand{\comz}{{\db C}}

\newcounter{pictyp}

\newenvironment{proof}{\medskip {\bf Proof : }}{\begin{flushright}
                       $\Diamond$ \end{flushright} \medskip}

{\Huge \bf Guide to the Video}
\newline
\newline

Part I and II of the video presentation follow rather closely
Chapter 1 /resp. Chapter 2 of the text. In the third part, however,
the viewer might find it difficult to locate the details of 
what is discussed in the main text. Therefore we have provided
this guide. It gives precise references to the main text (bold face)
and repeats the formulae
which were written on the black board for simple 
identification (underlined).
\newline
\newline
{\Large \bf Part III}
\newline \par
We recall the long time behavior of the driven lattice and its 
dependence on the driving frequency $\gamma$, which was described in
detail in Part I. The goal of this presentation is to explain our result
on the construction of $\frac{2 \pi}{\gamma}$ -- time periodic solutions
of the semi-infinite lattice, ({\bf cf Chapters 3-5}), 
which are in 
good agreement with the numerically observed time asymptotic motion of
the lattice. More precisely we construct solutions of
\be
\underline
{\ddot{x}_n = F(x_{n-1} - x_n) - F(x_n - x_{n+1}),\quad n \geq 1},
\label{gu3.5}
\ee
with
\be
\underline
{x_0(t) := \epsilon \sum_{m \in \ganz} b_m e^{i \gamma m t}} 
\label{gu3.10}
\ee
({\bf cf (3.1), (3.2)}). For \underline{$ \epsilon = 0$}
the lattice at rest \underline{$x_n(t) = c n$} solves (\ref{gu3.5}),
(\ref{gu3.10}). We ask whether these solutions persist for 
\underline{$\epsilon \neq 0$}? Our results can be summarized as follows.

Let $\gamma _k := \frac{2}{k} \sqrt{F'(-c)}, k \in \natz$, be the sequence 
of threshold values for the driving frequency.
\begin{itemize}
\item
In the case of arbitrary force functions $F$ ({\bf cf Section 3.2.2, General
Assumptions}) and for \underline{$ \gamma > \gamma _1$} or
\underline{$ \gamma _1 > \gamma > \gamma _2$} there exist $\epsilon _0 > 0$
such that for all $0 \leq \epsilon < \epsilon _0$ we can construct
$\frac{2 \pi}{\gamma}$ -- time periodic solutions of (\ref{gu3.5}),
(\ref{gu3.10}), satisfying $x_n (t) - c n = O (\epsilon)$
({\bf cf Theorems 3.56 and 4.24}).
\item
In the case of the Toda lattice \underline{($F(x) = e^x$)}, for all 
$\gamma \in \relz ^+ \setminus \{ \gamma _k: k \in \natz \}$ there 
exist $\epsilon _0 > 0$
such that for all $0 \leq \epsilon < \epsilon _0$ we can construct
$\frac{2 \pi}{\gamma}$ -- time periodic solutions of (\ref{gu3.5}),
(\ref{gu3.10}), satisfying $x_n (t) - c n = O (\epsilon)$
({\bf cf Theorem 5.13}).
\end{itemize}

It is instructive to consider the linear case \underline{$F(x) = x$}.
By separation of variables we find that
\be
x_n(t) = \epsilon \sum_{m \in \ganz} b_m z_m^n e^{i \gamma m t}
\label{gu3.15}
\ee
solves (\ref{gu3.5}), (\ref{gu3.10}) , if 
\be
\underline{
z_m^2 + \delta _m z_m +1 = 0, \mbox{ with } \delta _m = -2 + (\gamma m)^2.
}
\label{gu3.20}
\ee
We notice that a mode \underline{$ z_m^n e^{i \gamma m t}$ } corresponds 
to a travelling wave in the case  
\underline{$|\delta _m| \leq 2$}, (
\underline{$z_m = e^{i \beta _m}$} for $\beta _m 
\in \relz$). For \underline{$|\delta _m| > 2$} one can choose 
\underline{$|z_m| < 1$},
i.e. the corresponding mode is decaying exponentially in $n$.
({\bf cf (1.25)}). The solutions we will 
construct in {\bf Chapters 3-5 }
will be nonlinear versions of (\ref{gu3.15}).

Our first result in {\bf Chapter 3 (Lemma 3.31 and Theorem 3.38)}
will reduce the construction of the desired solutions of (\ref{gu3.5}),
(\ref{gu3.10}) to showing that there exists a sufficiently large parameter
family of travelling wave solutions of the doubly infinite lattice.
More precisely, in the case that the driving frequency satisfies
$\gamma _k > \gamma > \gamma _{k+1}$, we have to construct a $2 k$ --
parameter family of $k$-phase waves. This reduction is proved by separating
the equations for the ``travelling wave part'' and the 
``exponentially decaying part'' ({\bf cf Lemma 3.31}) and then by solving
for the exponentially decaying part ({\bf cf Theorem 3.38}) (see
also {\bf Section 3.1}).

In the case of the Toda lattice the multi-phase travelling waves 
can be constructed for arbitrary number of phases in terms of theta 
functions, using the well known $g$-gap solutions. Due to limited time
these solutions are not discussed any further in this presentation.
See {\bf Chapter 5}, {\bf Appendix A} and {\bf Appendix B}
for more detail.

We will spend the remaining time to show how one can obtain single --
phase waves in the case of arbitrary force functions $F$ and what the 
difficulties are in the construction of multi-phase waves (in
particular two-phase waves), which we have not yet overcome.
By a two-phase wave we denote a solution of the doubly
infinite lattice of the form 
\be
\underline{
x_n(t) = c n + \chi _2 (n \beta _1 + \gamma t, n \beta _2 + 2 \gamma t),
}
\label{gu3.25}
\ee
where $\chi _2$ is $2 \pi$-periodic in both arguments.
Functions of this form are $\frac{2 \pi}{\gamma}$ -- time periodic.
The spatial frequencies $\beta _1, \beta _2$ turn out to be approximately
given by the dispersion relation for the linearized lattice,
\be
\underline{
2 \cos \beta _k = 2 - \frac{(\gamma k)^2}{F'(-c)}.
}
\label{gu3.30}
\ee
{\bf Remarks:}
\begin{itemize}
\item[(1)]
This formula was given incorrectly in the video presentation.
\item[(2)]
As mentioned above, two-phase waves are considered in the case
$\gamma _2 > \gamma > \gamma _3$, which means that there exist real
solutions $\beta _k$ of (\ref{gu3.30}) only in the case $|k| \leq 2$.
\end{itemize}
We first dicuss the construction of single-phase waves ({\bf see
Chapter 4}). Expand 
\be
x_n(t) = c n + \chi _1 (n \beta + \gamma t)
\label{gu3.35}
\ee
in a Fourier series
\be
\underline{
\chi _1 (\theta) = \sum_{m \in \ganzz} r(m) e^{i m \theta}.
}
\label{gu3.40}
\ee
The resulting equation for the sequence $r$ can be written in the form
\be
\underline{
\Lambda (\beta) r + W(r) = 0,
}
\label{gu3.45}
\ee
where $\Lambda (\beta)$ is a linear operator depending on $\beta$
and $W(r)$ denotes the nonlinear term.
The operator $\Lambda (\beta)$ is diagonal and its entries are given by
\lfd
\beqa
\Lambda (\beta) (m,m) &=& 2 \cos \beta m - 2 + \frac{(\gamma m)^2}{F'(-c)}
\nonumber \\
&=& - 4 \sin ^2 \frac{\beta m}{2} + \frac{(\gamma m)^2}{F'(-c)}
\label{gu3.50}
\eeqa
(See {\bf Lemma 4.10} for a derivation of (\ref{gu3.45}),
(\ref{gu3.50})).
\newline
\newline
{\bf Remark:}

In the video presentation we have given a different formula for $\Lambda
(\beta)$ which can be obtained by dividing (\ref{gu3.45}) by
$- 4 \sin ^2 \frac{\beta m}{2}$, which then represents the equation for the 
sequence $\tilde{r}(m) := r(m)( e^{- i \beta m} -1)$. Furthermore we have 
assumed $F'(-c) = 1$, which can always be achieved by rescaling time.
Therefore we obtain
\be
\underline{
\Lambda (\beta) (m,m) = 1 - \frac{(\gamma m)^2}
{ 4 \sin ^2 \frac{\beta m}{2}}.
}
\label{gu3.55}
\ee

In the video we proceed to give a brief description of the Lyapunov --
Schmidt decomposition which is explained in {\bf the introduction of 
Chapter 4 and  after the proof of Remark 4.12}.
The final result is that we construct a two-parameter family of
single phase waves of the doubly infinite lattice.

We repeat this ansatz in the case of two-phase waves,
\be
x_n(t) = c n + \chi _2 (n \beta _1 + \gamma t, n \beta _2 + 2 \gamma t),
\label{gu3.60}
\ee
\be
\underline{
\chi _2 (\theta _1, \theta _2) = \sum_{m \in \ganzz ^2} r(m) e^{i <m, \theta>}.
}
\label{gu3.65}
\ee
The equation for 
$\tilde{r}(m) := r(m)( e^{- i <\beta, m>} -1)$ can be written as
\be
\underline{
\Lambda (\beta) r + W(r) = 0,
}
\label{gu3.70}
\ee
with
\be
\underline{
\Lambda (\beta) (m,m) = 1 - \frac{\gamma ^2 (m_1 + 2 m_2)^2}
{ 4 \sin ^2 \frac{<\beta, m>}{2}}.
}
\label{gu3.75}
\ee
Note that we have again set $F'(-c) = 1$ without loss of generality.
Proceeding analogously to the single-phase case, we define the projection
$P$, $(P r)(m) :=  {\bf 1}_{\{|m| \neq 1\}} r(m)$, ($|m| = |m_1| + |m_2|$),
$Q := I - P$. Equation (\ref{gu3.70}) is equivalent to the following 
two equations. Denote $\mu := Pr, \varphi := Qr$.
\lfd
\beqa
P \Lambda (\beta) \mu + P W(\mu + \varphi) &=& 0. \label{gu3.80} \\ 
\lfd
Q \Lambda (\beta) \varphi + Q W(\mu + \varphi) &=& 0. \label{gu3.85} 
\eeqa
The first step in the Lyapunov-Schmidt decomposition is to 
satisfy the infinite dimensional equation (\ref{gu3.80}) by choosing
$\mu \in $ Ran ($P$) as a function of $\beta$ and $\varphi \in$ Ran ($Q$)
for ($\beta, \varphi$) in a neighborhood of ($\beta ^{(0)}, 0$).
We choose $\beta ^{(0)}$ such that $Q \Lambda (\beta ^{(0)}) = 0$, i.e.
$\beta ^{(0)} = (\beta _1^{(0)}, \beta _2^{(0)})$, with
\lfd
\beqa
4 \sin ^2 \frac{\beta _1^{(0)}}{2}  &=& \gamma ^2, \label{gu3.90} \\
\lfd
4 \sin ^2 \frac{\beta _2^{(0)}}{2} &=& (2 \gamma)^2. \label{gu3.95} 
\eeqa
Note again that (\ref{gu3.90}) and (\ref{gu3.95}) have real solutions
$\beta _1^{(0)}, \beta _2^{(0)}$ because of the condition 
$\gamma _2 > \gamma > \gamma _3$. 

In order to solve (\ref{gu3.80}) by implicit function theory,
we have to investigate the invertibility of
$\Lambda (\beta) :$ Ran ($P$) $\rightarrow$ Ran ($P$). In contrast
to the construction of the single-phase waves it is here that we 
encounter a major obstacle. In fact, a small divisor problem 
occurs, as there exist sites $m \in \ganz$ for which $\Lambda (\beta)
(m,m)$ is
arbitrarily close to zero. We call such sites {\em singular sites}.

We will determine the location of 
the singular sites of $\Lambda (\beta ^{(0)})$ (cf Fig. \ref{fgu1} below). 
First of all one realizes that $\gamma _2 > \gamma > \gamma _3$
implies that $m \in \ganz ^2$ can only be a singular site if
$m_1 + 2 m_2 = \pm 1, \pm 2$. We find the following cases.
\begin{itemize}
\item
For $m \in \ganz $, satisfying $m ( -2 \beta _1^{(0)} + \beta _2^{(0)})
\approx 0$ (mod $2 \pi$):\newline
($\pm 1 -2m, m$), ($- 2 m, m \pm 1$) are singular sites. They form a cross.
\item
For $m \in \ganz $, satisfying $m ( -2 \beta _1^{(0)} + \beta _2^{(0)})
\approx 2 \beta _1^{(0)} $ (mod $2 \pi$):\newline
($- 1 -2m, m$) is a singular site.
\item
For $m \in \ganz $, satisfying $m ( -2 \beta _1^{(0)} + \beta _2^{(0)})
\approx - 2 \beta _1^{(0)} $ (mod $2 \pi$):\newline
($1 -2m, m$) is a singular site.
\item
For $m \in \ganz $, satisfying $m ( -2 \beta _1^{(0)} + \beta _2^{(0)})
\approx 2 \beta _2^{(0)} $ (mod $2 \pi$):\newline
($-2m, m - 1$) is a singular site.
\item
For $m \in \ganz $, satisfying $m ( -2 \beta _1^{(0)} + \beta _2^{(0)})
\approx - 2 \beta _2^{(0)} $ (mod $2 \pi$):\newline
($-2m, m + 1$) is a singular site.
\end{itemize}
\bef
\leavevmode \epsfysize=7cm
\epsfbox{fgu1.ps}
\caption{The location of the singular sites of $\Lambda(\beta ^{(0)})$}
\label{fgu1}
\eef

A guide of how to construct solutions of (\ref{gu3.80}) in such a 
situation is the work of W. Craig and C. E. Wayne on the existence
of small solutions of the nonlinear wave equation \cite{CW}, where a similar
small divisor problem occured. Their construction proceeds via a Newton 
iteration scheme. The crucial point is to control the small eigenvalues
of a truncation of 
the linearized operator. One recalls that in a Newton scheme it is necessary
to linearize around the last approximation $u_n$ in order to obtain 
the better approximation $u_{n+1}$ and therefore the linearized operator
is in general not in diagonal form with respect to the standard basis.
The main technical tool in \cite{CW} is a technique introduced by J. Fr\"ohlich and 
T. Spencer(\cite{FS}) for the construction of inverse operators 
which makes it possible to first examine one singular site at a time
and then to control the interaction between them.

If one tries to apply this procedure to our problem one recognizes 
that the singular sites are not as well seperated as in \cite{CW}. This causes
additional difficulties which we have not resolved yet. Recent work
of Surace (\cite{S1}, \cite{S2}) 
contains a version of the Fr\"ohlich-Spencer
technique which might be helpful in overcoming these difficulties.